\documentclass[12pt]{article}
\RequirePackage[OT1]{fontenc}
\RequirePackage{amsthm,amsmath}
\RequirePackage[authoryear]{natbib}

\usepackage{amsmath}
\usepackage{graphicx,psfrag,epsf}
\usepackage{enumerate}
\usepackage{natbib}
\usepackage{multirow}

\usepackage{mathrsfs}
\usepackage{amssymb,amsfonts}
\usepackage{bm}
\usepackage{dsfont}
\newtheorem{algorithm}{Algorithm}
\newtheorem{corollary}{Corollary}
\newtheorem{example}{Example}[section]
\newtheorem{lemma}{Lemma}[section]
\newtheorem{theorem}{Theorem}[section]
\usepackage{xcolor}
\usepackage{epstopdf}% To incorporate .eps illustrations using PDFLaTeX, etc.
\usepackage[caption=false]{subfig}% Support for small, `sub' figures and tables
\usepackage[latin1]{inputenc}
\usepackage{color}
\usepackage{xcolor}
\usepackage{lineno,hyperref}
\usepackage{multicol}
\usepackage{algorithm}% http://ctan.org/pkg/algorithms
\usepackage{algpseudocode}% http://ctan.org/pkg/algorithmicx
\usepackage{lipsum}% http://ctan.org/pkg/lipsum
\usepackage{appendix}
\usepackage{amsmath}
\usepackage{graphicx}
\usepackage{lineno}
\usepackage{array}
\usepackage{longtable}
\usepackage{natbib}
\usepackage{float}
\usepackage{comment}
\usepackage{url}
\usepackage{booktabs}
\usepackage{multirow}
\usepackage{natbib}% Citation support using natbib.sty

\begin{document}

\title{\large\bf From Poisson Observations to Fitted Negative Binomial Distribution}
\author{Yingying Yang$^1$,  Niloufar Dousti Mousavi$^2$, Zhou Yu$^3$, and Jie Yang$^1$\\
		University of Illinois at Chicago$^1$, University of Pittsburgh$^2$,\\ and Merck \& Co., Inc.$^3$}
	
\maketitle

\begin{abstract}
The negative binomial distribution has been widely used as a more flexible model than the Poisson distribution for count data. However, when the true data-generating process is Poisson, it is often challenging to distinguish it from a negative binomial distribution with extreme parameter values, and existing maximum likelihood estimation procedures for the negative binomial distribution may fail or produce unstable estimates. To address this issue, we develop a new algorithm for computing the maximum likelihood estimate of negative binomial parameters, which is more efficient and more accurate than existing methods. We further extend  negative binomial distributions with a new parameterization to cover Poisson distributions as a special class. We provide theoretical justifications showing that, when applied to a Poisson data, the estimated parameters of the extended negative binomial distribution can consistently recover the true Poisson distribution.
\end{abstract}


\section{Introduction}\label{sec:intro}

Count data arises naturally from many different scientific disciplines. When the events or subjects of interest occur or are positioned randomly and independently in space or time, the Poisson distribution is often a good model \citep{pielou1977, ross1985negative}. It has a probability mass function (PMF)
\begin{equation}\label{eq:PMF_Poisson}
f_\lambda(y) = \frac{\lambda^y}{y!}e^{-\lambda}, \>\>\> y=0, 1, 2, \ldots,
\end{equation}
where $\lambda>0$ is both the mean and variance of the distribution.  However, in practice the phenomenon of over-dispersion, that is, the variance is larger than the mean, is frequently observed \citep{bliss1953fitting, ross1985negative}. For those scenarios, a negative binomial distribution with two parameters is more appropriate.

Following Chapter~32 in \cite{forbes2011statistical},  a negative binomial (NB) distribution with $\boldsymbol\theta = (\nu, p)^T \in (0, \infty) \times (0,1)$ has a PMF
\begin{equation}\label{eq:NB_pmf}
f_{\boldsymbol\theta}(y) = \frac{\Gamma(\nu+y)}{\Gamma(\nu)\Gamma(y+1)} p^{\nu} (1-p)^y, \>\>\> y=0,1,2,\ldots
\end{equation}
with mean $\nu(1-p)/p$ and variance $\nu(1-p)/p^2$. When $\nu$ takes a positive-integer value, the PMF~\eqref{eq:NB_pmf} describes the distribution of the number of failures before the $\nu$th success in a sequence of independent Bernoulli trials, where each trial has a success probability of $p$. The number of failures here is commonly referred to as the Pascal variate \citep{forbes2011statistical}. 

Alternative parameterizations for NB distributions have been used in the literature as well, including $(\nu, \mu=\nu(1-p)/p)$ as in \cite{anscombe1949statistical} or \cite{bliss1953fitting}, $(\nu, P=(1-p)/p)$ as in \cite{johnson2005univariate}, and $(\nu, P=1-p)$ as in \cite{anscombe1950sampling} and \cite{aldirawi2022modeling}. For more detailed reviews on the history of negative binomial distributions, please see \cite{ross1985negative}, Section~5.1 in \cite{johnson2005univariate},  and references therein.

According to \eqref{eq:PMF_Poisson} and \eqref{eq:NB_pmf}, the Poisson distributions and NB distributions are two non-overlapped classes of distributions. A natural question is whether we can tell in practice that the data come from a Poisson or a negative binomial distribution. If the true distribution is a negative binomial distribution with some unknown $\nu>0$ and $p\in (0,1)$, one can eventually tell that the true distribution is not Poisson because the sample variance will be significantly larger than the sample mean as the sample size goes to infinity.

However, if the data come from a Poisson distribution with some unknown parameter $\lambda >0$, the situation is quite different. Even with a large sample size, it is intrinsically difficult to separate a Poisson law from a nearly equidispersed negative binomial distribution (see Section~\ref{sec:supp_plots} in the Supplementary Material). This phenomenon is mainly due to the mathematical fact that a Poisson distribution can be approximated well by an NB distribution. For example, as a direct conclusion of Theorem~5.13.1 in \cite{fisz1963probability}, $F_{{\rm NB}(\nu, p(\nu))}(x) \rightarrow F_{\lambda}(x)$ as $\nu \rightarrow \infty$ for each real number $x$, where $p(\nu) = \nu/(\nu+\lambda)$, $F_{{\rm NB}(\nu, p)}$ is the cumulative distribution function (CDF) of the NB distribution with parameters $\nu>0$ and $p \in (0,1)$, and $F_{\lambda}$ is the CDF of the Poisson distribution with parameter $\lambda>0$. 

When the sample mean and sample variance of a dataset are about the same, Poisson and NB models become nearly indistinguishable, and likelihood-based numerical tools behave poorly. In particular, if an NB model is fitted to a dataset that is generated by a Poisson model, existing algorithms for finding the maximum likelihood estimate (MLE) may fail to converge, yield unstable estimates, 
or end with a suboptimal solution (see Section~\ref{sec:numerical_studies} for examples). In the meantime, the procedure may cost more computational time.

To overcome these issues, in this paper we propose a new algorithm for finding the MLE of NB distributions, called an {\it adaptive profile maximization algorithm} (APMA, see Algorithm~\ref{algo:MLE_nu_p} in Section~\ref{sec:mle_NB}), which outperforms the state-of-the-art programs available for users (see  Section~\ref{sec:numerical_studies}).
In Section~\ref{sec:NB_given_Poisson}, we theoretically justify the asymptotic behaviors of the MLE for NB parameters given a Poisson random sample. In particular, we show that 
{\it (i)} under the constraint $0<\nu\le\nu_{\max}$~, with probability $1$, the MLE of $\nu$ is $\nu_{\max}$ for all large enough sample size $n$ (see Theorem~\ref{thm:nu_max}); and {\it (ii)} if the constraint is relaxed to $0<\nu\le\nu_{\max}(n)$ with $\lim_{n\to\infty} \nu_{\max}(n) = \infty$, then the MLE of $\nu$ goes to infinity with probability $1$ (see Theorem~\ref{thm:nu_max_n}). Inspired by these asymptotic properties, in Section~\ref{sec:other_parameterizations},
we recommend a new parameterization with $\mu\in [0, \infty)$ and $p\in (0,1]$ for an extended NB family, which 
include Poisson($\lambda$) as NB$(\mu=\lambda>0, p=1)$.
With the extended NB$(\mu, p)$ family, we show that if $Y_1, \ldots, Y_n$ are iid from Poisson$(\lambda)$, then $\hat\mu = \bar{Y}_n \to \lambda$ and $\hat{p} \to 1$ almost surely, as $n$ goes to infinity
(Theorem~\ref{thm:mu_p_poisson}). Therefore, our corresponding algorithm (see Algorithm~\ref{algo:MLE_mu_p} in Section~\ref{sec:new_parameterization}) can identify the true underlying Poisson distribution consistently. We conclude in Section~\ref{sec:conclusion}.

\section{Finding MLE for Negative Binomial Distribution}\label{sec:mle_NB}

In this section, we consider finding the MLE of $(\nu, p)$ given $Y_1, \ldots, Y_n$ iid $\sim$ NB$(\nu, p)$ under the constraint $\nu \in (0, \nu_{\rm max}]$ and $p\in (0,1]$, where $\nu_{\rm max}>0$ is pre-determined. In practice, when numerically maximizing the likelihood for negative binomial models, the estimate of $\nu$ may tend to be unbounded, especially when the true underlying distribution is close to a Poisson distribution. In this case, an extremely large $\nu$ often causes overflow or convergence failures. To mitigate this possible issue and to ensure a fair comparison with existing algorithms that also restrict $\nu$ to a finite range (e.g., the R package AZIAD in \cite{dousti2022r}), we set a finite upper bound $\nu_{\rm max}$ for $\nu$ in this section.

\subsection{Likelihood equation}\label{sec:MLE}

In this section, we summarize some key results in the literature on finding the MLE of negative binomial parameters taking the form of \eqref{eq:NB_pmf}.
Then
\[
\log f_{\boldsymbol\theta}(y) = \log\Gamma(\nu+y) - \log\Gamma(\nu) - \log\Gamma(y+1) + \nu\log p + y\log(1-p)\ .
\]
Taking partial derivatives, we have
\begin{eqnarray*}
\frac{\partial\log f_{\boldsymbol\theta}(y)}{\partial \nu} &=& \Psi(\nu+y)-\Psi(\nu)+\log p\ , \\
\frac{\partial\log f_{\boldsymbol\theta}(y)}{\partial p} &=& \frac{\nu}{p}-\frac{y}{1-p}\ ,
\end{eqnarray*}
where $\Psi(\cdot) = \Gamma'(\cdot) / \Gamma(\cdot)$ is known as the {\it digamma} function (see, e.g., Section~6.3 in \cite{abramowitz1964handbook}).

Suppose $Y_1, \ldots, Y_n$ are iid $\sim$ NB$(\nu, p)$. The log-likelihood $l(\boldsymbol{\theta}) = \sum_{i=1}^n \log f_{\boldsymbol{\theta}}(Y_i)$ and
\[
\frac{\partial l(\boldsymbol{\theta})}{\partial p} = \frac{n\nu}{p} - \frac{\sum_{i=1}^n Y_i}{1-p}
= \frac{n(\nu+\bar{Y}_n)}{p(1-p)} \left(\frac{\nu}{\nu+\bar{Y}_n}-p\right)
\left\{\begin{array}{ll}
>0 & \mbox{ if } p<\frac{\nu}{\nu+\bar{Y}_n}\\
=0 & \mbox{ if } p=\frac{\nu}{\nu+\bar{Y}_n}\\
<0 & \mbox{ if } p>\frac{\nu}{\nu+\bar{Y}_n}
\end{array}\right.\> ,
\]
where $\bar{Y}_n=n^{-1} \sum_{i=1}^n Y_i$~. That is, $\partial l(\boldsymbol{\theta})/\partial p = 0$, known as the {\it likelihood equation} of $p$, always implies $\hat{p}=\hat{\nu}/(\hat\nu+\bar{Y}_n)$ for the MLE $(\hat{\nu}, \hat{p})$.

\begin{example}\label{ex:p_equal_1} {\bf A degenerated NB distribution}\quad {\rm
The solution $\hat{p}=\hat{\nu}/(\hat\nu+\bar{Y}_n)$ for the likelihood equation of $p$ suggests that $\hat{p}=1$ if $\bar{Y}_n=0$, which is possible in practice. To extend the range of $\boldsymbol{\theta} = (\nu, p)^T$ from $(0, \infty)\times (0,1)$ to $(0, \infty)\times (0,1]$, we consider the limit of $\log f_{\boldsymbol{\theta}}(y)$ as $p\to 1$ given $\nu>0$. Since $\lim_{p\to 1-}\log f_{\boldsymbol{\theta}}(y) = 0$ for $y=0$; and $-\infty$ for $y>0$, which corresponds to $f_{\boldsymbol{\theta}}(y) = 1$ for $y=0$; and $0$ for $y>0$, then the limiting distribution satisfies $P(Y_i=0)=1$. In other words, if we regard this singleton distribution, which generates zero with probability $1$, as a degenerated NB distribution as $p\to 1$, then the parameter space for the NB distribution is $\boldsymbol{\Theta} = (0, \infty) \times (0,1]$. 
}\hfill{$\Box$}
\end{example}

\begin{example}\label{ex:all_zero} {\bf A trivial case of all zeros}\quad {\rm
A special case of the data occurs when $Y_1 = \cdots = Y_n = 0$, which is possible in practice especially with a moderate $n$. In this case, $\bar{Y}_n = 0$ as well. It can be verified that $l(\boldsymbol{\theta}) = n\nu \log p \leq 0$. To maximize $l(\boldsymbol{\theta})$ with $\nu\in (0, \infty)$ and $p\in (0, 1]$, we must have $\hat{p}=1$ and $\hat{\nu}$ can be any positive number. 
}\hfill{$\Box$}
\end{example}

As a conclusion, regardless $\bar{Y}_n=0$ or $\bar{Y}_n>0$, the MLE of $p$ given $\nu>0$ must take the form as follows
\begin{equation}\label{eq:p_as_nu}
\hat{p} =\hat{p}(\nu) = \frac{\nu}{\nu + \bar{Y}_n}
\end{equation}
(see, e.g., \cite{bliss1953fitting} or Section~5.8.2 in \cite{johnson2005univariate}).

To avoid the trivial case as in Example~\ref{ex:all_zero}, from now on, we assume $\bar{Y}_n > 0$, that is, $Y_i>0$ for at least one $i$. 

Since the best $p$ must be the function $\hat{p}(\nu)$ of $\nu$, maximizing $l(\boldsymbol{\theta}) = l(\nu, p)$ is equivalent to maximizing the profile log-likelihood function $h(\nu) = l(\nu, \hat{p}(\nu))$ for $\nu$ \citep{cox1989analysis, venzon1988method}, which has two different forms: 
\begin{eqnarray}
& & h(\nu) = l(\nu, \hat{p}(\nu))\nonumber\\
    &=& \sum_{i=1}^n \log \Gamma(\nu + Y_i) - \sum_{i=1}^n \log \Gamma(Y_i+1) - n\log \Gamma(\nu)\nonumber\\
    & & -\ n(\nu+\bar{Y}_n)\log(\nu + \bar{Y}_n) + n\nu \log \nu + n\bar{Y}_n \log \bar{Y}_n\label{eq:h(nu)_psi}\\
    &=& \sum_{y\in I} f_y\cdot \log \Gamma(\nu + y) - \sum_{y\in I} f_y\cdot \log \Gamma(y+1) - n\log \Gamma(\nu)\nonumber\\
    & & -\ n(\nu+\bar{Y}_n)\log(\nu + \bar{Y}_n) + n\nu \log \nu + n\bar{Y}_n \log \bar{Y}_n\ .\label{eq:h(nu)}
\end{eqnarray}
Here $I$ in \eqref{eq:h(nu)} is the range (or the collection of distinct values) of $\{Y_1, \ldots, Y_n\}$,  and $f_y = \#\{1\leq i\leq n \mid Y_i = y\}$ is the frequency of $y \in I$. According to \eqref{eq:p_as_nu}, $\partial l(\nu, \hat{p}(\nu))/\partial p = 0$. Then the profile score 
\begin{eqnarray*}
h'(\nu) &=& \frac{\partial l}{\partial \nu}(\nu, \hat{p}(\nu)) + \frac{\partial l}{\partial p}(\nu, \hat{p}(\nu)) \cdot \hat{p}'(\nu) = \frac{\partial l}{\partial \nu}(\nu, \hat{p}(\nu))\\
&=& \sum_{i=1}^n \left[\Psi(\nu + Y_i) - \Psi(\nu)\right] - n \log \left(1 + \frac{\bar{Y}_n}{\nu}\right)\\
&\stackrel{\triangle}{=}& n g(\nu)\ ,
\end{eqnarray*}
where (see also \cite{bliss1953fitting})
\begin{eqnarray}
g(\nu) 
&=& \frac{1}{n}\sum_{i=1}^n \left[\Psi(\nu + Y_i) - \Psi(\nu)\right] - \log \left( 1 + \frac{\bar{Y}_n}{\nu}\right) \label{eq:g(nu)_psi}\\
&=& \frac{1}{n} \sum_{y \in I\setminus \{0\}} f_y \left(\frac{1}{\nu} + \cdots + \frac{1}{\nu + y - 1}\right) - \log\left(1 + \frac{\bar{Y}_n}{\nu}\right) \label{eq:g(nu)}
\end{eqnarray}
by using the fact $\Psi(x + 1) - \Psi(x) = x^{-1}$ (see, e.g., Formula~6.3.5 in \cite{abramowitz1964handbook}). The same technique has also been used for calculating the Fisher information matrices related to negative binomial distributions \citep{yu2024trigamma}.

In the literature, $g(\nu)$, defined by \eqref{eq:g(nu)_psi} or \eqref{eq:g(nu)}, is known as the {\it efficient score function} of $\nu$ up to a constant $n$ \citep{bliss1953fitting}. Essentially, $g(\nu)=0$ is the likelihood equation of $\nu$ for negative binomial distributions. If $g(\nu)=0$ has a unique solution $\hat\nu \in (0, \infty)$, then $\hat\nu$ must be the MLE of $\nu$. 

In our notation, \cite{anscombe1950sampling} conjectured that $g(\nu)=0$ has a unique solution if and only if $S_n^2 > \bar{Y}_n$, where $S_n^2 = n^{-1}\sum_{i=1}^n (Y_i - \bar{Y}_n)^2$.  \cite{simonsen1976solution, simonsen1980correction} proved  \citeauthor{anscombe1950sampling}'s conjecture and obtained the following theorem. For more efforts made in solving this conjecture and relevant problems, please see \cite{wang1996estimation}, \cite{bandara2019computing}, and reference therein.

\begin{theorem}[\cite{simonsen1976solution, simonsen1980correction}]\label{thm:simonsen}
Denote $M=\max\{Y_1, \ldots, Y_n\}$.
\begin{itemize}
\item[(1)] If $M=1$ or $M\geq 2$ and $S_n^2 \leq  \bar{Y}_n$, then $g(\nu)=0$ has no solution.
\item[(2)] If $M\geq 2$ and $S_n^2 > \bar{Y}_n$, then $g(\nu)=0$ has a unique solution in $(0, \infty)$.
\end{itemize}    
\end{theorem}

Actually, if $M=1$, we always have $\bar{Y}_n>0$ and $S_n^2 = \bar{Y}_n (1-\bar{Y}_n) < \bar{Y}_n$~. The following theorem covers the case of $M=0$.

\begin{theorem}\label{thm:M=0}
Under the constraint $\nu\in (0, \nu_{\rm max}]$ and $p\in (0,1]$, the MLE of $p$ is $\hat{p}=1$ if and only if $M=0$. In this case, $\hat{\nu}$ can be any positive number in $(0, \nu_{\rm max}]$.    
\end{theorem}
The proof of Theorem~\ref{thm:M=0}, as well as other proofs, is relegated to Section~\ref{sec:proofs} of the Supplementary Material.

If the data come from a Poisson distribution, the probability of $S_n^2 > \bar{Y}_n$ is roughly between $0.4$ and $0.5$ (see Table~\ref{tab:probability_estimate}). In this case, the MLE of $\nu$ can only be obtained by numerical algorithms.

\begin{table}[hbt]
\caption{$P(S_n^2 > \bar{Y}_n)$ estimated by $1,000$ random samples of size $n$ from Poisson($\lambda$)}\label{tab:probability_estimate}
\begin{center}
\begin{tabular}{@{}lrrrrr@{}}
\toprule
& \multicolumn{5}{@{}c}{\textbf{n}}  \\\cmidrule{2-6}
$\bm{\lambda}$&\textbf{50}&\textbf{500}&\textbf{5,000}&\textbf{50,000}&\textbf{500,000}\\
\midrule
\textbf{1}& 0.408 & 0.500 & 0.497 & 0.476 & 0.484\\
\textbf{3}& 0.441 & 0.445 & 0.461 & 0.503 & 0.487\\
\textbf{5}& 0.406 & 0.495 & 0.490 & 0.476 & 0.501\\
\textbf{10}& 0.408 & 0.487 & 0.496 & 0.496 & 0.469\\
\bottomrule
\end{tabular}
\end{center}
\end{table}

\subsection{Efficient score function} \label{sec:g_nu}

According to Theorem~\ref{thm:simonsen}, the efficient score function $g(\nu)$ does not have a solution in $(0, \infty)$ if $\bar{Y}_n>0$ and $S_n^2 \leq  \bar{Y}_n$~. In this section, we reveal more properties of $g(\nu)$ and also show that the MLE $\hat\nu = \nu_{\rm max}$ under the constraint $\nu \in (0, \nu_{\rm max}]$, if $g(\nu)=0$ does not have a solution in $(0, \infty)$.

Due to \eqref{eq:g(nu)_psi} or \eqref{eq:g(nu)}, the first-order derivative of $g(\nu)$ is 
\begin{eqnarray}
g'(\nu) &=& \frac{1}{n} \sum_{i=1}^n \Psi_1 (\nu+Y_i) - \Psi_1(\nu) + \frac{1}{\nu} - \frac{1}{\nu + \bar{Y}_n} \nonumber \\
&=& - \frac{1}{n} \sum_{y \in I\setminus \{0\}} f_y \left[\frac{1}{\nu^2} + \cdots + \frac{1}{(\nu + y - 1)^2}\right] + \frac{\bar{Y}_n}{\nu (\nu+\bar{Y}_n)} \ , \label{eq:g'(nu)} 
\end{eqnarray}
where $\Psi_1(\cdot) = \Psi'(\cdot)$ is known as the {\it trigamma} function (see, e.g., Section~6.4 in \cite{abramowitz1964handbook}).

\begin{lemma}\label{lem:g(nu)}
Suppose $\bar{Y}_n>0$. Then $\lim_{\nu\rightarrow 0^+} g(\nu) = \infty$, $\lim_{\nu\rightarrow\infty} g(\nu) = 0$, $\lim_{\nu\rightarrow 0^+} g'(\nu) = -\infty$, and  $\lim_{\nu\rightarrow\infty} g'(\nu) = 0$. Furthermore, $\lim_{\nu\rightarrow 0^+} \nu g(\nu) = 1-f_0/n > 0$, $\lim_{\nu\rightarrow\infty} {\nu}^2 g(\nu) = \frac{1}{2} \left(\bar{Y}_n - S^2_n\right)$, $\lim_{\nu\rightarrow 0^+} \nu^2 g'(\nu) = f_0/n - 1 < 0$, and $\lim_{\nu\rightarrow\infty} {\nu}^3 g'(\nu) = S^2_n - \bar{Y}_n$, where $S^2_n = n^{-1}\sum_{i=1}^n (Y_i - \bar{Y}_n)^2$.
\end{lemma}

\begin{example}\label{ex:all_c} {\bf A trivial case of constant observations}\quad {\rm
Suppose $Y_1 = \cdots = Y_n = c \geq 1$. 
Then $\bar{Y}_n=c$, $I=\{c\}$, $f_c=n$, and thus 
\begin{eqnarray*}
g'(\nu) &=& -\frac{1}{\nu^2} - \cdots - \frac{1}{(\nu+c-1)^2} + \frac{c}{\nu(\nu+c)}\\
& < & -\frac{1}{\nu(\nu+1)} - \cdots - \frac{1}{(\nu+c-1)(\nu+c)} + \frac{c}{\nu(\nu+c)}\\
&=& -\frac{1}{\nu} + \frac{1}{\nu+1} -\cdots -\frac{1}{v+c-1} + \frac{1}{\nu+c} + \frac{1}{\nu} - \frac{1}{\nu+c}\\
&=& 0
\end{eqnarray*}
for all $\nu > 0$. According to Lemma~\ref{lem:g(nu)}, $\lim_{\nu\rightarrow 0^+} g(\nu) = \infty$ and $\lim_{\nu\rightarrow\infty} g(\nu) = 0$, we must have $g(\nu) > 0$ for all $\nu > 0$ since $g(\nu)$ is strictly decreasing for all $\nu>0$ in this case.
According to the definition of $g(\nu)$, it implies that $\partial l(\boldsymbol{\theta})/\partial \nu = ng(\nu) > 0$ for all $\nu>0$. Then the MLE of $\nu$ must be the maximum value that $\nu$ is allowed to take.
That is, the MLE of $\nu$ under the constraint $\nu \in (0, \nu_{\rm max}]$ is $\hat\nu = \nu_{\rm max}$~.
}\hfill{$\Box$}
\end{example}

Example~\ref{ex:all_c} provides a case when $g(\nu) = 0$ does not have a solution. In practice, to find the MLE of $\nu$, instead of solving $g(\nu)=0$, we may minimize $g^2(\nu)$, which always has a solution.

\begin{theorem}\label{thm:min_g^2_positive}
Suppose $\bar{Y}_n>0$. If $\min_{\nu \in (0, \nu_{\rm max}]} g^2(\nu) > 0$, then the MLE of $\nu$ under the constraint $\nu \in (0, \nu_{\rm max}]$ is $\hat{\nu} = \nu_{\rm max}$~.
\end{theorem}

\subsection{A new algorithm for negative binomial distributions}\label{sec:algorithm}

In this section, we propose a new algorithm to find the MLE of NB distributions more accurately and more efficiently, by incorporating the previous theoretical results. 

Instead of the moment estimate $\bar{Y}_n^2/(S^2 - \bar{Y}_n)$ widely used for initializing $\nu$ \citep{fisher1941negative, bliss1953fitting, johnson2005univariate}, where $S^2 = \frac{1}{n-1} \sum_{i=1}^n (Y_i - \bar{Y}_n)^2$ is the sample variance, we propose the following estimate of $\nu$ as an initial value 
\begin{equation}\label{eq:moment_estimate_nu}
\hat\nu^{(0)} = \min\left\{ \nu_{\rm max}, \> \frac{\bar{Y}_n^2}{\max\{\varepsilon, \> S^2 - \bar{Y}_n\}}\right\}\ ,
\end{equation}
where $\nu_{\rm max} > 0$ is a predetermined upper bound for $\nu$ (e.g., $\nu_{\rm max} = 10^4$), and $\varepsilon > 0$ is a predetermined small positive value (e.g., $\varepsilon = 10^{-3}$). Note that $\hat\nu^{(0)} \in [0, \nu_{\rm max}]$ even if $S^2 \leq \bar{Y}_n$~, which is possible in practice. If we further have $\bar{Y}_n>0$, then $\hat\nu^{(0)} \in (0, \nu_{\rm max}]$.

We propose the following  algorithm, called an {\it adaptive profile maximization algorithm} (APMA, see Algorithm~\ref{algo:MLE_nu_p}), for finding the MLE of NB parameters. We also provide its R code in Section~\ref{sec:r_code} of the Supplementary Material.

\begin{algorithm}
\caption{Adaptive Profile Maximization Algorithm (APMA)
%: Finding MLE for Negative Binomial Distribution
}
\label{algo:MLE_nu_p}

Input: Data $Y_1, \ldots, Y_n \in \{0, 1, 2, \ldots\}$, predetermined threshold values $\nu_{\rm max}>0$, $\varepsilon>0$, and $\delta \in (0,1)$ (e.g., $\nu_{\rm max} = 10^4$, $\varepsilon = 10^{-3}$, $\delta=0.1$)

Output: MLE $(\hat\nu, \hat{p})$ for NB($\nu, p$)

\begin{algorithmic}[1]
\State Calculate the sample mean $\bar{Y}_n = \frac{1}{n} \sum_{i=1}^n Y_i$ and the sample variance $S^2 = \frac{1}{n-1} \sum_{i=1}^n (Y_i - \bar{Y}_n)^2$.
\State If $\bar{Y}_n=0$, go to Step~$7$ with $\hat\nu=1$ and $\hat{p}=1$; otherwise calculate the initial value $\hat\nu^{(0)}$ of $\nu$ according to \eqref{eq:moment_estimate_nu}.
\State If the ratio of distinct values in the sample $|I|/n < \delta$ (a few-distinct-values scenario), use \eqref{eq:h(nu)} and \eqref{eq:g(nu)} for computing $h(\nu)$ and $g(\nu)$, respectively; otherwise (a many-distinct-values scenario), use \eqref{eq:h(nu)_psi} and \eqref{eq:g(nu)_psi} instead. 
\State Maximize $h(\nu)$ over $(\varepsilon, \nu_{\max}]$ starting at $\hat{\nu}^{(0)}$ using the L-BFGS-B algorithm \citep{byrd1995limited}, a quasi-Newton method 
allowing box constraints. The analytic gradient is $h'(\nu) = ng(\nu)$, and the optimization is constrained within $(\varepsilon, \nu_{\max}]$. The procedure terminates when either the maximum number of iterations (e.g., 500) is reached or the convergence tolerance of the L-BFGS-B algorithm is satisfied. The optimization can be implemented, for example, using the \texttt{optim} function in R.
\State Let $\nu_* =$ optimizer output. If $h(\nu_{\rm max}) > h(\nu_*)$ (see \eqref{eq:h(nu)_psi} or \eqref{eq:h(nu)}), we let $\hat\nu=\nu_{\rm max}$; otherwise, $\hat\nu=\nu_*$~. Return a warning message when $\hat\nu = \nu_{\rm max}$~. 
\State Calculate $\hat{p} = \hat\nu/(\hat\nu + \bar{Y}_n)$.
\State Report $\hat\nu$ and $\hat{p}$.
\end{algorithmic}
\end{algorithm}

\section{Negative Binomial Estimates for Poisson Sample}\label{sec:NB_given_Poisson}

In this section, we derive two asymptotic results of the MLE of NB parameters given a Poisson random sample of size $n$: With probability $1$, {\it (i)} the MLE of $\nu$ under the constraint $0<\nu\leq \nu_{\rm max}$ is $\nu_{\rm max}$ for all large enough $n$ (see Theorem~\ref{thm:nu_max}); and {\it (ii)} if we adjust the constraint to be $0<\nu\leq \nu_{\rm max}(n)$ with $\lim_{n\rightarrow \infty}\nu_{\rm max}(n) = \infty$, then the MLE of $\nu$ approaches $\infty$, as $n$ goes to infinity (see Theorem~\ref{thm:nu_max_n}).

From Section~\ref{sec:g_nu}, we know that the efficient score function $g(\nu)$ (see  \eqref{eq:g(nu)_psi} or \eqref{eq:g(nu)}) plays an important role in finding the MLE of $\nu$ for an NB distribution.

In general, if $Y_1, \ldots, Y_n$ are iid from a probability distribution on $\{0, 1, 2, \ldots \}$, with CDF $F$ and mean $\mu = E(Y_1) < \infty$, we define  
\begin{eqnarray}
G_F(\nu) 
&=& E \Psi(\nu + Y_1) - \Psi(\nu) - \log \left( 1 + \frac{\mu}{\nu}\right) \label{eq:G_F(nu)_psi}\\
&=& \sum_{y=1}^\infty P(Y_1=y) \left(\frac{1}{\nu} + \cdots + \frac{1}{\nu + y - 1}\right) - \log\left(1 + \frac{\mu}{\nu}\right) \ ,\label{eq:G_F(nu)}
\end{eqnarray}
which is actually the limiting function of $g(\nu)$ as $n$ goes to infinity (see Lemma~\ref{lem:G_F(nu)_1-F} and Theorem~\ref{thm:G(nu)} below).

\begin{lemma}\label{lem:G_F(nu)_1-F}
Given $\mu = E(Y_1) < \infty$, $G_F(\nu)$ is a continuous function on $\nu \in (0, \infty)$. Furthermore,
\begin{equation}\label{eq:G_F(nu)_1-F}
    G_F(\nu) = \sum_{y=0}^\infty \frac{1-F(y)}{\nu+y}  - \log\left(1 + \frac{\mu}{\nu}\right)\ ,
\end{equation}  
where $F(y) = P(Y_1 \leq y)$ is the CDF of $Y_1$~.
\end{lemma}

\begin{theorem}\label{thm:G(nu)}
Suppose $Y_1, \ldots, Y_n$ are iid $\sim F$ on $\{0, 1, 2, \ldots \}$ with $\mu = E(Y_1) < \infty$. Then as $n\rightarrow \infty$, $g(\nu)\rightarrow G_F(\nu)$ almost surely.    
\end{theorem}

Using Kolmogorov's Strong Law of Large Numbers (see, e.g., Theorem 7.5.1 in \cite{resnick1999probpath}), we provide a rigorous proof of Theorem~\ref{thm:G(nu)} in Section~\ref{sec:proofs} of the Supplementary Material. It is actually applicable for discrete distributions other than Poisson distributions as well, as long as the distribution is on nonnegative integers. For example, in Lemma~\ref{lem:G_NB(nu)}, we show that $G_F(\nu)\equiv 0$ if $F$ is an NB CDF.

\begin{lemma}\label{lem:G_NB(nu)}
    Suppose $Y_1, \ldots, Y_n$ are iid $\sim$ NB$(\nu, p)$ with $\nu>0$ and $0<p<1$. Then $G_{NB(\nu, p)}(\nu)=0$, for all $\nu>0$.
\end{lemma}

Suppose $Y_1, \ldots, Y_n$ are iid $\sim$ Poisson$(\lambda)$ with $\lambda > 0$. For $y \in \{0, 1, 2, \ldots \}$, we denote $f_\lambda (y) = P(Y_1=y) = \frac{\lambda^y}{y!}e^{-\lambda}$, $F_\lambda (y) = P(Y_1 \leq y)$, and 
\begin{eqnarray}
G_\lambda(\nu) &=& \sum_{y=1}^\infty \frac{\lambda^y}{y!}e^{-\lambda} \left(\frac{1}{\nu} + \cdots + \frac{1}{\nu + y - 1}\right) - \log\left(1 + \frac{\lambda}{\nu}\right)\nonumber\\
&=& \sum_{y=0}^\infty \frac{1-F_{\lambda}(y)}{\nu+y}  - \log\left(1 + \frac{\lambda}{\nu}\right) \label{eq:G_lambda_nu}
\end{eqnarray}
as the almost sure limit of $g(\nu)$ given a random sample of Poisson$(\lambda)$.  

Given $\lambda > 0$ and $\nu>0$, we let $p = \nu/(\nu+\lambda) \in (0,1)$. Then the mean of NB$(\nu, p)$ is $\nu(1-p)/p = \nu \lambda/\nu = \lambda$, the same as Poisson$(\lambda)$'s. For such an NB$(\nu, p)$, we denote $f_{NB(\nu, p)}(y)$ and $F_{NB(\nu, p)}(y)$ as the corresponding PMF and CDF, respectively. According to Lemma~\ref{lem:G_NB(nu)}, 
\begin{eqnarray*}
G_{NB(\nu, p)}(\nu) &=& \sum_{y=0}^\infty \frac{1-F_{NB(\nu, p)}(y)}{\nu+y}  - \log\left(1 + \frac{\mu}{\nu}\right)\\
&=& \sum_{y=0}^\infty \frac{1-F_{NB(\nu, p)}(y)}{\nu+y}  - \log\left(1 + \frac{\lambda}{\nu}\right)\\
&=&0\ .
\end{eqnarray*}
Comparing $G_{NB(\nu, p)}(\nu)$ with \eqref{eq:G_lambda_nu}, if we can show that 
\begin{equation}\label{eq:F_lambda_F_NB}
\sum_{y=0}^\infty \frac{1-F_{\lambda}(y)}{\nu+y} > \sum_{y=0}^\infty \frac{1-F_{NB(\nu, p)}(y)}{\nu+y}\ ,
\end{equation}
then $G_{\lambda}(\nu) > G_{NB(\nu, p)}(\nu) =0$, for each $\nu>0$.

To prove \eqref{eq:F_lambda_F_NB}, we first denote the continuous version of $f_\lambda(y)$ and $f_{NB(\nu, p)}(y)$ by
\begin{eqnarray}
    f_\lambda(x) &=& \frac{\lambda^x}{\Gamma(x+1)}e^{-\lambda}\label{eq:f_lambda_x}\ ,\\
    f_{NB(\nu, p)}(x) &=& \frac{\Gamma(\nu+x)}{\Gamma(x+1)\Gamma(\nu)} p^\nu (1-p)^x = \frac{\Gamma(\nu+x)}{\Gamma(x+1)\Gamma(\nu)} \cdot \frac{\nu^\nu \lambda^x}{(\nu+\lambda)^{\nu+x}}\ ,\label{eq:f_NB_x}
\end{eqnarray}
respectively, where $x\in [0, \infty)$ and $p=\nu/(\nu+\lambda)$. Note that \eqref{eq:f_lambda_x} and \eqref{eq:f_NB_x} provide the PMFs of Poisson$(\lambda)$ and NB$(\nu, p)$, respectively, if we take $x$ in $\{0, 1, 2, \ldots \}$ only.

\begin{lemma}\label{lem:compare_f_p_f_NB_log}
We let $r(x) = \log f_{NB(\nu, p)}(x) - \log f_\lambda(x)$, where $x\in [0, \infty)$ and $p=\nu/(\nu+\lambda)$. Then its first derivative
\[
r'(x) \left\{
\begin{array}{cl}
<0 & \mbox{ if } 0 < x < x_*\\
=0 & \mbox{ if } x = x_*\\
>0 & \mbox{ if } x > x_*
\end{array}\right.\> ,
\]
where $x_* = \Psi^{-1}\left(\log(\nu+\lambda)\right) - \nu > \lambda$. Furthermore, $r(0) > 0$ and $\lim_{x\rightarrow \infty} r(x) = \infty$.
\end{lemma}

Lemma~\ref{lem:compare_f_p_f_NB_log} tells us that $f_{NB(\nu, p)}(x) < \log f_\lambda(x)$, for $0 < x < x_*$~, and $f_{NB(\nu, p)}(x) > \log f_\lambda(x)$, for $x > x_*$~. As the next step, we compare $F_{NB(\nu, p)}(y)$ and $F_\lambda(y)$ in Theorem~\ref{thm:compare_f_p_f_NB}.

\begin{theorem}\label{thm:compare_f_p_f_NB}
Given $\lambda > 0$ and $\nu>0$, we let $p=\nu/(\nu+\lambda) \in (0,1)$, $d(y) = f_{NB(\nu, p)}(y) - f_\lambda(y)$, and $D(y) = \sum_{k=0}^y d(k) = F_{NB(\nu, p)}(y) - F_\lambda(y)$, $y\in \{0, 1, 2, \ldots \}$. Then there exist three positive integers $K_1$~, $K_*$~, and $K_2$~, such that, $0 < K_1 < K_* < K_2 < \infty$, 
\[
d(y) \left\{
\begin{array}{cl}
>0 & \mbox{ if } 0 \leq y < K_1 \mbox{ or } y > K_2\\
\geq 0 & \mbox{ if } y = K_1 \mbox{ or } K_2\\
<0 & \mbox{ if } K_1 < y < K_2
\end{array}\right.
,\>\>\>
D(y) \left\{
\begin{array}{cl}
>0 & \mbox{ if } 0 \leq y < K_*\\
\geq 0 & \mbox{ if } y = K_*\\
<0 & \mbox{ if } y > K_*
\end{array}\right.\> .
\]
Furthermore, $D(y)$ strictly increases before $K_1$ and attains its maximum  at $y=K_1$ with $D(K_1)>0$, then strictly decreases between $K_1$ and $K_2$ and attains its minimum at $y=K_2$ with $D(K_2) < 0$, and then strictly increases after $K_2$ and approaches to $0$ as $y$ goes to $\infty$.
\end{theorem}

\begin{lemma}\label{lem:sum_D(y)}
Given $D(y)$ as defined in Theorem~\ref{thm:compare_f_p_f_NB}, we must have $\sum_{y=0}^\infty D(y) = 0$.
\end{lemma}

Now we are ready to show \eqref{eq:F_lambda_F_NB} and thus $G_\lambda(\nu) > 0$ for all $\nu>0$ in Theorem~\ref{thm:G_lambda_nu>0}.

\begin{theorem}\label{thm:G_lambda_nu>0}
    For any $\lambda > 0$ and $\nu>0$, we must have $G_\lambda (\nu) > 0$.
\end{theorem}

According to Theorem~\ref{thm:G(nu)}, given $Y_1, \ldots, Y_n$ iid $\sim$  Poisson($\lambda$),  we have $g(\nu) \rightarrow G_\lambda(\nu)$ almost surely, as $n$ goes to infinity. Now Theorem~\ref{thm:G_lambda_nu>0} guarantees that $G_\lambda(\nu)>0$ for all $\nu>0$. It is straightforward to derive $g(\nu)>0$ for large enough $n$ as in Lemma~\ref{lem:g_lambda_nu>0}.  

\begin{lemma}\label{lem:g_lambda_nu>0}
Suppose $Y_1, \ldots, Y_n$ are iid $\sim$ Poisson($\lambda$) with $\lambda >0$. Then for each $\nu>0$, there exists a constant $c>0$, such that, with probability $1$, $g(\nu) > c$ for all large enough $n$.
\end{lemma}

For each $\nu>0$, Lemma~\ref{lem:g_lambda_nu>0} guarantees $g(\nu)>c$ for some $c>0$ and large enough $n$. Nevertheless, it does not guarantee that $g(\nu)>0$ for all $\nu>0$ and large enough $n$. To show that $g(\nu)>0$ for all $\nu \in (0, \nu_{\rm max}]$ and large enough $n$, we split $(0, \nu_{\rm max}]$ into $(0, \nu_0]$ and $(\nu_0, \nu_{\rm max}]$, where $\nu_0$ is some positive number that is close enough to $0$. Then we bound $g(\nu)$ from below for large enough $n$. 

First, we find a lower bound for $g(\nu)$ on $\nu\in (0, \nu_0]$. Actually, according to \eqref{eq:g(nu)}, if $Y_1, \ldots, Y_n$ are iid $\sim$ Poisson($\lambda$) with $\lambda > 0$, then
\begin{eqnarray*}
& & g(\nu)\\ 
&=& \frac{1}{n} \sum_{y \in I\setminus \{0\}} f_y \left(\frac{1}{\nu} + \cdots + \frac{1}{\nu + y - 1}\right) - \log\left(1 + \frac{\bar{Y}_n}{\nu}\right)\\
&=& \frac{1}{n} \cdot \frac{\sum_{y=1}^\infty f_y}{\nu}  + \frac{1}{n} \sum_{y \in I\setminus \{0,1\}} f_y \left(\frac{1}{\nu+1} + \cdots + \frac{1}{\nu + y - 1}\right) - \log\left(1 + \frac{\bar{Y}_n}{\nu}\right)\\
&\geq & \frac{1}{\nu} \cdot \frac{\#\{i:Y_i>0\}}{n} - \log\left(1 + \frac{\bar{Y}_n}{\nu}\right)\\
&\rightarrow & \frac{1-e^{-\lambda}}{\nu} - \log\left(1 + \frac{\lambda}{\nu}\right) 
\end{eqnarray*}
almost surely, as $n$ goes to $\infty$.

\begin{lemma}\label{lem:g_nu_0_nu0}
Given $\lambda >0$ and $c>0$, there exists a $\nu_0>0$, such that, 
\[
h(\nu, a, b) := \frac{a}{\nu} - \log\left(1+\frac{b}{\nu}\right) \geq c
\]
for all $0 < \nu \leq \nu_0$~, $a\in \left[\frac{1}{2}(1-e^{-\lambda}), \frac{3}{2}(1-e^{-\lambda})\right]$, and $b\in \left[\frac{\lambda}{2}, \frac{3\lambda}{2}\right]$.
\end{lemma}

Note that $\nu_0$ in Lemma~\ref{lem:g_nu_0_nu0} does not depend on $Y_i$'s. Since $\#\{i:Y_i>0\}/n \rightarrow 1-e^{-\lambda}$ and $\bar{Y}_n \rightarrow \lambda$ almost surely, as $n$ goes to $\infty$, there exists an event 
\begin{equation}\label{eq:A_set}
A = \{\omega \in \Omega \mid \lim_{n\rightarrow \infty} \frac{\#\{i:Y_i(\omega)>0\}}{n} = 1-e^{-\lambda} \mbox{ and } \lim_{n\rightarrow \infty} \bar{Y}_n(\omega) = \lambda\}
\end{equation} 
with probability $1$, where $\Omega$ is the sample space.
Given $\omega \in A$, for large enough $n$, $\#\{i:Y_i>0\}/n \in \left[\frac{1}{2}(1-e^{-\lambda}), \frac{3}{2}(1-e^{-\lambda})\right]$, $\bar{Y}_n\in \left[\frac{\lambda}{2}, \frac{3\lambda}{2}\right]$, and thus $g(\nu)(\omega) \geq c$ for all $\nu \in (0, \nu_0]$. That is, with probability $1$, we obtain a common lower bound $c$ of $g(\nu)$, for all $\nu\in (0, \nu_0]$, with large enough $n$.

Secondly, we construct a lower bound for $g(\nu)$, $\nu \in (\nu_0, \nu_{\rm max}]$, with large enough $n$, which relies on Lemma~\ref{lem: upperbound_gr_gnu}.

\begin{lemma}\label{lem: upperbound_gr_gnu}
For any given $\nu_0>0$, there exists an $M>0$, which does not depend on the data, such that $|g'(\nu)| \leq M$ for all $\nu>\nu_0$~.
\end{lemma}

With the aids of Theorem~\ref{thm:G_lambda_nu>0}, Lemmas~\ref{lem:g_lambda_nu>0}, \ref{lem:g_nu_0_nu0}, and \ref{lem: upperbound_gr_gnu}, we are ready to derive the first major theoretical result of this paper as described in Theorem~\ref{thm:nu_max}.

\begin{theorem}\label{thm:nu_max}
Suppose $Y_1, \ldots, Y_n$ are iid $\sim$ Poisson($\lambda$) with $\lambda >0$, then for any $\nu_{\rm max}>0$, with probability 1, the MLE  of $(\nu, p)$ for a negative binomial distribution under the constraint $0<\nu\leq \nu_{\rm max}$ is $\left(\nu_{\rm max}, \frac{\nu_{\rm max}}{\nu_{\rm max} + \bar{Y}_n}\right)$ for all large enough $n$.
\end{theorem}

In practice, we may allow $\nu_{\rm max}$ to increase along with $n$. As a direct conclusion of the proof of Theorem~\ref{thm:nu_max}, we have the following theorem.

\begin{theorem}\label{thm:nu_max_n}
Suppose $Y_1, \ldots, Y_n$ are iid $\sim$ Poisson($\lambda$) with $\lambda >0$. Let $\nu_{\rm max}(n)$ be a sequence of positive real numbers such that $\nu_{\rm max}(n) \rightarrow \infty$ as $n$ goes to $\infty$. Denote $(\hat{\nu}_n, \hat{p}_n)$ as the MLE of $(\nu, p)$ for a negative binomial distribution under the constraint $0<\nu\leq \nu_{\rm max}(n)$. Then with probability $1$, $\hat\nu_n \rightarrow \infty$, as $n$ goes to $\infty$.
\end{theorem}

With the aid of Theorem~\ref{thm:nu_max_n}, we can derive the limiting behaviors of the Kolmogorov-Smirnov test statistic \citep{massey1951kolmogorov} for NB distributions as well (see Section~\ref{sec:supp_thms} in the Supplementary Material), which is a classical goodness-of-fit procedure but out of the scope of this paper.

\section{Other Parameterizations for NB Distributions }\label{sec:other_parameterizations}

In this section, we first review other parameterizations in the literature for negative binomial distributions, and derive the corresponding asymptotic properties of their MLEs as direct conclusions of the invariance property of MLE \citep{zehna1966invariance, pal1992invariance}. We further propose a new parameterization of negative binomial distributions, which extends the NB distributions to cover Poisson distributions as a special class. Under the new parameterization, the MLEs from the extended NB distribution satisfy $\hat{\mu}=\bar{Y}_n$ and $\hat{p}\to 1$ for Poisson data, and thus can consistently recover the true Poisson distribution (see Theorem~\ref{thm:mu_p_poisson}).

\subsection{Asymptotic results with other parameterizations in the literature}\label{sec:other_parameterizations_literature}

\cite{anscombe1949statistical} and \cite{bliss1953fitting} adopted $\nu$ and $\mu = \nu(1-p)/p$ as the parameters for a negative binomial distribution with the PMF
\begin{equation}\label{eq:pmf_nu_mu}
f_{\nu, \mu}(y) =  \frac{\Gamma(\nu+y)}{\Gamma(\nu)\Gamma(y+1)} \left(\frac{\nu}{\nu+\mu}\right)^{\nu} \left(\frac{\mu}{\nu+\mu}\right)^y, \>\>\> y=0,1,2,\ldots\ .
\end{equation}

\begin{theorem}\label{thm:nu_mu}
Given nonnegative integers $Y_1, \ldots, Y_n$, let $(\hat{\nu}, \hat{\mu})$ denote an MLE of $(\nu, \mu)$ for the negative binomial distribution described as in \eqref{eq:pmf_nu_mu} with constraints $\nu\in (0, \nu_{\rm max}]$ and $\mu \in [0, \infty)$. Then (i) $\hat{\mu}=0$ if and only if $\bar{Y}_n = 0$; and (ii) if $Y_1, \ldots, Y_n$ are iid $\sim$ Poisson($\lambda$) with $\lambda >0$, then with probability 1, $\hat{\nu} = \nu_{\rm max}$ and $\hat{\mu} = \bar{Y}_n$ for all large enough $n$.
\end{theorem}

\cite{johnson2005univariate} used $\nu$ and $P = (1-p)/p$ for a negative binomial distribution with the PMF
\begin{equation}\label{eq:pmf_nu_P}
f_{\nu, P}(y) =  \frac{\Gamma(\nu+y)}{\Gamma(\nu)\Gamma(y+1)} \left(\frac{1}{1+P}\right)^{\nu} \left(\frac{P}{1+P}\right)^y, \>\>\> y=0,1,2,\ldots\ .
\end{equation}

\begin{theorem}\label{thm:nu_P}
Given nonnegative integers $Y_1, \ldots, Y_n$, let $(\hat{\nu}, \hat{P})$ denote an MLE of $(\nu, P)$ for the negative binomial distribution described as in \eqref{eq:pmf_nu_P} with constraints $\nu\in (0, \nu_{\rm max}]$ and $P \in [0, \infty)$. Then (i) $\hat{P}=0$ if and only if $\bar{Y}_n = 0$; and (ii) if $Y_1, \ldots, Y_n$ are iid $\sim$ Poisson($\lambda$) with $\lambda >0$, then with probability 1, $\hat{\nu} = \nu_{\rm max}$ and $\hat{P} = \bar{Y}_n/\nu_{\rm max}$ for all large enough $n$.
\end{theorem}

\cite{anscombe1950sampling} and \cite{aldirawi2022modeling} used $\nu$ and $P = 1-p$ for a negative binomial distribution with the PMF
\begin{equation}\label{eq:pmf_nu_1-p}
f_{\nu, P}(y) =  \frac{\Gamma(\nu+y)}{\Gamma(\nu)\Gamma(y+1)} (1-P)^{\nu} P^y, \>\>\> y=0,1,2,\ldots\ .
\end{equation}

\begin{theorem}\label{thm:nu_1-p}
Given nonnegative integers $Y_1, \ldots, Y_n$, let $(\hat{\nu}, \hat{P})$ denote an MLE of $(\nu, P)$ for the negative binomial distribution described as in \eqref{eq:pmf_nu_1-p} with constraints $\nu\in (0, \nu_{\rm max}]$ and $P \in [0, 1)$. Then (i) $\hat{P}=0$ if and only if $\bar{Y}_n = 0$; and (ii) if $Y_1, \ldots, Y_n$ are iid $\sim$ Poisson($\lambda$) with $\lambda >0$, then with probability 1, $\hat{\nu} = \nu_{\rm max}$ and $\hat{P} = \bar{Y}_n/(\nu_{\rm max}+\bar{Y}_n)$ for all large enough $n$.
\end{theorem}

So far all the four parameterizations for negative binomial distributions, namely \eqref{eq:NB_pmf}, \eqref{eq:pmf_nu_mu}, \eqref{eq:pmf_nu_P}, \eqref{eq:pmf_nu_1-p}, adopted $\nu$, whose original meaning is the given number of successes in a sequence of independent Bernoulli trials, while the negative binomial random variable is the observed number of failures before the $\nu$th success. Although we knew that a Poisson$(\lambda)$ distribution can be regarded as the limiting distribution as $\nu\to \infty$  and $p\to 1$ with $\nu(1-p)/p = \lambda$, the current parameterizations do not allow us to simply claim that $\hat{\nu}=\infty$, if the data come from a Poisson distribution.

\subsection{A new parameterization and extended NB distributions}\label{sec:new_parameterization}

In this section, we propose a new parameterization, namely $\mu \in [0, \infty)$ and $p\in (0,1]$ for an extended negative binomial distribution. It includes both the degenerated case (see Example~\ref{ex:p_equal_1}) with $\mu=0$ and Poisson distributions with $\mu>0$ and $p=1$. More specifically, {\it (i)} for $\mu > 0$ and $p\in (0,1)$, its PMF is defined as
\begin{equation}\label{eq:pmf_mu_p}
f_{\mu, p}(y) =  \frac{\Gamma\left(\frac{\mu p}{1-p}+y\right)}{\Gamma\left(\frac{\mu p}{1-p}\right)\Gamma(y+1)} p^{\frac{\mu p}{1-p}} (1-p)^y, \>\>\> y=0,1,2,\ldots\ ;
\end{equation}
{\it (ii)} for $p\in (0,1)$, as the limit when $\mu\to 0^+$, its PMF is
\begin{equation}\label{eq:pmf_mu_0_p}
f_{0,p}(y) = \lim_{\mu\to 0^+} f_{\mu, p}(y) = \left\{
\begin{array}{cl}
1 & \mbox{ if }y=0\\
0 & \mbox{ otherwise}
\end{array}\right.\ ,
\end{equation}
namely the degenerated case described in Example~\ref{ex:p_equal_1}; {\it (iii)} for $\mu>0$, as the limit when $p\to 1^-$, its PMF is 
\begin{equation}\label{eq:pmf_mu_p_1}
f_{\mu, 1}(y) =  \lim_{p\to 1^-} f_{\mu, p}(y) =\frac{\mu^y}{y!} e^{-\mu}, \>\>\> y=0,1,2,\ldots\ ,
\end{equation}
namely the Poisson$(\mu)$ distribution; and {\it (iv)} as the limit when $\mu\to 0^+$ and $p\to 1^-$, its PMF is
\begin{equation}\label{eq:pmf_mu_0_p_1}
f_{0,1}(y) = \lim_{\begin{smallmatrix}\mu\to 0^+\\ p\to 1^-\end{smallmatrix}} f_{\mu, p}(y) = \left\{
\begin{array}{cl}
1 & \mbox{ if }y=0\\
0 & \mbox{ otherwise}
\end{array}\right.\ ,
\end{equation}
namely the degenerated case in Example~\ref{ex:p_equal_1} again.

Given a nonnegative-integer-valued dataset $\{Y_1, \ldots, Y_n\}$, the log-likelihood function of $(\mu, p) \in [0, \infty)\times (0,1]$ is $l(\mu, p) = \sum_{i=1}^n \log f_{\mu, p}(Y_i)$. Recall that the function $g(\nu)$ defined in \eqref{eq:g(nu)_psi} or \eqref{eq:g(nu)} plays a central role in finding the MLE of $\nu$. It can be verified that for $\mu>0$ and $p\in (0,1)$, 
\begin{eqnarray}
\frac{\partial l(\mu, p)}{\partial \mu} &=& \frac{np}{1-p} \left\{g\left(\frac{\mu p}{1-p}\right) + \log\left[p + \frac{\bar{Y}_n}{\mu}(1-p)\right]\right\}\ ,\label{eq:l_mu(mu,p)}\\
\frac{\partial l(\mu, p)}{\partial p} &=& \frac{n\mu}{(1-p)^2} \left\{g\left(\frac{\mu p}{1-p}\right) + \log\left[p + \frac{\bar{Y}_n}{\mu}(1-p)\right]\right\} + \frac{n}{1-p}(\mu - \bar{Y}_n)\ .\label{eq:l_p(mu,p)}
\end{eqnarray}

According to the definition of $f_{\mu, p}(y)$ under different scenarios, we have the following lemma and theorem for an MLE $(\hat{\mu}, \hat{p}) \in [0, \infty)\times (0,1]$.

\begin{lemma}\label{lem:continuity_mu_p}
For each $y\in \{0,1,2,\ldots\}$, $f_{\mu, p}(y)$ as defined by equations \eqref{eq:pmf_mu_p}, \eqref{eq:pmf_mu_0_p}, \eqref{eq:pmf_mu_p_1}, and \eqref{eq:pmf_mu_0_p_1} is a continuous function of $\mu\in [0, \infty)$ and $p\in (0, 1]$, so does the log-likelihood $l(\mu, p)$.    
\end{lemma}

\begin{theorem}\label{thm:mu_p_MLE_existence}
There must exist a $\hat{\mu}\in [0, \infty)$ and a $\hat{p}\in (0,1]$ that maximizes $l(\mu, p)$ under the constraint $(\mu, p) \in [0, \infty)\times (0,1]$. If $\bar{Y}_n=0$, then we must have $\hat{\mu}=0$ while $\hat{p}$ can be any number in $(0,1]$; if $\bar{Y}_n > 0$, then we must have $\hat{\mu}=\bar{Y}_n$~.   
\end{theorem}

\begin{theorem}\label{thm:mu_p_poisson}
Given nonnegative integers $Y_1, \ldots, Y_n$, let $(\hat{\mu}, \hat{p})$ denote an MLE of $(\mu, p)$ for the extended negative binomial distribution described via \eqref{eq:pmf_mu_p}, \eqref{eq:pmf_mu_0_p}, \eqref{eq:pmf_mu_p_1}, and \eqref{eq:pmf_mu_0_p_1} with parameters $\mu\in [0, \infty)$ and $p\in (0, 1]$. Then (i) $\hat{\mu}=0$ if and only if $\bar{Y}_n = 0$; (ii) we always have $\hat{\mu} = \bar{Y}_n$; and (iii) if $Y_1, \ldots, Y_n$ are iid $\sim \text{Poisson}(\lambda)$ with $\lambda >0$, then $\lim_{n\to\infty}\hat{p} = 1$ with probability $1$.  
\end{theorem}

According to Theorem~\ref{thm:mu_p_poisson}, by adopting the extended NB distribution with parameters $(\mu, p) \in [0, \infty)\times (0,1]$, if the data come from Poisson$(\lambda)$ with $\lambda > 0$, then $\hat\mu = \bar{Y}_n$~, which converges to $\lambda$ almost surely, and $\hat{p}$ converges to $1$ almost surely. That is, the MLE leads to the true Poisson distribution consistently. The major difference between the new parameterization and the old ones is that the new one contains Poisson distributions as a special class.

Similarly to Algorithm~\ref{algo:MLE_nu_p}, we propose the following Algorithm~\ref{algo:MLE_mu_p} to find the MLE $(\hat{\mu}, \hat{p})$ for the extended negative binomial distribution, which finds $\hat{p}$ with the aid of $\hat\nu$.

\begin{algorithm}
\caption{Adaptive Profile Maximization Algorithm for Extended NB($\mu, p$)
%: Finding MLE for Negative Binomial Distribution
}
\label{algo:MLE_mu_p}

Input: Data $Y_1, \ldots, Y_n \in \{0, 1, 2, \ldots\}$, predetermined threshold values $\nu_{\rm max}>0$, $\varepsilon>0$, and $\delta \in (0,1)$ (e.g., $\nu_{\rm max} = 10^4$, $\varepsilon = 10^{-3}$, $\delta=0.1$)

Output: MLE $(\hat\mu, \hat{p})$ for NB($\mu, p$)

\begin{algorithmic}[1]
\State Calculate the sample mean $\bar{Y}_n = \frac{1}{n} \sum_{i=1}^n Y_i$ and the sample variance $S^2 = \frac{1}{n-1} \sum_{i=1}^n (Y_i - \bar{Y}_n)^2$.
\State If $\bar{Y}_n=0$, go to Step~$7$ with $\hat\mu=0$ and $\hat{p}=1$; if $\bar{Y}_n \geq \frac{n-1}{n} S^2$, go to Step~$7$ with $\hat\mu = \bar{Y}_n$ and $\hat{p}=1$; otherwise calculate the initial value $\hat\nu^{(0)}$ of $\nu$ according to \eqref{eq:moment_estimate_nu}.
\State If the ratio of distinct values in the sample $|I|/n < \delta$ (a few-distinct-values scenario), use \eqref{eq:h(nu)} and \eqref{eq:g(nu)} for computing $h(\nu)$ and $g(\nu)$, respectively; otherwise (a many-distinct-values scenario), use \eqref{eq:h(nu)_psi} and \eqref{eq:g(nu)_psi} instead. 
\State Maximize $h(\nu)$ over $(\varepsilon, \nu_{\max}]$ starting at $\hat{\nu}^{(0)}$ using the L-BFGS-B algorithm \citep{byrd1995limited}, a quasi-Newton method allowing box constraints. The analytic gradient is $h'(\nu) = ng(\nu)$, and the optimization is constrained within $(\varepsilon, \nu_{\max}]$. The procedure terminates when either the maximum number of iterations (e.g., 500) is reached or the convergence tolerance of the L-BFGS-B algorithm is satisfied. The optimization can be implemented, for example, using the \texttt{optim} function in R.
\State Let $\nu_* =$ optimizer output. If $h(\nu_{\rm max}) > h(\nu_*)$ (see \eqref{eq:h(nu)_psi} or \eqref{eq:h(nu)}), we let $\hat\nu=\nu_{\rm max}$; otherwise, $\hat\nu=\nu_*$~. Return a warning message when $\hat\nu = \nu_{\rm max}$~.  
\State Calculate $\hat{p} = \hat\nu/(\hat\nu + \bar{Y}_n)$.
\State Report $\hat\mu = \bar{Y}_n$ and $\hat{p}$.
\end{algorithmic}
\end{algorithm}

Note that in Algorithm~\ref{algo:MLE_mu_p}, we still keep a threshold $\nu_{\rm max}>0$ as an upper bound for $\nu$ as an intermediate parameter. Different from Algorithm~\ref{algo:MLE_nu_p}, if $\bar{Y}_n\geq \frac{n-1}{n}S^2$, we conclude with a Poisson distribution directly without referring to $\nu_{\rm max}$~. In Algorithm~\ref{algo:MLE_mu_p}, $\nu_{\rm max}$ plays an effective role only when $\bar{Y}_n$ is close to but still less than $\frac{n-1}{n}S^2$. In this case, $\hat\nu$ tends to be unbounded and the optimization procedure becomes unstable. By setting a predetermined large number $\nu_{\rm max}$~, we stabilize the procedure. In the meantime, we also leave the users an option to enlarge $\nu_{\rm max}$ when it is attained.

\section{Numerical Studies}\label{sec:numerical_studies}

In this section, we use simulation studies to numerically confirm our findings.

\subsection{Comparison study on finding the MLE of NB distribution}\label{sec:comparison_MLE}

In this section, we compare the proposed Adaptive Profile Maximization Algorithm (APMA, Algorithm~\ref{algo:MLE_nu_p}) with seven other approaches that are currently available for finding the MLE of NB parameters, including the R function {\tt new.mle} in package AZIAD (version 0.0.2, \cite{dousti2022r}) recently developed for 27 probability distributions including NB, the R function {\tt fitdistr} in package MASS (version 7.3-61, \cite{venables2002modern}), a commonly used approach for finding MLE in R (\citealp{R-base}), and the R function {\tt fitdist} in package fitdistrplus (version 1.2-1, \cite{delignette2015fitdistrplus}) supporting a variety of classical optimization methods, namely ``Nelder-Mead'' \citep{nelder1965simplex}, ``BFGS'' (a quasi-Newton method), ``CG'' (a conjugate gradients method), ``L-BFGS-B'' (BFGS with box constraints), ``SANN'' (simulated annealing), and ``Brent''. We skip ``Brent'' in our comparison study since it is restricted to one-dimensional problems and not feasible for NB distributions. 

For this comprehensive simulation study, we consider 25 parameter pairs $(\nu, p)$ for negative binomial distributions, including $\nu = 0.01, 0.1, 1, 10, 100$ and $p = 0.99, 0.9, 0.5, 0.1,$ $0.01$. The complete results for all parameter pairs are reported in Section~\ref{sec:supp_tables} of the Supplementary Material. To avoid redundancy in the main text and to illustrate the algorithmic differences more clearly, we report in this section the results of five representative $(\nu, p)$ pairs with $p=0.9$ and $\nu = 0.01, 0.1, 1, 10, 100$. The choice of $p=0.9$ provides a Poisson-like scenario where the mean and variance are of comparable magnitude, making it a natural benchmark for assessing the performance when negative binomial data resemble a Poisson case. Compared with the more extreme case of $p=0.99$, where our algorithm achieves even more improvements, $p=0.9$ is less degenerate and therefore more representative. For smaller values of $p$ (e.g., $p=0.5, 0.1, 0.01$), our method still outperforms competing approaches, although its superiority is less striking.

Varying $\nu$ across $\{0.01, 0.1, 1, 10, 100\}$ further allows us to demonstrate the performance under different mean levels. Specifically, for $p=0.9$, the corresponding means (variances) are approximately $0.0011\,(0.0012)$, $0.0111\,(0.0123)$, $0.1111\,(0.1235)$, $1.1111$ $(1.2346)$, and $11.1111\,(12.3457)$, respectively. These five $(\nu, p)$ pairs, corresponding to a range of mean values from very small to moderately large, are then used to simulate $B=100$ NB$(\nu, p)$ random samples of size $n=100$ or $1,000$, respectively. For each random sample, we calculate the MLE of $(\nu, p)$ using APMA (Algorithm~\ref{algo:MLE_nu_p}), AZIAD, MASS ({\tt fitdistr}), as well as Nelder-Mead, BFGS, CG, L-BFGS-B, and SANN by utilizing {\tt fitdist} in R package fitdistrplus.

Table~\ref{tab:Error_rate} shows the relative frequency out of 100 simulated datasets from NB$(\nu, p=0.9)$ when an algorithm fails to return an MLE, namely the failure rate. For example, a failure occurs when the optimization encounters a non-finite finite-difference issue, preventing it from convergence and causing the procedure to terminate without producing an MLE. The results in Table~\ref{tab:Error_rate} indicate that APMA and AZIAD consistently succeed in returning MLEs across all values of $\nu$ and sample sizes. For the other six algorithms, failures are observed in various degrees depending on $\nu$ and the algorithms. When the sample size is small ($n=100$), methods such as MASS, Nelder-Mead, BFGS, CG, L-BFGS-B, and SANN show substantial failure rates for some values of $\nu$, with the most pronounced failures occurring for $\nu = 0.01$. Increasing the sample size to $n=1,000$ generally reduces the failure rates, but several methods including CG and L-BFGS-B still fail in a non-negligible rate across multiple $\nu$ values. In all five scenarios considered, failures occur when the sample mean is close to the sample variance, suggesting that this condition poses difficulties for these methods.

\begin{table}[hbt]
\caption{Failure rate out of 100 samples for which algorithms fail to return an MLE}
\label{tab:Error_rate}
\footnotesize
\begin{center}
\begin{tabular}{@{}lcccccccc}
\toprule
& \multicolumn{8}{@{}c}{\textbf{Simulations from NB$(\nu, p=0.9)$ with sample size} $n=100$}  \\\cmidrule{2-9}
$\bm{\nu}$&\textbf{APMA}&\textbf{AZIAD}&\textbf{MASS}&\textbf{Nelder-Mead}&\textbf{BFGS}&\textbf{CG}&\textbf{L-BFGS-B}&\textbf{SANN}\\
\midrule
\textbf{0.01}& 0   &  0   &  0.86   &  0.86   &  0.86   &  0.86   &  0.86   &  0.86 \\ 
\textbf{0.1}& 0   &  0   &  0.34   &  0.34   &  0.34   &  0.34   &  0.34   &  0.34 \\ 
\textbf{1}& 0   &  0   &  0.00   &  0.00   &  0.23   &  0.81   &  0.35   &  0.00 \\ 
\textbf{10}& 0   &  0   &  0.23   &  0.00   &  0.37   &  1.00   &  0.55   &  0.00 \\ 
\textbf{100}& 0   &  0   &  0.20   &  0.00   &  0.00   &  1.00   &  0.16   &  0.00 \\ 
\bottomrule

\toprule
& \multicolumn{8}{@{}c}{\textbf{Simulations from NB$(\nu, p=0.9)$ with sample size} $n=1,000$}  \\\cmidrule{2-9}
$\bm{\nu}$&\textbf{APMA}&\textbf{AZIAD}&\textbf{MASS}&\textbf{Nelder-Mead}&\textbf{BFGS}&\textbf{CG}&\textbf{L-BFGS-B}&\textbf{SANN}\\
\midrule
\textbf{0.01}& 0   &  0   &  0.71   &  0.74   &  0.71   &  0.71   &  0.98   &  0.84 \\ 
\textbf{0.1}& 0   &  0   &  0.00   &  0.00   &  0.00   &  0.01   &  0.49   &  0.01 \\ 
\textbf{1}& 0   &  0   &  0.00   &  0.00   &  0.10   &  0.87   &  0.14   &  0.00 \\ 
\textbf{10}& 0   &  0   &  0.00   &  0.00   &  0.11   &  0.99   &  0.19   &  0.00 \\ 
\textbf{100}& 0   &  0   &  0.46   &  0.00   &  0.00   &  1.00   &  0.03   &  0.00 \\ 
\bottomrule
\end{tabular}
\end{center}
Note: For example, the value ``0.86'' indicates that the corresponding algorithm fails to find the MLE in 86 out of 100 simulations.
\normalsize
\end{table}

To compare the efficiency or accuracy of the algorithms, we adopt the maximum likelihood ratio, that is, $L(\hat\nu, \hat{p})/L(\hat\nu_{\rm APMA}, \hat{p}_{\rm APMA}) = \exp\{l(\hat\nu, \hat{p})-l(\hat\nu_{\rm APMA}, \hat{p}_{\rm APMA})\}$, which is computed only when $(\hat\nu, \hat{p})$ is successfully obtained. Table~\ref{tab:ratio_maximum_likelihood} reports the average of these ratios taken over the successful runs. The results show that, even conditional on returned parameter estimates, AZIAD and Nelder-Mead achieve likelihood values comparable to those of our proposed algorithm, whereas the other methods generally do not. To examine the performance of AZIAD in more detail, a comparison between APMA and AZIAD is provided subsequently in Section~\ref{sec:compare_APMA_and_AZIAD}. In contrast, the performance of Nelder-Mead is influenced by its ability to explore an unbounded parameter space, which may lead to larger estimates of $\nu$ and likelihood values comparable to those of our algorithm. When Nelder-Mead is able to converge to a solution, it may yield estimates of $\nu$ exceeding our pre-specified upper limit of $\nu_{\max}=10{,}000$; however, the resulting estimates are typically numerically unstable. For the dataset considered in this section, the likelihood continues to increase as $\nu$ grows larger. We impose a finite upper bound $\nu_{\max}$ solely to ensure numerical stability and computational feasibility, which can be adjusted if needed.

\begin{table}[hbt]
\caption{Average ratio of maximum likelihood over 100 simulations with $p=0.9$}
\label{tab:ratio_maximum_likelihood}
\footnotesize
\begin{center}
\begin{tabular}{@{}lcccccccc}
\toprule
& \multicolumn{8}{@{}c}{\textbf{Average ratio of maximum likelihood with} $n=100$}  \\\cmidrule{2-9}
$\bm{\nu}$&\textbf{APMA}&\textbf{AZIAD}&\textbf{MASS}&\textbf{Nelder-Mead}&\textbf{BFGS}&\textbf{CG}&\textbf{L-BFGS-B}&\textbf{SANN}\\
\midrule
\textbf{0.01}& 1   &  1   &  1   &  1   &  1   &  1   &  1   &  1 \\
\textbf{0.1}& 1   &  1   &  0.91   &  1   &  0.91   &  0.91   &  0.91   &  0.94 \\ 
\textbf{1}& 1   &  1   &  0.94   &  1   &  0.93   &  0.99   &  0.92   &  0.85 \\ 
\textbf{10}& 1   &  1   &  1   &  1   &  0.95   &  NA   &  1   &  0.88 \\ 
\textbf{100}& 1   &  0.99   &  0.97   &  1   &  0.86   &  NA   &  1   &  0.16 \\ 
\bottomrule

\toprule
& \multicolumn{8}{@{}c}{\textbf{Average ratio of maximum likelihood with} $n=1,000$}  \\\cmidrule{2-9}
$\bm{\nu}$&\textbf{APMA}&\textbf{AZIAD}&\textbf{MASS}&\textbf{Nelder-Mead}&\textbf{BFGS}&\textbf{CG}&\textbf{L-BFGS-B}&\textbf{SANN}\\
\midrule
\textbf{0.01}& 1   &  1   &  0.69   &  1   &  0.69   &  0.69   &  0.02   &  0.57 \\
\textbf{0.1}& 1   &  1   &  0.70   &  1   &  0.70   &  0.70   &  0.49   &  0.72 \\ 
\textbf{1}& 1   &  1   &  0.46   &  1   &  0.39   &  0.93   &  0.40   &  0.61 \\ 
\textbf{10}& 1   &  1   &  0.98   &  1   &  0.95   &  1   &  0.99   &  0.89 \\ 
\textbf{100}& 1   &  0.99   &  0.99   &  1   &  0.52   &  NA   &  1   &  $<0.01$ \\ 
\bottomrule
\end{tabular}
\end{center}
Note: (1) The ratio is the likelihood of an algorithm's MLE divided by the likelihood of the APMA's MLE, averaged over simulations when the algorithm successfully returns an estimate; (2) ``NA'' indicates that the corresponding algorithm fails to find an MLE in all 100 simulations.
\normalsize
\end{table}

In terms of computational time (see Table~\ref{tab:Average_computational_time}), APMA (Algorithm~\ref{algo:MLE_nu_p}) demonstrates a significant advantage compared to all other methods. All algorithms are implemented through routines that rely on the \texttt{optim} function in R and are run under its default convergence tolerance (\texttt{reltol = 1e-8}), ensuring that the reported time reflects comparable accuracy requirements across different methods. APMA and AZIAD both employ the L-BFGS-B method with parameter bounds \texttt{lower = 1e-4} and \texttt{upper = 1e4}. MASS employs the BFGS method, and the other algorithms (e.g., Nelder--Mead, CG, SANN) use the optimization methods indicated by their names, all under their default settings and sharing the same convergence tolerance. Although AZIAD does nearly as well as APMA in terms of the failure rates and the average ratios of maximum likelihood, it costs much more computational time than APMA. Overall, we recommend our APMA (Algorithm~\ref{algo:MLE_nu_p}).

\begin{table}[hbt]
\caption{Average computational time (sec) over 100 simulations with $p=0.9$, $n=1,000$}
\label{tab:Average_computational_time}
\footnotesize
\begin{center}
\begin{tabular}{@{}lcccccccc@{}}
\toprule
$\bm{\nu}$&\textbf{APMA}&\textbf{AZIAD}&\textbf{MASS}&\textbf{Nelder-Mead}&\textbf{BFGS}&\textbf{CG}&\textbf{L-BFGS-B}&\textbf{SANN}\\
\midrule
\textbf{0.01}& 0.0008   &  0.0294   &  0.0017   &  0.0341   &  0.0078   &  0.0125   &  0.0153   &  1.6548 \\ 
\textbf{0.1}& 0.0010   &  0.0355   &  0.0042   &  0.0205   &  0.0062   &  0.0230   &  0.0195   &  1.1223 \\ 
\textbf{1}& 0.0006   &  0.0113   &  0.0051   &  0.0155   &  0.0102   &  0.0739   &  0.0126   &  1.2124 \\ 
\textbf{10}& 0.0005   &  0.0125   &  0.0110   &  0.0139   &  0.0170   &  0.1052   &  0.0150   &  2.4412 \\ 
\textbf{100}& 0.0015   &  0.0349   &  0.1868   &  0.0609   &  0.0491   &  0.3469   &  0.0650   &  6.2008 \\ 
\bottomrule
\end{tabular}
\end{center}
\normalsize
\end{table}

\subsection{Further comparison between APMA and AZIAD}\label{sec:compare_APMA_and_AZIAD}

From Tables~\ref{tab:Error_rate} and \ref{tab:ratio_maximum_likelihood}, we observe that APMA and AZIAD outperform others in terms of failure rates and average maximum likelihood ratios, when the data are simulated from negative binomial distributions. 
In this section, we further compare the efficiency of APMA and AZIAD when they encounter Poisson samples. 

\begin{table}[hbt]
\caption{Average ratio of maximum likelihood obtained by AZIAD against APMA on $100$ random samples from Poisson($\lambda$) of size $n$}
\label{tab:Likelihood_Comparison_Algo1_Nil}
\begin{center}
\begin{tabular}{@{}lrrrr@{}}
\toprule
& \multicolumn{4}{@{}c}{\textbf{n}}  \\\cmidrule{2-5}
$\bm{\lambda}$&\textbf{50}&\textbf{500}&\textbf{5,000}&\textbf{50,000}\\
\midrule
\textbf{1}  & 0.996 & 0.978 & 0.871 & 0.476  \\
\textbf{3}  & 0.991 & 0.945 & 0.729 & 0.166   \\
\textbf{5}  & 0.985 & 0.930 & 0.679 & 0.061   \\
\textbf{10} & 0.980 & 0.885 & 0.443 & 0.007  \\
\bottomrule
\end{tabular}
\end{center}
\end{table}

From Table~\ref{tab:Likelihood_Comparison_Algo1_Nil}, we can see that APMA outperforms AZIAD in terms of finding higher maximum likelihood when the data come from a Poisson distribution. The advantage of APMA over AZIAD is greater as the sample size $n$ or the sample mean becomes more pronounced.

\subsection{Performance across three types of data regimes}\label{sec:real_data}

Based on our experiences, APMA (Algorithm~\ref{algo:MLE_nu_p}) performs particularly well in three types of data regimes. In this section, we examine the performance of all eight algorithms compared in Section~\ref{sec:comparison_MLE} on representative datasets from existing R packages, which serve as examples and complement the evidence from numerical experiments. For each regime, we then explain why our algorithm is particularly effective given the corresponding data characteristics.

The first type of examples involves nearly equidispersed count data, where the sample mean and variance are of comparable magnitude. This situation can arise under certain parameter combinations of the negative binomial distribution, typically when $p$ is large, e.g., the simulated data used in Section~\ref{sec:comparison_MLE}. In such cases, the MLE of the negative binomial distribution often becomes numerically unstable. We illustrate this with the \texttt{prussian} dataset \citep{bortkiewicz1898} from R package \texttt{pscl} (version 1.5.9, \cite{pscl2024package}), which records the annual death counts caused by horse kicks in 14 cavalry corps of the Prussian army from 1875 to 1894. This dataset contains 280 observations, with a sample mean of 0.7 and a sample variance of 0.7627. In this example, the CG algorithm fails to return an MLE. Among all the eight methods, only APMA and L-BFGS-B achieve the highest log-likelihood with the shortest computational time, and produce an estimated mean that matches the sample mean perfectly (see Table~\ref{tab:estimation_prussian} in Section~\ref{sec:supp_tables} of the Supplementary Material).

The advantage of APMA for handling nearly equidispersed count data lies in solving the NB MLE over a bounded interval $(\varepsilon, \nu_{\max}]$. Such kind of data may arise either from a genuine NB model with a large $\nu$ or from a Poisson model. If the data are essentially Poisson, the estimate $\hat{\nu}$ tends to be unbounded, according to the theoretical results in Sections~\ref{sec:NB_given_Poisson} and~\ref{sec:other_parameterizations}. When $\nu$ is very large, an unrestricted optimization may fail to converge, return unstable estimates, or yield suboptimal solutions with lower likelihood values. It can also be computationally costly. By setting a large $\nu_{\max}$~, APMA ensures stable and fast estimation for this type of data. If its solution reaches $\nu_{\max}$~, the algorithm will issue a warning to the users.

The second type of examples involves  extremely sparse data, typically with a small $\nu$ and a $p$ close to $1$. The data in this regime consist almost entirely of zeros, whose sample mean and sample variance are both small. As shown in Section~\ref{sec:supp_tables} of the Supplementary Material, this pattern is evident for settings with small $\nu$ (e.g., $0.01,0.1$) and large $p$ (e.g., $0.99,0.9$). For these cases, all routines except APMA and AZIAD may fail to produce stable parameter estimates.

For illustration purposes, we examine the numbers of insurance claims in the \texttt{dataCar} dataset \citep{dejong2008} from R package \texttt{insuranceData} (version 1.0, \cite{insurance2014package}), whose ratio of zeros is as high as 93.19\% out of 67,856 observations. For this dataset, only APMA, AZIAD, and Nelder-Mead achieve the highest log-likelihood, and APMA requires far less computational time than all the others (see Table~\ref{tab:estimation_dataCar} in Section~\ref{sec:supp_tables} of the Supplementary Material). It shows clear advantage of our algorithm with sparse data. The good performance of APMA in this regime stems from its  branching strategy that treats sparse and dense cases with different formulas (see Step~3 of Algorithm~\ref{algo:MLE_nu_p}), which improves computational efficiency significantly under sparsity.

The third type of examples involve overdispersed data with either very small or very large counts, which causes additional difficulties in computing the MLE of an NB distribution. As shown in Tables~\ref{tab:Error_rate_supplementary_100} and \ref{tab:Error_rate_supplementary_1000} of Section~\ref{sec:supp_tables} in the Supplementary Material, when the variance substantially exceeds the mean (e.g., $p \leq 0.5$) and the scale of the counts is small (e.g., $\nu=0.01,0.1$), the L-BFGS-B algorithm exhibits a high failure rate. When the scale is extremely large (e.g., $\nu=100$), the CG algorithm consistently fails. We use two empirical datasets to illustrate these issues. The \texttt{epil} dataset \citep{thall1990} from the \texttt{MASS} package contains seizure counts from an epilepsy clinical trial and represents a typical case of overdispersed count data. It has a sample mean of $8.26$, a sample variance of $152.68$, and a sample size of $236$. In this case, the L-BFGS-B algorithm fails. The \texttt{UKDriverDeaths} dataset \citep{harvey1986} from R package \texttt{datasets} (version 4.6.0, \cite{R-base}) records monthly counts of car drivers killed or seriously injured in the United Kingdom from 1969 to 1984. It has a sample mean of $1,670$, a sample variance of $83,875$, and a sample size of $192$. In this case, the CG algorithm fails to converge. For both examples, APMA provides stable parameter estimates, highlighting its robustness in highly overdispersed regimes. 

Taken together, these examples confirm that APMA succeeds in data regimes where other competing methods may struggle. Its particular strength lies in Poisson-like settings, in line with our theoretical results developed earlier, while its ability to handle sparse or highly overdispersed data illustrates its further robustness that enhances its practical significance.

\section{Conclusion}\label{sec:conclusion}

For certain datasets, traditional algorithms often encounter difficulties in computing the MLE of the NB distribution. To address this issue, we theoretically investigate the asymptotic properties of the MLE and the resulting distribution when the data come from a Poisson distribution. As a solution, we propose a new algorithm for finding the MLE of NB distributions, and introduce a new reparameterization that naturally incorporates the Poisson distribution as a special case.

Our algorithm demonstrates improved solvability, accuracy, and efficiency in scenarios where existing methods frequently fail, while the reparameterization of NB distributions provides a clearer interpretation of model parameters. Together, these contributions enhance both the theoretical understanding and practical applicability of negative binomial models, providing a solid foundation for future methodological developments.

\clearpage
\setcounter{page}{1}
\def\thepage{S\arabic{page}}

\centerline{\large\bf From Poisson Observations to Fitted Negative Binomial Distribution}
\vspace{.25cm}
 \centerline{Yingying Yang$^1$, Niloufar Dousti Mousavi$^2$, Zhou Yu$^3$, and Jie Yang$^1$} 
\vspace{.4cm}
 \centerline{\it University of Illinois at Chicago$^1$, University of Pittsburgh$^2$, and Merck \& Co., Inc.$^3$}
\vspace{.55cm}
 \centerline{\bf Supplementary Material}
\vspace{.55cm}
%\fontsize{9}{11.5pt plus.8pt minus .6pt}\selectfont
%\noindent
%CONTENT OF A BRIEF NOTE.........\\
%CONTENT OF THE BRIEF NOTE.\\
%CONTENT OF THE BRIEF NOTE.\\
\par

\setcounter{section}{0}
\setcounter{equation}{0}
\def\theequation{S.\arabic{equation}}
\def\thesection{S\arabic{section}}

\setcounter{figure}{0}
\def\thefigure{S.\arabic{figure}}
\setcounter{table}{0}
\def\thetable{S.\arabic{table}}

%\fontsize{12}{14pt plus.8pt minus .6pt}\selectfont

\begin{itemize}
\item[{\bf S1}] {\bf Proofs:} We provide proofs for lemmas and theorems in the main text.
\item[{\bf S2}] {\bf A simulation study on using KS test for NB against Poisson distributions:}
We use the Kolmogorov-Smirnov (KS) test to illustrate that it is difficult to separate a Poisson law from a nearly equidispersed NB distribution.
\item[{\bf S3}] {\bf Asymptotic properties of KS test statistic:} We theoretically derive the asymptotic properties of the KS test statistic for testing NB against Poisson distributions, given that the data come from a Poisson distribution.
\item[{\bf S4}] {\bf More tables for Section~\ref{sec:numerical_studies}:} We provide more tables for supporting the discussion in Section~\ref{sec:numerical_studies} (Numerical Studies). 
\item[{\bf S5}] {\bf R code for Algorithm~\ref{algo:MLE_nu_p} and Algorithm~\ref{algo:MLE_mu_p}:} We provide R code for implementing Algorithms~\ref{algo:MLE_nu_p} \& \ref{algo:MLE_mu_p}, as well as applications to three examples.
\end{itemize}

\section{Proofs}\label{sec:proofs}

\medskip\noindent
{\bf Proof of Theorem~\ref{thm:M=0}:}
If $M=0$, as mentioned in Example~\ref{ex:all_zero}, $l(\nu, p) = n\nu\log p\leq 0$. To maximize $l(\nu, p)$ with $\nu \in (0, \nu_{\rm max}]$ and $p\in (0, 1]$, we must have $\hat{p}=1$. In this case, $\hat{\nu}$ can be any positive number in $(0, \nu_{\rm max}]$.

On the other hand, if $M>0$, then $\bar{Y}_n >0$. For any fixed $\nu>0$, the best $p$ must take the form of $\nu/(\nu+\bar{Y}_n)$ as in \eqref{eq:p_as_nu}. Then $\hat{p} \leq \nu_{\rm max}/(\nu_{\rm max} + \bar{Y}_n) < 1$.
\hfill{$\Box$}

\medskip\noindent
{\bf Proof of Lemma~\ref{lem:g(nu)}:}
Since $\bar{Y}_n>0$, then there exists a $y\in I\setminus \{0\}$, such that, $f_y > 0$, or equivalently, $f_0<n$. It can be verified that
\[
\nu g(\nu) = \frac{1}{n}\sum_{y \in I\setminus \{0\}} f_y + \frac{\nu}{n} \sum_{y \in I\setminus \{0,1\}} f_y\left(\frac{1}{\nu+1} + \cdots + \frac{1}{\nu+y-1}\right) - \nu\log\left(1+\frac{\bar{Y}_n}{\nu}\right).
\]
Then $\lim_{\nu\to 0^+} \nu g(\nu) = \frac{1}{n}\sum_{y \in I\setminus \{0\}} f_y = 1-\frac{f_0}{n} >0$.
As a direct conclusion, 
$\lim_{\nu\rightarrow 0^+} g(\nu) = \infty$. Similarly, since
\[
\nu^2 g'(\nu) = -\frac{1}{n}\sum_{y \in I\setminus \{0\}} f_y - \frac{\nu^2}{n} \sum_{y \in I\setminus \{0,1\}} f_y\left[\frac{1}{(\nu+1)^2} + \cdots + \frac{1}{(\nu+y-1)^2}\right] + \frac{\nu\bar{Y}_n}{\nu+\bar{Y}_n},
\]
and thus 
$\lim_{\nu \rightarrow 0^+} \nu^2 g'(\nu)  = -\frac{1}{n}\sum_{y \in I\setminus \{0\}} f_y = \frac{f_0}{n}-1 <0$, then $\lim_{\nu\rightarrow 0^+} g'(\nu) = -\infty$. 

On the other hand, as $\nu\rightarrow \infty$, by applying Taylor's theorem with the Peano form of the reminder (see, e.g., Section~7.9 in \cite{apostol1991calculus}), we have 
\begin{eqnarray*}
\frac{1}{\nu+k} &=& \frac{\frac{1}{\nu}}{1+\frac{k}{\nu}} = \frac{1}{\nu} - \frac{k}{\nu^2} + o(\nu^{-2}) \\
\log\left(1+\frac{\bar{Y}_n}{\nu}\right) &=&  \frac{\bar{Y}_n}{\nu} - \frac{\bar{Y}_n^2}{2\nu^2} + o(\nu^{-2}) \\
\frac{1}{(\nu + k)^2} &=& \frac{\frac{1}{\nu^2}}{\left(1+\frac{k}{\nu}\right)^2} = \frac{1}{\nu^2} - \frac{2k}{\nu^3} + o(\nu^{-3}) \\
\frac{\bar{Y}_n}{\nu (\nu+\bar{Y}_n)} &=& \frac{\bar{Y}_n}{\nu^2} - \frac{\bar{Y}_n^2}{\nu^3} + o(\nu^{-3})
\end{eqnarray*}
Then, as $\nu\rightarrow \infty$, by letting $S^2_n = n^{-1}\sum_{i=1}^n (Y_i - \bar{Y}_n)^2 = n^{-1}\sum_{i=1}^n Y_i^2 - \bar{Y}_n^2$,
\begin{eqnarray*}
g(\nu) &=& \frac{1}{n} \sum_{y \in I\setminus \{0\}} f_y \left(\frac{1}{\nu} + \cdots + \frac{1}{\nu + y - 1}\right) - \log\left(1 + \frac{\bar{Y}_n}{\nu}\right)\\
&=& \frac{1}{n} \sum_{y \in I\setminus \{0\}} f_y \left[\frac{y}{\nu} - \frac{y(y-1)}{2\nu^2}\right] - \frac{\bar{Y}_n}{\nu} + \frac{\bar{Y}_n^2}{2\nu^2} + o(\nu^{-2})\\
&=& \frac{1}{n} \sum_{y \in I\setminus \{0\}} \frac{yf_y}{\nu} - \frac{1}{n} \sum_{y \in I\setminus \{0\}} \frac{y^2f_y}{2\nu^2} + \frac{1}{n} \sum_{y \in I\setminus \{0\}} \frac{yf_y}{2\nu^2} - \frac{\bar{Y}_n}{\nu} + \frac{\bar{Y}_n^2}{2\nu^2} + o(\nu^{-2})\\
&=& \frac{\bar{Y}_n}{\nu} - \frac{1}{2\nu^2} \left(S^2_n + \bar{Y}_n^2\right) + \frac{\bar{Y}_n}{2\nu^2} - \frac{\bar{Y}_n}{\nu} + \frac{\bar{Y}_n^2}{2\nu^2} + o(\nu^{-2})\\
&=& \frac{1}{2\nu^2}\left(\bar{Y}_n - S^2_n\right) + o(\nu^{-2}),
\end{eqnarray*}
\begin{eqnarray*}
g'(\nu) &=& - \frac{1}{n} \sum_{y \in I\setminus \{0\}} f_y \left[\frac{1}{\nu^2} + \cdots + \frac{1}{(\nu + y - 1)^2}\right] + \frac{\bar{Y}_n}{\nu (\nu+\bar{Y}_n)} \\
&=& - \frac{1}{n} \sum_{y \in I\setminus \{0\}} f_y \left[\frac{y}{\nu^2} - \frac{y(y-1)}{\nu^3}\right] + \frac{\bar{Y}_n}{\nu^2} - \frac{\bar{Y}_n^2}{\nu^3} + o(\nu^{-3}) \\
&=& - \frac{1}{\nu^2} \sum_{y \in I\setminus \{0\}} \frac{yf_y}{n} + \frac{1}{\nu^3} \sum_{y \in I\setminus \{0\}} \frac{y^2f_y}{n} - \frac{1}{\nu^3} \sum_{y \in I\setminus \{0\}} \frac{yf_y}{n} + \frac{\bar{Y}_n}{\nu^2} - \frac{\bar{Y}_n^2}{\nu^3} + o(\nu^{-3}) \\
&=& - \frac{\bar{Y}_n}{\nu^2} + \frac{1}{\nu^3} \left(S^2_n + \bar{Y}_n^2\right) - \frac{\bar{Y}_n}{\nu^3} + \frac{\bar{Y}_n}{\nu^2} - \frac{\bar{Y}_n^2}{\nu^3} + o(\nu^{-3})\\
&=& \frac{1}{\nu^3}\left(S^2_n - \bar{Y}_n\right) + o(\nu^{-3}).
\end{eqnarray*}
In other words, we obtain $\lim_{\nu\rightarrow\infty} {\nu}^2 g(\nu) = \frac{1}{2} \left(\bar{Y}_n - S^2_n\right)$ and $\lim_{\nu\rightarrow\infty} {\nu}^3 g'(\nu) = S^2_n - \bar{Y}_n$, which imply $\lim_{\nu\rightarrow\infty} g(\nu) = 0$ and $\lim_{\nu\rightarrow\infty} g'(\nu) = 0$, respectively.
\hfill{$\Box$}

\medskip\noindent
{\bf Proof of Theorem~\ref{thm:min_g^2_positive}:}
Let $\nu_* = {\rm argmin}_{\nu \in (0, \nu_{\rm max}]} g^2(\nu)$, then $g(\nu_*)\neq 0$. According to Lemma~\ref{lem:g(nu)}, $\lim_{\nu\rightarrow 0^+} g(\nu) = \infty$. With the constraint $\nu \in (0, \nu_{\rm max}]$, we must have $\nu_*>0$. 

We claim that $g(\nu_*)>0$. Actually, if $g(\nu_*) < 0$, then there exists a $\nu_{**} \in (0, \nu_*)$, such that, $g(\nu_{**})=0$, which violates $\min_{\nu \in (0, \nu_{\rm max}]} g^2(\nu) > 0$.

The fact $g(\nu_*)>0$ implies $g(\nu) \geq g(\nu_*)>0$ for all $\nu \in (0, \nu_{\rm max}]$. That is, $h(\nu)$ increases all the time and thus the MLE of $\nu$ is $\nu_{\rm max}$~.
\hfill{$\Box$}

\medskip\noindent
{\bf Proof of Lemma~\ref{lem:G_F(nu)_1-F}:} 
For any $\nu_0>0$, once $\nu>\nu_0/2$, 
\begin{eqnarray*}
& &\left|\sum_{y=1}^\infty P(Y_1=y) \left(\frac{1}{\nu} + \cdots + \frac{1}{\nu + y - 1}\right) - \sum_{y=1}^\infty P(Y_1=y) \left(\frac{1}{\nu_0} + \cdots + \frac{1}{\nu_0 + y - 1}\right)\right|\\
&\leq & \sum_{y=1}^\infty P(Y_1=y) \left[\left|\frac{1}{\nu}-\frac{1}{\nu_0}\right| + \cdots + \left|\frac{1}{\nu + y-1} - \frac{1}{\nu_0 + y-1}\right|\right]\\
&=& \sum_{y=1}^\infty P(Y_1=y) \left[\frac{|\nu - \nu_0|}{\nu\nu_0} + \cdots + \frac{|\nu-\nu_0|}{(\nu+y-1)(\nu_0+y-1)}\right]\\
&\leq & \sum_{y=1}^\infty P(Y_1=y)\cdot y \cdot \frac{|\nu-\nu_0|}{\nu_0^2/2}\\
&=& \frac{2 E(Y_1)}{\nu_0^2}\cdot |\nu-\nu_0|\\
&\rightarrow & 0\ ,
\end{eqnarray*}
as $\nu$ goes to $\nu_0$~. Since $\log(1+\mu/\nu)$ is a continuous function of $\nu$, $G_F(\nu)$ is a continuous function of $\nu$ on $(0, \infty)$, according to \eqref{eq:G_F(nu)}.

Furthermore, 
\begin{eqnarray*}
G_F(\nu) 
&=& \sum_{y=1}^\infty P(Y_1=y) \left(\frac{1}{\nu} + \cdots + \frac{1}{\nu + y - 1}\right) - \log\left(1 + \frac{\mu}{\nu}\right)\\
&=& \sum_{y=1}^\infty P(Y_1=y) \sum_{k=0}^{y-1} \frac{1}{\nu+k} - \log\left(1 + \frac{\mu}{\nu}\right)\\
&=& \sum_{k=0}^\infty \sum_{y=k+1}^\infty \frac{P(Y_1=y)}{\nu+k} - \log\left(1 + \frac{\mu}{\nu}\right)\\
&=& \sum_{k=0}^\infty \frac{1-P(Y_1\leq k)}{\nu+k} - \log\left(1 + \frac{\mu}{\nu}\right)\\
&=& \sum_{k=0}^\infty \frac{1-F(k)}{\nu+k} - \log\left(1 + \frac{\mu}{\nu}\right)\ .
\end{eqnarray*}
Thus the conclusion follows.
\hfill{$\Box$}

\medskip\noindent
{\bf Proof of Theorem~\ref{thm:G(nu)}:} 
According to Kolmogorov's Strong Law of Large Numbers (see, e.g., Theorem 7.5.1 in \cite{resnick1999probpath}), since $E(|Y_1|) = E(Y_1) = \mu < \infty $, then $\bar{Y}_n = \frac{1}{n} \sum_{i=1}^n Y_i \rightarrow \mu$ almost surely, as $n\rightarrow \infty$. 

Similarly, for each $y \in \{0, 1, 2, \ldots\}$, $E(|\mathds{1}_{\{Y_1 = y\}}|) = P(Y_1=y) \leq 1 < \infty $. Then
\begin{eqnarray*}
\frac{f_y}{n} &=& \frac{\#\{1\leq i\leq n\mid Y_i=y\}}{n} \quad =\quad \frac{1}{n} \sum_{i=1}^{n} \mathds{1}_{\{Y_i = y\}} \\
&\longrightarrow& E \mathds{1}_{\{Y_1 = y\}} = P(Y_1=y)
\end{eqnarray*}
almost surely, as $n\rightarrow \infty$. Furthermore, for each integer $M>0$, $E(|Y_1\mathds{1}_{\{Y_1 >M\}}|) = E(Y_1\mathds{1}_{\{Y_1 >M\}}) \leq E(Y_1) = \mu < \infty $. Then
\begin{eqnarray*}
\sum_{y=M+1}^\infty \frac{yf_y}{n} &=& \sum_{y=M+1}^\infty \frac{\#\{1\leq i\leq n\mid Y_i=y\}\cdot y}{n} \quad =\quad \frac{1}{n} \sum_{i=1}^{n} Y_i \mathds{1}_{\{Y_i > M\}} \\
&\longrightarrow& E (Y_1\mathds{1}_{\{Y_1 > M\}})
\end{eqnarray*}
almost surely, as $n\rightarrow \infty$. Then there exists an event $A$ with probability $1$, such that, for each $\omega \in A$, 
\begin{equation}\label{eq:limit_fy/n}
\frac{f_y(\omega)}{n} = \frac{\#\{1\leq i\leq n\mid Y_i(\omega)=y\}}{n} \longrightarrow P(Y_1=y)
\end{equation}
for each $y\in \{0, 1, 2, \ldots \}$, $\bar{Y}_n(\omega) \rightarrow \mu$ and thus $\log\left(1+\frac{\bar{Y}_n(\omega)}{\nu}\right) \rightarrow \log\left(1+\frac{\mu}{\nu}\right)$, and
\begin{equation}\label{eq:limit_sum_M+1}
\sum_{y=M+1}^\infty \frac{yf_y(\omega)}{n} \rightarrow E (Y_1\mathds{1}_{\{Y_1 > M\}})
\end{equation}
for each $M>0$ as the convergence of a sequence of real numbers. 

Now we only need to show that 
\begin{equation}\label{eq:limt_sum_fy_v}
\sum_{y=1}^{\infty} \frac{f_y(\omega)} {n} \left(\frac{1}{\nu} + \cdots + \frac{1}{\nu + y - 1} \right)  \longrightarrow \sum_{y=1}^{\infty} P(Y_1=y) \left(\frac{1}{\nu} + \cdots + \frac{1}{\nu + y - 1} \right) 
\end{equation}
for each $\omega\in A$.
That is, we want to show that, for any $\epsilon > 0$, there exists an $N$, such that, 
\begin{equation}\label{eq:lim_sum_fy-sum_p}
\left|
\sum_{y=1}^{\infty} \frac{f_y(\omega)} {n} \left(\frac{1}{\nu} + \cdots + \frac{1}{\nu + y - 1} \right)  - \sum_{y=1}^{\infty} P(Y_1=y) \left(\frac{1}{\nu} + \cdots + \frac{1}{\nu + y - 1} \right)\right| < \epsilon 
\end{equation}
for all $n>N$.

First, there exists an $M>0$, such that,
\begin{equation}\label{eq:sum_P_M+1<e/3}
\sum_{y=M+1}^{\infty} P(Y_1=y) \left(\frac{1}{\nu} + \cdots + \frac{1}{\nu + y - 1} \right) < \frac{\epsilon}{3}
\end{equation}
and 
\[
E (Y_1\mathds{1}_{\{Y_1 > M\}}) < \frac{\nu\epsilon}{6}\ .
\]
Then for each $y=1, \ldots, M$, there exists an $N_y>0$, such that, 
\[
\left|\frac{f_y(\omega)}{n} - P(Y_1=y)\right| < \frac{\epsilon}{3} \cdot \frac{\nu}{M} \cdot \frac{1}{2^y} \ ,
\]
as long as $n>N_y$~. Then for each $n>\max\{N_1, \ldots, N_M\}$, 
\begin{eqnarray}
    & & \left|
\sum_{y=1}^{M} \frac{f_y(\omega)} {n} \left(\frac{1}{\nu} + \cdots + \frac{1}{\nu + y - 1} \right)  - \sum_{y=1}^{M} P(Y_1=y) \left(\frac{1}{\nu} + \cdots + \frac{1}{\nu + y - 1} \right)\right| \nonumber\\
& \leq & \sum_{y=1}^{M} \left|\frac{f_y(\omega)} {n} - P(Y_1=y)\right| \left(\frac{1}{\nu} + \cdots + \frac{1}{\nu + y - 1} \right)\nonumber\\
& < & \sum_{y=1}^{M} \frac{\epsilon}{3} \cdot \frac{\nu}{M} \cdot \frac{1}{2^y} \cdot \frac{M}{\nu}\nonumber\\
& < & \frac{\epsilon}{3} \ .
\label{eq:sum_1_M<e/3}
\end{eqnarray}
We still need to show that 
\[
\sum_{y=M+1}^{\infty} \frac{f_y(\omega)} {n} \left(\frac{1}{\nu} + \cdots + \frac{1}{\nu + y - 1} \right) < \frac{\epsilon}{3}
\]
for large enough $n$. Actually, there exists an $N_\nu > 0$ such that
\[
\left|\sum_{y=M+1}^\infty \frac{yf_y(\omega)}{n} - E (Y_1\mathds{1}_{\{Y_1 > M\}})\right| < \frac{\nu\epsilon}{6}
\]
for all $n > N_\nu$~. Then
\begin{eqnarray}
& & \sum_{y=M+1}^{\infty} \frac{f_y(\omega)} {n} \left(\frac{1}{\nu} + \cdots + \frac{1}{\nu + y - 1} \right) \nonumber\\
&<& \sum_{y=M+1}^{\infty} \frac{f_y(\omega)}{n}\frac{y}{\nu} \nonumber\\
&<& \frac{1}{\nu} \left(E\left(Y_1 {\mathbf 1}_{Y_1>M}\right) + \frac{\nu\epsilon}{6}\right)\nonumber\\
&<& \frac{1}{\nu} \left(\frac{\nu\epsilon}{6} + \frac{\nu\epsilon}{6}\right)\nonumber\\
&=& \frac{\epsilon}{3}\ . \label{eq:sum_M+1<e/3}
\end{eqnarray}

Combining \eqref{eq:sum_P_M+1<e/3}, \eqref{eq:sum_1_M<e/3}, and \eqref{eq:sum_M+1<e/3}, we obtain that \eqref{eq:lim_sum_fy-sum_p} is true for all $n >$ $\max\{N_1, \ldots, $ $N_M, N_\nu\}$. That is, \eqref{eq:limt_sum_fy_v} is true for all $\omega \in A$. Combining with $\log\left(1+\frac{\bar{Y}_n(\omega)}{\nu}\right) \rightarrow \log\left(1+\frac{\mu}{\nu}\right)$ for all $\omega\in A$, we obtain $g(\nu) \rightarrow G_F(\nu)$ almost surely.
\hfill{$\Box$}

\medskip\noindent
{\bf Proof of Lemma~\ref{lem:G_NB(nu)}:}
According to, for example, Example~3 in \cite{aldirawi2022modeling}, $E \Psi(\nu + Y_1) = \Psi(\nu) - \log(p)$, if $Y_1 \sim {\rm NB}(\nu, p)$ (please note that the parameters in their paper are different from ours, with $p$ replaced with $1-p$). Then according to \eqref{eq:G_F(nu)_psi}, \begin{eqnarray*}
G_{NB(\nu, p)}(\nu) &=& E \Psi(\nu + Y_1) - \Psi(\nu) - \log \left( 1 + \frac{\mu}{\nu}\right)\\
&=&  \Psi(\nu) - \log(p) - \Psi(\nu) - \log \left(1+\frac{1-p}{p}\right)\\
&=& 0\ ,
\end{eqnarray*}
where $\mu = \nu(1-p)/p$ for NB$(\nu, p)$. In other words, 
\[
G_{NB(\nu, p)}(\nu) = \sum_{y=0}^\infty \frac{1-F_{NB(\nu, p)}(y)}{\nu+y}  - \log\left(1 + \frac{\mu}{\nu}\right)= 0,
\] 
where $F_{NB(\nu, p)}(y) = P(Y_1 \leq y)$ with $Y_1 \sim {\rm NB}(\nu, p)$.
\hfill{$\Box$}

\medskip\noindent
{\bf Proof of Lemma~\ref{lem:compare_f_p_f_NB_log}:}
Plugging in \eqref{eq:f_lambda_x} and \eqref{eq:f_NB_x}, we have for $x>0$,
\begin{eqnarray*}
r(x) &=& \log \Gamma(\nu+x) - x\log(\nu+\lambda) - \log\Gamma(\nu) + \lambda - \nu\log\left(1+\frac{\lambda}{\nu}\right)\ ,\\
r'(x) &=& \Psi(\nu+x) - \log(\nu+\lambda) \ .
\end{eqnarray*}
According to, for example, \cite{abramowitz1964handbook} or \cite{alzer2017harmonic}, $\Psi(x)$ is strictly increasing and differentiable on $(0, \infty)$. Therefore, $r'(x)$ is strictly increasing and differentiable on $(0, \infty)$. 

According to \cite{alzer1997some}, $\log x - \frac{1}{x} < \Psi(x) < \log x - \frac{1}{2x}$ for all $x>0$. Then for $0 < x < \lambda$, $r'(x) < \Psi(\nu + \lambda) - \log(\nu + \lambda) < 0$, and $\lim_{x\rightarrow \infty} r'(x) = \infty$. The conclusion about $r'(x)$ follows, with $x_* = \Psi^{-1}\left(\log(\nu+\lambda)\right) - \nu$ well defined and $x_* > \lambda$.

Since $x > \log(1+x)$ for all $x>0$, then $r(0) = \lambda - \nu\log\left(1+\frac{\lambda}{\nu}\right) > \lambda - \nu\cdot \frac{\lambda}{\nu} = 0$.

On the other hand, according to Lemma~1 of \cite{minc1964some} or \cite{alzer1997some}, 
\[
0 < \log\Gamma(x) - \left[(x-\frac{1}{2})\log x - x + \frac{1}{2}\log (2\pi)\right] < \frac{1}{x}
\]
for all $x>1$. It can be verified that $\lim_{x\rightarrow\infty} \frac{r(x)}{x\log(\nu+x)} = 1$, which implies $\lim_{x\rightarrow \infty} r(x) = \infty$.
\hfill{$\Box$}

\medskip\noindent
{\bf Proof of Theorem~\ref{thm:compare_f_p_f_NB}:}
According to Lemma~\ref{lem:compare_f_p_f_NB_log}, $r(x)$ attains its minimum at $x_*>\lambda$. 

We first show that $r(x_*) < 0$. Otherwise, if $r(x_*) \geq 0$, then $r(x) \geq 0$ for all $x\geq 0$ and $r(x) > 0$ for all $x\neq x_*$~. Then equivalently, we have $f_{NB(\nu, p)}(x) \geq f_\lambda(x)$ for all $x\geq 0$ and $f_{NB(\nu, p)}(x) > f_\lambda(x)$ for all $x\neq x_*$~, which implies $\sum_{y=0}^\infty f_{NB(\nu, p)}(y) > \sum_{y=0}^\infty f_\lambda(y)$. However, $\sum_{y=0}^\infty f_{NB(\nu, p)}(y) = \sum_{y=0}^\infty f_\lambda(y) = 1$. The contradiction implies that we must have $r(x_*) < 0$.

Note that the sign of $r(y)$ is exactly the same as the sign of $d(y)$ for all $y\in \{0, 1, 2, \ldots \}$. We let $K_1$ be the largest possible integer such that $d(y)\geq 0$ for all $0\leq y\leq K_1$, and $K_2$ be the smallest possible integer such that $d(y)\geq 0$ for all $y\geq K_2$~. Note that we must have $K_2-K_1\geq 2$ and $d(y)<0$ for each $K_1 < y < K_2$ (otherwise, we have the contradiction $\sum_{y=0}^\infty f_{NB(\nu, p)}(y) > \sum_{y=0}^\infty f_\lambda(y)$ again).

According to Lemma~\ref{lem:compare_f_p_f_NB_log}, $r(0)>0$ and thus $d(0)=D(0)>0$. Since $d(y)>0$ for $0\leq y<K_1$ and $y>K_2$, then $D(y)$ strictly increases on $0\leq y<K_1$ and $y>K_2$; while $d(y) < 0$ for $K_1 < y < K_2$ implies that $D(y)$ strictly decreases on $K_1 < y < K_2$~. Since $\lim_{y\rightarrow \infty} D(y) = 1-1 = 0$, then $D(y)$ must first increase and attain its maximum $D(K_1) > 0$, then strictly decrease and attain its minimum $D(K_2) < 0$, and then strictly increase after $K_2$ and approach to $0$ as $y$ goes to $\infty$. We let $K_*$ be the largest possible integer such that $D(y)\geq 0$. Then we must have $K_1 < K_* < K_2$~, $D(y) > 0$ for all $0\leq y<K_*$~, and $D(y) < 0$ for all $y > K_*$~.
\hfill{$\Box$}

\medskip\noindent
{\bf Proof of Lemma~\ref{lem:sum_D(y)}:}
First of all, for any random variable $Y$ which takes values in $\{0, 1, 2, \ldots \}$ with a CDF $F(y)$ and a finite expectation $E(Y)$, we must have $\sum_{y=0}^\infty [1-F(y)] = E(Y)$. Actually,
\begin{eqnarray*}
\sum_{y=0}^\infty [1-F(y)] &=& \sum_{y=0}^\infty \sum_{k=y+1}^\infty P(Y=k) = \sum_{k=1}^\infty \sum_{y=0}^{k-1} P(Y=k)\\
&=& \sum_{k=1}^\infty k\cdot P(Y=k)\\
&=& E(Y)\ .
\end{eqnarray*}
Since $p=\nu/(\nu+\lambda)$, the means of NB$(\nu, p)$ and Poisson$(\lambda)$ are the same. That is, 
\[
\sum_{y=0}^\infty [1-F_{NB(\nu, p)}(y)] = \sum_{y=0}^\infty [1-F_\lambda(y)] = \lambda\ .
\]
Then 
\[
\sum_{y=0}^\infty D(y) = \sum_{y=0}^\infty [F_{NB(\nu, p)}(y) - F_\lambda (y)] = \sum_{y=0}^\infty [1-F_\lambda(y)] - \sum_{y=0}^\infty [1-F_{NB(\nu, p)}(y)] = 0\ .
\]
\hfill{$\Box$}

\medskip\noindent
{\bf Proof of Theorem~\ref{thm:G_lambda_nu>0}:}
Given $\lambda > 0$ and $\nu>0$, we let $p = \nu/(\nu+\lambda) \in (0,1)$. Then the mean of NB$(\nu, p)$ is $\mu = \nu(1-p)/p = \nu\cdot \lambda/\nu = \lambda$. According to Lemmas~\ref{lem:G_F(nu)_1-F} and \ref{lem:G_NB(nu)}, $G_{NB(\nu, p)}(\nu) = \sum_{y=0}^\infty \frac{1-F_{NB(\nu, p)}(y)}{\nu+y} - \log \left(1+\frac{\lambda}{\nu}\right)=0$. Since $G_\lambda(\nu) = \sum_{y=0}^\infty \frac{1-F_\lambda (y)}{\nu+y} - \log \left(1+\frac{\lambda}{\nu}\right)$, we only need to show that $\sum_{y=0}^\infty \frac{1-F_\lambda (y)}{\nu+y} > \sum_{y=0}^\infty \frac{1-F_{NB(\nu, p)}(y)}{\nu+y}$. 

According to Theorem~\ref{thm:compare_f_p_f_NB} and Lemma~\ref{lem:sum_D(y)}, we have $D(y) \geq 0$ for $0\leq y\leq K_*$~, $D(y)<0$ for $y>K_*$~, and $\sum_{y=0}^{K_*} D(y) = -\sum_{y=K_*+1}^\infty D(y) > 0$. Then 
\begin{eqnarray*}
\sum_{y=0}^\infty \frac{1-F_\lambda (y)}{\nu+y} - \sum_{y=0}^\infty \frac{1-F_{NB(\nu, p)}(y)}{\nu+y}  
&=& \sum_{y=0}^\infty \frac{F_{NB(\nu, p)} (y) - F_\lambda (y)}{\nu+y} = \sum_{y=0}^\infty \frac{D(y)}{\nu+y}\\
&=& \sum_{y=0}^{K_*} \frac{D(y)}{\nu+y} + \sum_{y=K_*+1}^\infty \frac{D(y)}{\nu+y}\\
&>& \sum_{y=0}^{K_*} \frac{D(y)}{\nu+K_*} + \sum_{y=K_*+1}^\infty \frac{D(y)}{\nu+K_*+1}\\
&=& \sum_{y=0}^{K_*} D(y) \cdot \left(\frac{1}{\nu+K_*} - \frac{1}{\nu+K_*+1}\right)\\
&=& \sum_{y=0}^{K_*}D(y) \cdot \frac{1}{(\nu+K_*)(\nu+K_*+1)}\\
&>& 0 \ .
\end{eqnarray*}
\hfill{$\Box$}

\medskip\noindent
{\bf Proof of Lemma~\ref{lem:g_lambda_nu>0}:}
According to Theorem~\ref{thm:G(nu)}, if $Y_1, \ldots, Y_n$ are iid $\sim$ Poisson$(\lambda)$ with $\lambda > 0$, then there exists an event $A$ with probability $1$, such that, for each $\omega \in A$, $g(\nu)\rightarrow G_\lambda(\nu) > 0$, as $n$ goes to $\infty$. We let $c=G_\lambda(\nu)/2 > 0$. Then for each $\omega \in A$, there exists an $N(\omega)>0$, such that, $|g(\nu) - G_\lambda(\nu)| < c$ and thus $g(\nu) > G_\lambda(\nu) - c = c$, for all $n\geq N(\omega)$.
\hfill{$\Box$}

\medskip\noindent
{\bf Proof of Lemma~\ref{lem:g_nu_0_nu0}:}
First of all, for any $a\in \left[\frac{1}{2}(1-e^{-\lambda}), \frac{3}{2}(1-e^{-\lambda})\right]$, $b\in \left[\frac{\lambda}{2}, \frac{3\lambda}{2}\right]$,
\[
h(\nu,a,b) \geq \frac{1-e^{-\lambda}}{2\nu} - \log\left(1+\frac{3\lambda}{2\nu}\right)\ .
\]
Since $\lim_{\nu \rightarrow 0^+} \frac{1-e^{-\lambda}}{2\nu} = \infty$ and
\[
\lim_{\nu \rightarrow 0^+} \frac{\log\left(1 + \frac{3\lambda}{2\nu}\right)}{\frac{1-e^{-\lambda}}{2\nu}} 
= \lim_{\nu \rightarrow 0^+} \frac{\frac{3\lambda}{2}}{(1+\frac{3\lambda}{2\nu})\frac{1-e^{-\lambda}}{2}} = 0\ ,
\]
then $\lim_{\nu \rightarrow 0^+} \frac{1-e^{-\lambda}}{2\nu} - \log\left(1 + \frac{3\lambda}{2\nu}\right) = \infty$. Therefore, there exists a $\nu_0>0$, such that, $\frac{1-e^{-\lambda}}{2\nu} - \log\left(1 + \frac{3\lambda}{2\nu}\right) > c$ for all $0 < \nu \leq \nu_0$~, which leads to the conclusion.
\hfill{$\Box$}

\medskip\noindent
{\bf Proof of Lemma~\ref{lem: upperbound_gr_gnu}:}
According to \eqref{eq:g'(nu)}, for each $\nu > \nu_0$,
\begin{eqnarray*}
|g'(\nu)| &=& \left|-\frac{1}{n} \sum_{y \in I\setminus \{0\}} f_y \left[\frac{1}{\nu^2} + \cdots + \frac{1}{(\nu + y - 1)^2}\right]  + \frac{1}{\nu} - \frac{1}{\nu + \bar{Y}_n}\right|\\
&\leq & \frac{1}{n} \sum_{y \in I\setminus \{0\}} f_y \left[\frac{1}{\nu^2} + \cdots + \frac{1}{(\nu + y - 1)^2}\right] + \frac{2}{\nu_0}\\
&\leq & \frac{1}{n} \sum_{y \in I\setminus \{0\}} f_y \cdot \sum_{k=0}^\infty \frac{1}{(\nu+k)^2} + \frac{2}{\nu_0}\\
&\leq & 1\cdot \sum_{k=0}^\infty \frac{1}{(\nu_0 + k)^2} + \frac{2}{\nu_0}\\
&= & \frac{1}{{\nu_0}^2} + \sum_{k=1}^\infty \frac{1}{(\nu_0 + k)^2} + \frac{2}{\nu_0}\\
&\leq & \frac{1}{{\nu_0}^2} + \sum_{k=1}^\infty \frac{1}{k^2} + \frac{2}{\nu_0}\\
&= & \frac{1}{{\nu_0}^2} + \frac{\pi}{6} + \frac{2}{\nu_0}\ .
\end{eqnarray*}
Let $M = \frac{1}{{\nu_0}^2} + \frac{2}{\nu_0} + \frac{\pi}{6} > 0$. Then $|g'(\nu)| \leq M < \infty$ for all $\nu>\nu_0$~.
\hfill{$\Box$}

\medskip\noindent
{\bf Proof of Theorem~\ref{thm:nu_max}:} 
Given $\lambda > 0$ and $\nu_{\rm max}>0$, we set $c=1$ in Lemma~\ref{lem:g_nu_0_nu0} and obtain $\nu_0 \in (0, \nu_{\rm max})$ as in  Lemma~\ref{lem:g_nu_0_nu0} and event $A$ as defined in \eqref{eq:A_set}. Then for each $\omega \in A$, let $N_0(\omega)$ be a large enough integer, such that for all $n\geq N_0(\omega)$, $\frac{\#\{i:Y_i(\omega)>0\}}{n} \in \left[\frac{1}{2}(1-e^{-\lambda}), \frac{3}{2}(1-e^{-\lambda})\right]$, $\bar{Y}_n(\omega)\in \left[\frac{\lambda}{2}, \frac{3\lambda}{2}\right]$, and thus $g(\nu)(\omega)\geq 1$ for all $0 < \nu \leq \nu_0$~.

According to Lemma~\ref{lem:G_F(nu)_1-F}, $G_\lambda(\nu)$ is continuous on $(0, \infty)$. According to Theorem~\ref{thm:G_lambda_nu>0}, $G_\lambda(\nu)>0$ for each $\nu>0$. Therefore, for the given $0 < \nu_0 < \nu_{\rm max} < \infty$, there exists a $c \in (0, 1)$ such that $G_\lambda (\nu) > c$ for all $\nu \in [\nu_0, \nu_{max}]$. On the other hand, according to Lemma~\ref{lem: upperbound_gr_gnu}, there exists an $M>0$, which does not depend on $\omega\in \Omega$, such that $|g'(\nu)| \leq M < \infty$ for all $\nu>\nu_0$ and all $\omega \in \Omega$.

We let $K = 4(\nu_{\rm max} - \nu_0)M/c + 1$, then $(\nu_{\rm max} - \nu_0)M/K < c/4$. We denote the grid points, $\nu_i = \nu_0 + (\nu_{\rm max} - \nu_0)\cdot i/K$, for $i=1, \ldots, K$. According to Theorem~\ref{thm:G(nu)}, for each $i$, $g(\nu_i) \rightarrow G_\lambda(\nu_i)$ almost surely. We denote
\[
A_K = \{\omega \in A \mid  \lim_{n\rightarrow \infty} g(\nu_i)(\omega) = G(\nu_i), i=1, \ldots, K\}\ ,
\] 
which has probability $1$.
For each $\omega \in A_K$~, there exist $N_i(\omega)$, $i=1, \ldots, K$, such that, $g(\nu_i)(\omega) \geq \frac{c}{2}$ for all $n\geq N_i(\omega)$.  We denote $N(\omega)=\max\{N_0(\omega), N_1(\omega), \ldots, N_K(\omega)\}$.

For each $\nu\in (\nu_0, \nu_{\rm max}]$, there exists a $\nu_i$~, such that, $\nu \in (\nu_{i-1}, \nu_i]$. Then 
\[
|g(\nu)(\omega) - g(\nu_i)(\omega)| \leq M\cdot \frac{\nu_{\rm max}-\nu_0}{K} < M\cdot \frac{c}{4M} = \frac{c}{4}\ ,
\]
which implies $g(\nu)(\omega) > g(\nu_i)(\omega) - \frac{c}{4} \geq \frac{c}{2} - \frac{c}{4} =  \frac{c}{4}$~, for all $n\geq N(\omega)$.

After all, for each $\omega \in A_K$ and each $n\geq N(\omega)$, $g(\nu)(\omega) > \frac{c}{4} > 0$ for all $0<\nu\leq \nu_{\rm max}$~. According to Theorem~\ref{thm:min_g^2_positive}, under the constraint $\nu \in (0, \nu_{\rm max}]$, the MLE of $\nu$ for a negative binomial distribution is $\hat{\nu}(\omega)=\nu_{\rm max}$~. The MLE of $p$ is obtained via \eqref{eq:p_as_nu}.
\hfill{$\Box$}

\medskip\noindent
{\bf Proof of Theorem~\ref{thm:nu_max_n}:}
According to the proof of Theorem~\ref{thm:nu_max}, for each positive integer $k$, there exists an event $A_k$ with probability $1$, such that, for each $\omega \in A_k$~, there exist an integer $N(\omega, k)$ and a lower bound $c_k>0$, such that, $g(\nu)>c_k>0$ for all $n > N(\omega, k)$ and all $\nu\in (0, \nu_{\rm max}(k)]$.

We let $A_0 = \cap_{k=1}^\infty A_k$~. Then $P(A_0)=1$ as well. Since $\nu_{\rm max}(n) \rightarrow \infty$ as $n$ goes to infinity, then for any $M>0$, there exists an integer $K_M>0$, such that, $\nu_{\rm max}(k) > M$ for all $k\geq K_M$~. For each $\omega \in A_0$~, given any $M>0$, we let $N_M = \max_{1\leq k\leq N_M} N(\omega, k)$. Then $g(\nu)(\omega) > \min_{1\leq k\leq K_M} c_k > 0$, for all $\nu \in (0, M]$ and all $n > N_M(\omega)$, which implies that the MLE $\hat\nu_n(\omega) > M$. That is, $\hat\nu_n(\omega)\rightarrow \infty$, as $n$ goes to infinity.
\hfill{$\Box$}

\medskip\noindent
{\bf Proof of Theorem~\ref{thm:nu_mu}:}
If $\bar{Y}_n=0$, then $Y_1=\cdots = Y_n=0$ and the log-likelihood $l(\nu, \mu) = n\nu\log \frac{\nu}{\nu+\mu}\leq 0$, which is maximized only when $\mu=0$. That is, $\hat{\mu}=0$ while $\hat{\nu}$ can be any number in $(0, \nu_{\rm max}]$.

If $\bar{Y}_n > 0$, then $Y_i>0$ for at least one $i$. It can be verified that $\lim_{\mu\to 0^+}l(\nu, \mu) = -\infty$ for any fixed $\nu>0$. Thus we must have $\hat{\mu}>0$. Actually, it can be verified that 
\[
\frac{\partial l (\nu, \mu)}{\partial\mu} = 
\frac{n\nu}{\mu(\nu+\mu)}\cdot (\bar{Y}_n - \mu),
\]
which implies  $\hat{\mu}=\bar{Y}_n>0$. Combining the cases of $\bar{Y}_n=0$ and $\bar{Y}_n > 0$, we conclude that $\hat{\mu}=0$ if and only if $\bar{Y}_n=0$.

Recall the invariance property of MLE (see, e.g., \cite{zehna1966invariance} or \cite{pal1992invariance}), that is, given a function $\tau$ and an MLE $\hat{\boldsymbol{\theta}}$ of $\boldsymbol{\theta}$,  $\tau(\hat{\boldsymbol{\theta}})$ is always an MLE of $\tau(\boldsymbol{\theta})$. In this case, $\mu = \nu(1-p)/p$ builds up a one-to-one correspondence from $\{(\nu, p) \mid \nu \in (0, \nu_{\rm max}], p \in (0, 1)\}$ to $\{(\nu, \mu) \mid \nu \in (0, \nu_{\rm max}], \mu \in (0, \infty)\}$, which will pass along the uniqueness of MLE. According to the proof of Theorem~\ref{thm:nu_max}, there exists an event $A$ of probability $1$, such that, for each $\omega \in A$, there exists an $N(\omega)>0$, $\bar{Y}_n(\omega) \geq \lambda/2 >0$, $(\hat{\nu}=\nu_{\rm max},\ \hat{p} = \nu_{\rm max}/(\nu_{\rm max} + \bar{Y}_n))$ is the unique MLE of $(\nu, p)$ for all $n\geq N(\omega)$. Therefore,    $(\hat{\nu}=\nu_{\rm max},\ \hat{\mu} = \hat{\nu}(1-\hat{p})/\hat{p} = \bar{Y}_n)$ is the unique MLE of $(\nu, \mu)$ for all $n\geq N(\omega)$.
\hfill{$\Box$}

\medskip\noindent
{\bf Proof of Theorem~\ref{thm:nu_P}:}
If $\bar{Y}_n=0$, then $Y_1=\cdots = Y_n=0$ and the log-likelihood $l(\nu, P) = -n\nu\log (1+P) \leq 0$, which is maximized only when $P=0$. That is, $\hat{P}=0$ while $\hat{\nu}$ can be any number in $(0, \nu_{\rm max}]$.

If $\bar{Y}_n > 0$, it can be verified that $\lim_{P\to 0^+}l(\nu, P) = -\infty$ for any fixed $\nu>0$. Thus we must have $\hat{P}>0$. Actually, it can be verified that 
\[
\frac{\partial l (\nu, P)}{\partial P} = 
\frac{n\nu}{P(1+P)} \left(\frac{\bar{Y}_n}{\nu} - P\right),
\]
which implies  $\hat{P}=\bar{Y}_n/\hat{\nu} > 0$. Combining the cases of $\bar{Y}_n=0$ and $\bar{Y}_n > 0$, we conclude that $\hat{P}=0$ if and only if $\bar{Y}_n=0$.

Note that in this case, $P = (1-p)/p$ builds up a one-to-one correspondence from $\{(\nu, p) \mid \nu \in (0, \nu_{\rm max}], p \in (0, 1)\}$ to $\{(\nu, P) \mid \nu \in (0, \nu_{\rm max}], P \in (0, \infty)\}$. According to the proof of Theorem~\ref{thm:nu_max}, there exists an event $A$ of probability $1$, such that, for each $\omega \in A$, there exists an $N(\omega)>0$, $\bar{Y}_n(\omega) \geq \lambda/2 >0$, $(\hat{\nu}=\nu_{\rm max},\ \hat{p} = \nu_{\rm max}/(\nu_{\rm max} + \bar{Y}_n))$ is the unique MLE of $(\nu, p)$ for all $n\geq N(\omega)$. Therefore,    $(\hat{\nu}=\nu_{\rm max},\ \hat{P} = (1 - \hat{p})/\hat{p} = \bar{Y}_n/\nu_{\rm max})$ is the unique MLE of $(\nu, P)$ for all $n\geq N(\omega)$.
\hfill{$\Box$}

\medskip\noindent
{\bf Proof of Theorem~\ref{thm:nu_1-p}:}
If $\bar{Y}_n=0$, then $Y_1=\cdots = Y_n=0$ and the log-likelihood $l(\nu, P) = n\nu\log (1-P) \leq 0$, which is maximized only when $P=0$. That is, $\hat{P}=0$ while $\hat{\nu}$ can be any number in $(0, \nu_{\rm max}]$.

If $\bar{Y}_n > 0$, it can be verified that $\lim_{P\to 0^+}l(\nu, P) = -\infty$ for any fixed $\nu>0$. Thus we must have $\hat{P}>0$. Actually, it can be verified that 
\[
\frac{\partial l (\nu, P)}{\partial P} = 
\frac{n(\nu+\bar{Y}_n)}{P(1-P)} \left(\frac{\bar{Y}_n}{\nu+\bar{Y}_n} - P\right),
\]
which implies  $\hat{P}=\bar{Y}_n/(\hat{\nu} + \bar{Y}_n) > 0$. Combining the cases of $\bar{Y}_n=0$ and $\bar{Y}_n > 0$, we conclude that $\hat{P}=0$ if and only if $\bar{Y}_n=0$.

In this case, $P = 1-p$ builds up a one-to-one correspondence from $\{(\nu, p) \mid \nu \in (0, \nu_{\rm max}], p \in (0, 1)\}$ to $\{(\nu, P) \mid \nu \in (0, \nu_{\rm max}], P \in (0, 1)\}$. According to the proof of Theorem~\ref{thm:nu_max}, there exists an event $A$ of probability $1$, such that, for each $\omega \in A$, there exists an $N(\omega)>0$, $\bar{Y}_n(\omega) \geq \lambda/2 >0$, $(\hat{\nu}=\nu_{\rm max},\ \hat{p} = \nu_{\rm max}/(\nu_{\rm max} + \bar{Y}_n))$ is the unique MLE of $(\nu, p)$ for all $n\geq N(\omega)$. Therefore,    $(\hat{\nu}=\nu_{\rm max},\ \hat{P} = 1 - \hat{p} = \bar{Y}_n/(\nu_{\rm max} + \bar{Y}_n)$ is the unique MLE of $(\nu, P)$ for all $n\geq N(\omega)$.
\hfill{$\Box$}

\medskip\noindent
{\bf Proof of Lemma~\ref{lem:continuity_mu_p}:}
We only need to confirm the definitions in \eqref{eq:pmf_mu_0_p}, \eqref{eq:pmf_mu_p_1}, and \eqref{eq:pmf_mu_0_p_1} as the corresponding limits of $f_{\mu, p}(y)$ defined in \eqref{eq:pmf_mu_p} with more detail.

To verify \eqref{eq:pmf_mu_0_p}, we fix a $p\in (0,1)$. According to \eqref{eq:pmf_mu_p},  for $\mu>0$, $f_{\mu, p}(0) = p^{\frac{\mu p}{1-p}}$, which leads to $\lim_{\mu\to 0^+} f_{\mu, p}(0) = 1$. For each $y>0$, $\lim_{\mu\to 0^+} f_{\mu, p}(y) = 0$, as $\lim_{x\to 0^+} \Gamma(x) = \infty$. 

To verify \eqref{eq:pmf_mu_p_1}, we fix a $\mu \in (0, \infty)$. It can be verified that $\log f_{\mu, p}(0) = \frac{\mu p}{1-p} \log p = \frac{\mu p}{1-p} \log [1 - (1-p)]$ goes to $-\mu$, as $p\to 1^-$. That is, $\lim_{p\to 1^-} f_{\mu, p}(0) = e^{-\mu}$. For each $y>0$, it can be verified that, as $p\to 1^-$,
\[
f_{\mu, p}(y) = \frac{1}{y!} \cdot  \prod_{j=0}^{y-1}\left(\frac{\mu p}{1-p} + j\right) (1-p)^y \cdot p^{\frac{\mu p}{1-p}} \longrightarrow  \frac{1}{y!} \cdot \mu^y \cdot e^{-\mu} \ .
\]

As for \eqref{eq:pmf_mu_0_p_1}, when $\mu\to 0^+$ and $p\to 1^-$, it can be verified that $f_{\mu, p}(0) = p^{\frac{\mu p}{1-p}} = \exp\left\{\frac{\mu p}{1-p}\log[1-(1-p)]\right\}\to  1$.
For each $y>0$, it can be verified that $f_{\mu, p}(y) = \frac{1}{y!} \cdot  \prod_{j=0}^{y-1}[\mu p + j(1-p)] \cdot p^{\frac{\mu p}{1-p}} \to \frac{1}{y!} \cdot 0 \cdot 1 = 0$.
\hfill{$\Box$}

\medskip\noindent
{\bf Proof of Theorem~\ref{thm:mu_p_MLE_existence}:}
Due to Lemma~\ref{lem:continuity_mu_p}, we only need to show that the maximum of $l(\mu, p)$ {\it (1)} won't occur when $\mu\to\infty$; and {\it (2)} can always occur at some $p\in (0,1]$.

{\it Case (i):} $\bar{Y}_n=0$, which implies $Y_1=\cdots = Y_n=0$. Fixing any $p\in (0,1)$, $l(\mu, p) = \frac{n\mu p}{1-p} \log p < 0$ when $\mu>0$, and $l(0, p)\equiv 0$, which implies that the best $\mu$ is $0$, while $p$ can be any number in $(0,1)$. Fixing $p=1$, $l(\mu,1) = -n\mu < 0$ when $\mu>0$, and $l(0, 1) = 0$, which implies that the best $\mu$ is $0$ again. Combining the cases of $p\in (0,1)$ and $p=1$, we  conclude that $\bar{Y}_n=0$ implies that the best $\mu$ must be $0$, while $p$ can be any number in $(0,1]$.

{\it Case (ii):} $\bar{Y}_n>0$, which implies that $Y_i>0$ for at least one $i$. Fixing any $p\in (0,1)$, it can be verified that $l(0, p) = \lim_{\mu\to 0^+} l(\mu, p) = -\infty$ and $\lim_{\mu\to \infty} l(\mu, p) = -\infty$, which imply that the best $\mu$ exists in $(0, \infty)$ in this case. When $p=1$, $l(0,1) = \lim_{\mu\to 0^+} l(\mu, 1) = -\infty$ and $\lim_{\mu\to\infty} l(\mu, 1) = -\infty$, which imply that the best $\mu$ exists in $(0, \infty)$ as well. On the other hand, fixing any $\mu \in (0, \infty)$, $\lim_{p\to 0^+} l(\mu, p) = -\infty$ as $\lim_{x\to 0^+} \Gamma(x)=\infty$. Then the best $p$ exists in $(0,1]$. As a conclusion, $\bar{Y}_n>0$ implies that an MLE $(\hat{\mu}, \hat{p})$ must exist in $(0, \infty)\times (0, 1]$. If $\hat{p}=1$, it can be verified that $l(\mu, 1)$ is maximized at $\hat{\mu}=\bar{Y}_n$ only. If $\hat{p}<1$, then we must have $\partial l/\partial \mu = \partial l/\partial p = 0$ at $(\hat{\mu}, \hat{p})$. From \eqref{eq:l_mu(mu,p)} and \eqref{eq:l_p(mu,p)}, we obtain $\hat{\mu}=\bar{Y}_n$ as well. Combining the cases of $\hat{p}=1$ and $\hat{p}<1$, we always have $\hat{\mu}=\bar{Y}_n$~. 
\hfill{$\Box$}

\medskip\noindent
{\bf Proof of Theorem~\ref{thm:mu_p_poisson}:}
As a direct conclusion of Theorem~\ref{thm:mu_p_MLE_existence}, $\hat{\mu}=0$ if and only if $\bar{Y_n}=0$.

Suppose $Y_1, \ldots, Y_n$ are iid $\sim$ Poisson($\lambda$) with $\lambda >0$. Fixing a $c>0$, for any $M>0$, according to the proof of Theorem~\ref{thm:nu_max}, there exists an event $A_M$ of probability $1$, such that, for any $\omega \in A_M$, {\it (1)} there exists an $N_M(\omega)>0$; {\it (2)} $\bar{Y}_n(\omega) \in \left[\frac{\lambda}{2},\ \frac{3}{2}\lambda\right]$ for all $n\geq N_M(\omega)$; and {\it (3)} $g(\nu)(\omega) > c/4 > 0$ for all $0 < \nu \leq M$ and all $n\geq N_M(\omega)$. 

According to Theorem~\ref{thm:mu_p_MLE_existence}, we always have $\hat{\mu} = \bar{Y}_n$~. For all $n\geq N_M(\omega)$, $\hat{\mu}(\omega) = \bar{Y}_n(\omega) \geq \lambda/2 > 0$, which implies 
\[
\frac{\partial l}{\partial p}\left(\bar{Y}_n(\omega), p\right) = \frac{n\bar{Y}_n(\omega)}{(1-p)^2} \cdot g\left(\frac{\bar{Y}_n(\omega) p}{1-p}\right) > 0\ ,
\]
for all $0<p\leq \frac{M}{M+\bar{Y}_n(\omega)}$~. We then conclude $\hat{p}(\omega) \geq \frac{M}{M+\bar{Y}_n(\omega)} \geq \frac{M}{M+\frac{3}{2}\lambda}$ for all $n\geq N_M(\omega)$.

We let $A_\infty = \cap_{M=1}^\infty A_M$, which also has probability $1$. Then for each $\omega \in A_\infty$, $\hat{p}\to 1$, as $n\to\infty$.
\hfill{$\Box$}

\section{A simulation study on using KS test for NB against Poisson distributions}\label{sec:supp_plots}

In this section, we use the Kolmogorov-Smirnov (KS) test \citep{massey1951kolmogorov}, a classical goodness-of-fit procedure, to illustrate that it is difficult to separate a Poisson law from a nearly equidispersed NB distribution. It supports the relevant statement in Section~\ref{sec:intro}.

The original KS test was proposed to assess whether a given random sample $\{Y_1, \ldots,$ $ Y_n\}$ originates from a continuous cumulative distribution function (CDF) $F_{\boldsymbol\theta}(y)$, characterized by specific model parameter(s) $\boldsymbol\theta$. The KS test statistic is defined as $D_n = \sup_{y} |F_n(y) - F_{\boldsymbol\theta}(y)|$, where $F_n(y) = n^{-1} \sum_{i=1}^n {\mathbf 1}_{(-\infty, y]}(Y_i)$ represents the empirical distribution function (EDF) of the sample. \cite{dimitrova2020computing} extended the applicability of the KS test to encompass a broader range of distributions $F_{\boldsymbol\theta}(y)$, including discrete and mixed distributions with known $\boldsymbol\theta$. When $\boldsymbol\theta$ is unknown, \cite{lilliefors1967kolmogorov,lilliefors1969kolmogorov} suggested replacing $\boldsymbol\theta$ with its estimate $\hat{\boldsymbol\theta}$ when calculating $D_n$~. Adjusted $p$-values or thresholds for test statistics can be estimated by parametric bootstrapping \citep{henze1996empirical}, nonparametric bootstrapping \citep{aldirawi2019identifying}, or nested bootstrapping \citep{dousti2022r} under those circumstances. R packages, such as iZID \citep{wang2020identifying} and AZIAD \citep{dousti2022r}, have been developed for such purposes.

For illustrative purposes, we adopt the parametric bootstrapping proposed by \cite{henze1996empirical} to calculate the critical numbers or thresholds of the KS test. More specifically, for a given sample size $n$, we generate $1{,}000$ independent random samples from Poisson($\lambda = 10$). For each sample we calculate the maximum likelihood estimate (MLE) $\hat{\boldsymbol\theta} = (\hat\nu, \hat{p})$, the KS test statistic $D_n = \sup_y|F_n(y) - F_{\hat{\boldsymbol{\theta}}}(y)|$, the critical number $d_n$ such that $P(D_n \geq d_n)=0.05$ based on the parametric bootstrapping \citep{henze1996empirical}. 
Since, in numerical computation, the MLE of $\nu$ may either diverge or converge when the data are Poisson, we impose a large upper bound, $\nu_{\max}=10{,}000$, to stabilize the estimation (see Section~\ref{sec:mle_NB} for more discussions). Table~\ref{table:Poisson-NB} reports the median and average of these quantities across the $1{,}000$ replications  (see below for the corresponding histograms and scatter plots as well), together with the empirical power, i.e., the ratio of rejecting $H_0$ (NB) given Poisson samples.
There is a clear trend: {\it (i)} both $D_n$ and $d_n$ vanish as $n$ goes to infinity; and {\it (ii)} the empirical power is poor (around $0.05$) and not getting better as $n$ goes to infinity. We repeat this simulation study with $\nu_{\rm max}=50,000$ (results not shown here), and our conclusions are similar.

\begin{table}[hbt]
\caption{KS tests for NB ($H_0$) against Poisson ($H_1$) over 1,000 simulations given $Y_1, \ldots, Y_n$ iid $\sim$ Poisson($\lambda = 10$)}
\label{table:Poisson-NB}
\begin{center}
\begin{tabular}{@{}lrrrr@{}}
\toprule
Sample Size $n$ & $50$   & $500$  & $5000$  & $10000$    \\ 
\midrule
$\hat{\nu}$ (median) & 10000 & 10000 & 10000 & 10000 \\
$\hat{\nu}$ (average)& 6002 & 5638 & 6094 & 6187 \\
$\hat{p}$ (median) & 0.9990 & 0.9990 & 0.9990 & 0.9990 \\
$\hat{p}$ (average) & 0.9384 & 0.9767 & 0.9921 & 0.9939 \\
Test Statistic $D_n$ (median) & 0.0630 & 0.0201 & 0.0063 & 0.0046 \\
Test Statistic $D_n$ (average) & 0.0654 & 0.0209 & 0.0066 & 0.0047 \\
Critical Number $d_n$ (median) & 0.1031 & 0.0327 & 0.0104 & 0.0073 \\
Critical Number $d_n$ (average) & 0.1031 & 0.0328 & 0.0104 & 0.0073 \\
Empirical Power at Poisson($10$) & 0.056 & 0.055 & 0.067 & 0.053 \\
\bottomrule
\end{tabular}	
\end{center}
\end{table} 

To better display the sampling distributions of the key quantities listed in Table~\ref{table:Poisson-NB}, we provide in this section histograms and scatter plots of them as well. Recall that, for each sample size $n \in \{50, 500, 5,000, 10,000\}$, we generate 1,000 independent Poisson($10$) samples; and for each sample, we compute the corresponding MLE $(\hat{\nu}, \hat{p})$, the KS test statistic $D_n$~, and the critical number or threshold $d_n$~. 

\begin{figure}[htbp]
    \centering
    \includegraphics[width=0.7\textwidth]{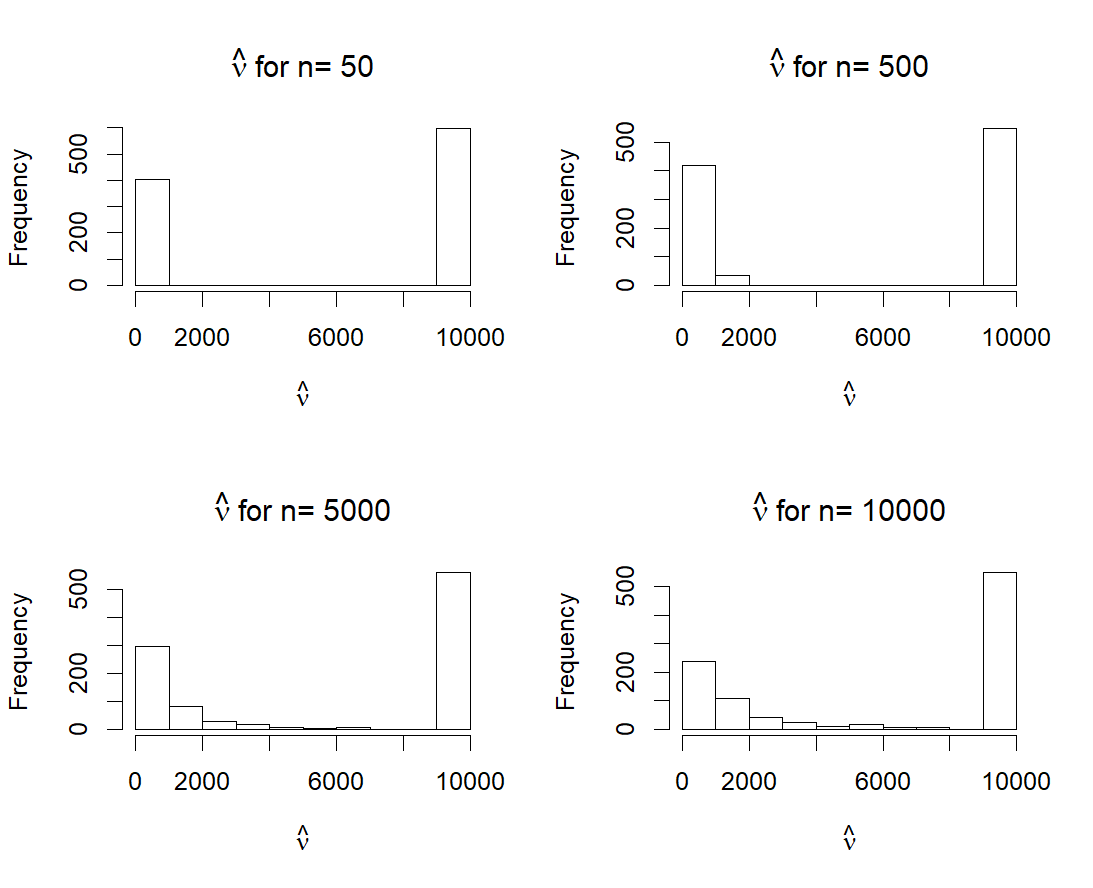} \\[1ex]
    \includegraphics[width=0.7\textwidth]{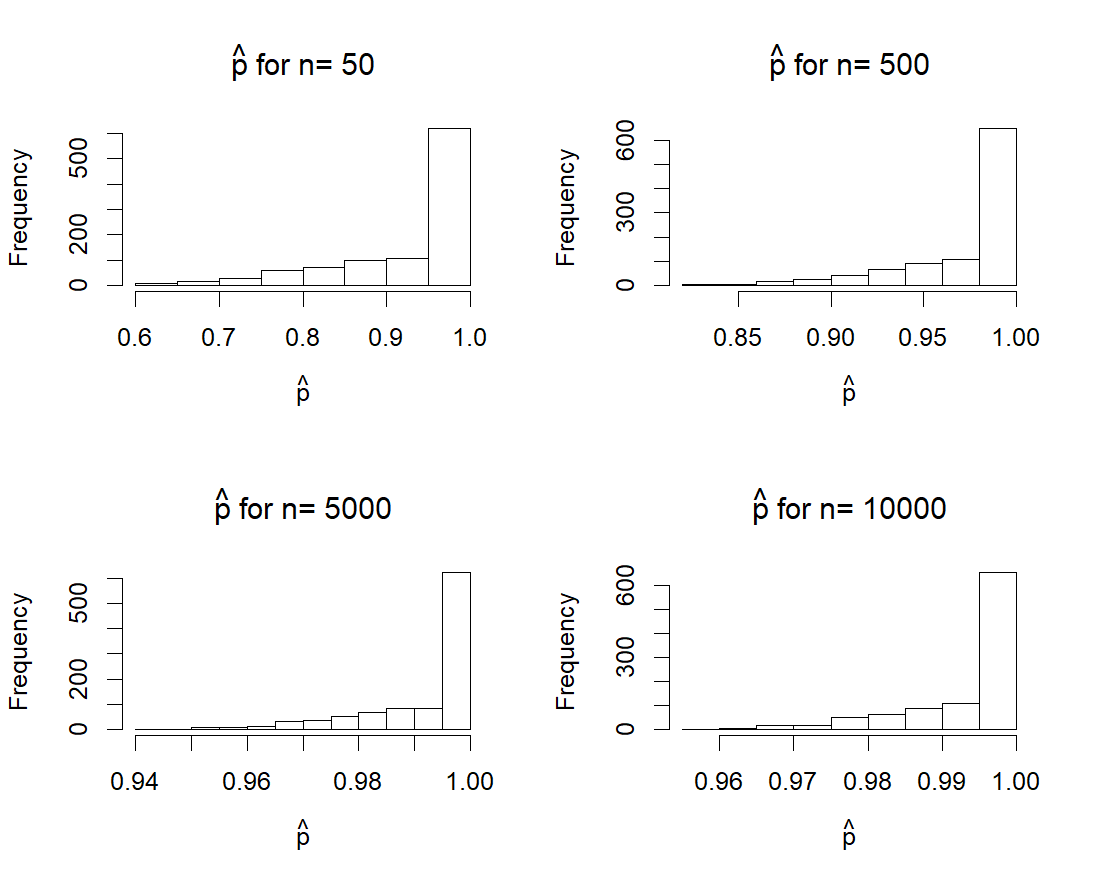} \\[1ex]
    \caption{Histograms of the MLEs $\hat{\nu}$ (top) and $\hat{p}$ (bottom) across 1,000 Poisson$(10)$ random samples with various sample sizes.}
    \label{fig:hist_nuhat_phat}
\end{figure}

\begin{figure}[htbp]
    \centering
    \includegraphics[width=0.7\textwidth]{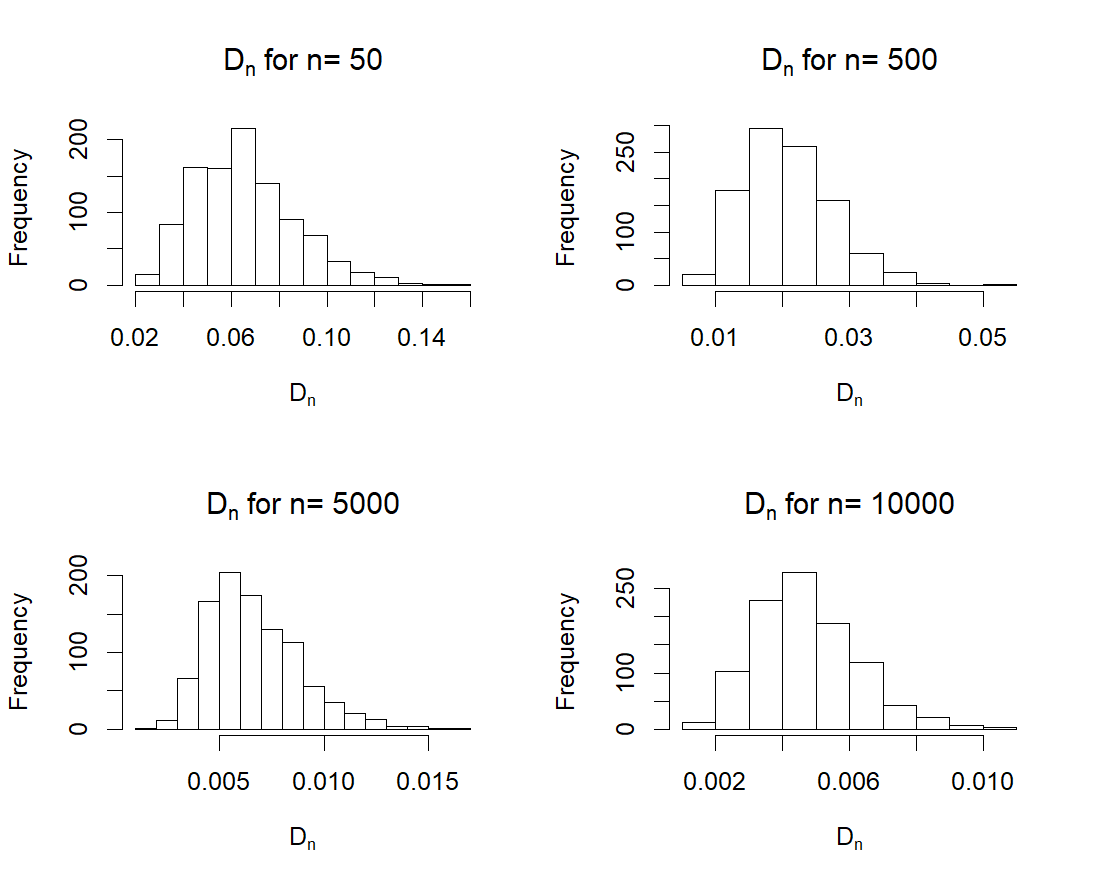} \\[1ex]
    \includegraphics[width=0.7\textwidth]{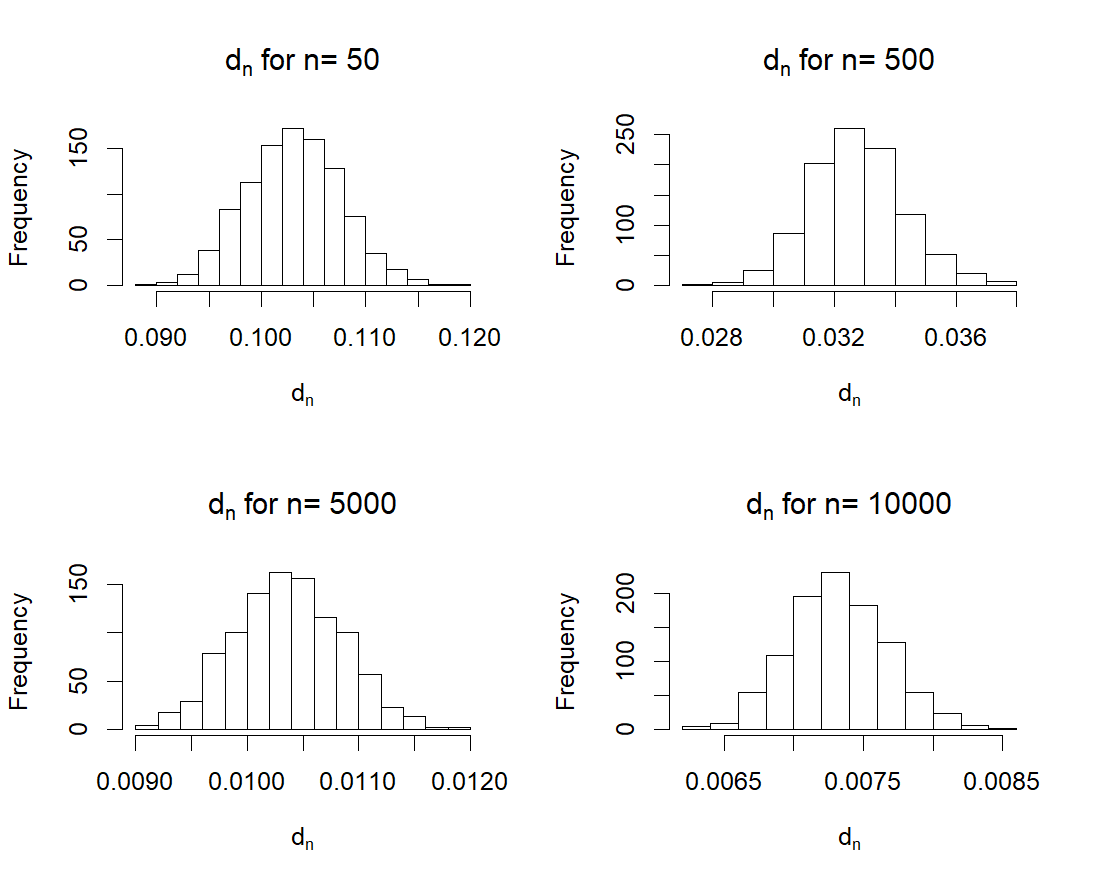} \\[1ex]
    \caption{Histograms of the KS test statistic $D_n$ (top) and the critical number $d_n$ (bottom) across 1,000 Poisson($10$) random samples with various sample sizes.}
    \label{fig:hist_Dn_dn}
\end{figure}

\begin{figure}[htbp]
    \centering
    \includegraphics[width=0.7\textwidth]{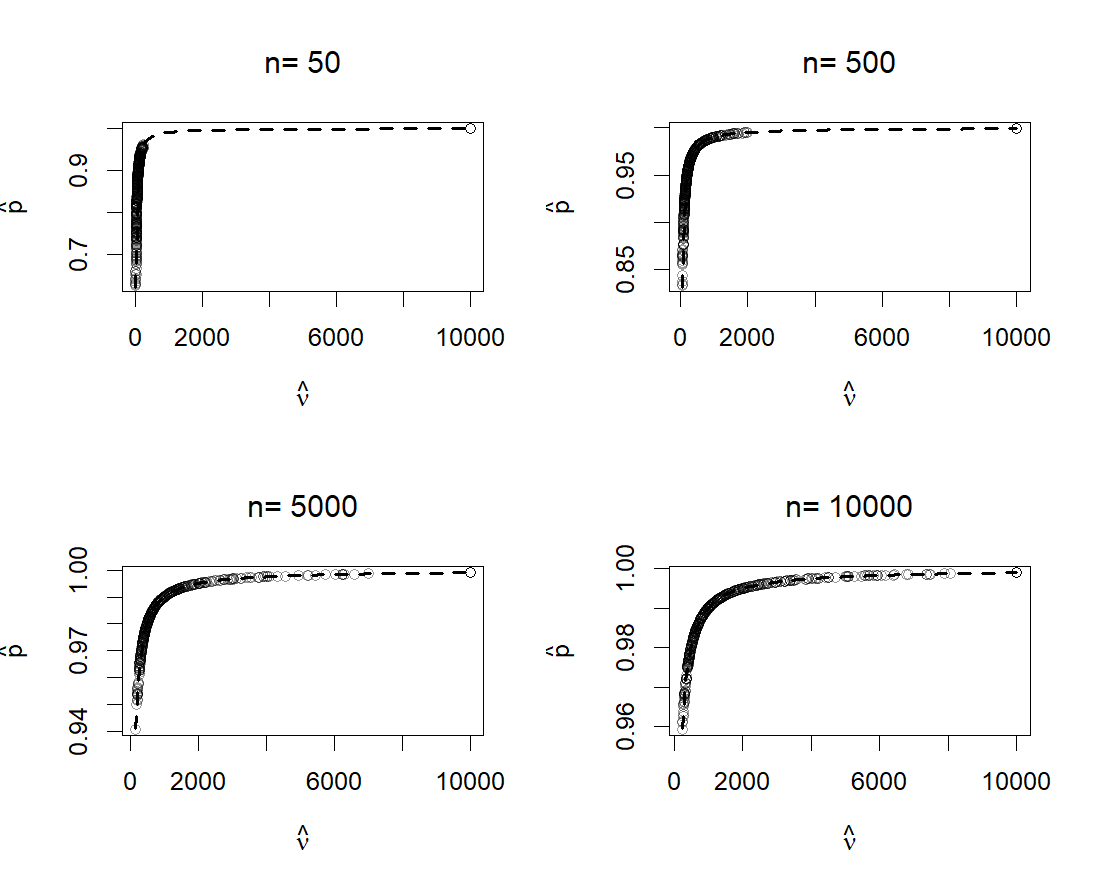} \\[1ex]
    \includegraphics[width=0.7\textwidth]{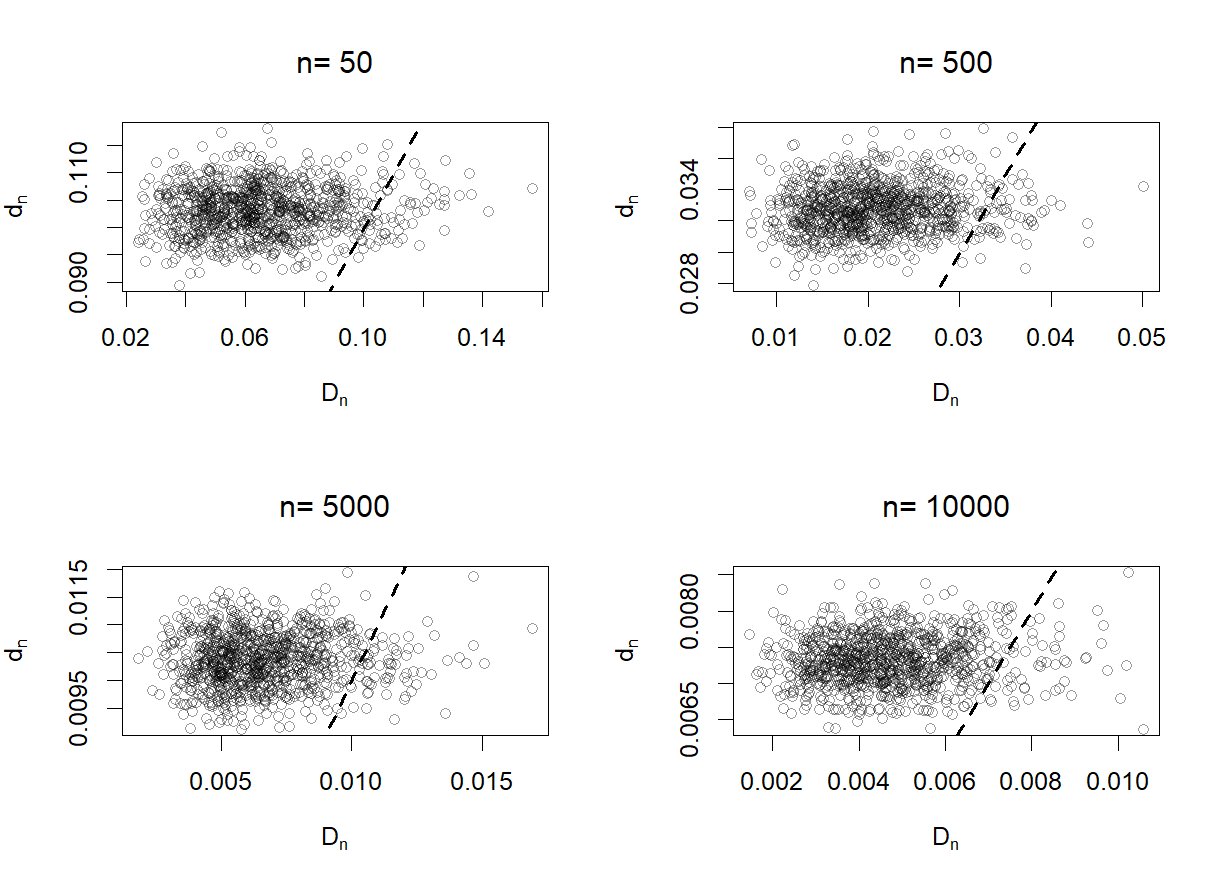}
    \caption{Scatter plots of $(\hat{\nu}, \hat{p})$ (top) along with its theoretical curve $p = \nu/(\nu+10)$, and scatter plots of $(D_n, d_n)$ (bottom) along with a reference line $d_n=D_n$~, based on 1,000 Poisson($10$) random samples.}
    \label{fig:scatter}
\end{figure}

Figure~\ref{fig:hist_nuhat_phat} shows the histograms of $\hat{\nu}$ and $\hat{p}$, 
revealing their concentrations near the boundary values. Figure~\ref{fig:hist_Dn_dn} displays the histograms of $D_n$ and $d_n$, illustrating their convergence towards zero when the sample size increases. Figure~\ref{fig:scatter} presents scatter plots of $(\hat{\nu},\hat{p})$ compared to their theoretical curve $p=\nu/(\nu+10)$, and scatter plots of $(D_n,d_n)$ with a reference line $d_n = D_n$~. These plots together provide a more detailed view of the variability and limiting behaviors of these estimators, test statistics, or thresholds across independent samples.

\section{Asymptotic properties of KS test statistic}\label{sec:supp_thms}

In this section, with the aid of Theorem~\ref{thm:nu_max_n}, we theoretically derive the asymptotic properties of the KS test statistic $D_n$ for testing NB$(\nu, p)$ ($H_0$) against Poisson distribution ($H_1$), given that the data come from a Poisson distribution. Our conclusion (see Theorem~\ref{thm: KStest}) confirms the limiting behavior of $D_n$ displayed in Table~\ref{table:Poisson-NB} of Section~\ref{sec:supp_plots}. We also obtain the limiting behavior of $D_n$ for testing extended NB$(\mu, p)$ distributions (see Theorem~\ref{thm: KStest_extended_NB}) proposed in Section~\ref{sec:new_parameterization}.

Given $Y_1, \ldots, Y_n$ that are iid $\sim {\rm Poisson}(\lambda)$, we denote  $\hat{\boldsymbol{\theta}} = (\hat{\nu}_n, \hat{p}_n)^T$ as the MLE of $\boldsymbol{\theta} = (\nu, p)^T$ for a negative binomial distribution under the constraint $0<\nu\leq \nu_{\rm max}(n)$, where $\lim_{n\rightarrow \infty}\nu_{\rm max}(n) = \infty$. 
With the aid of Theorem~\ref{thm:nu_max_n}, we obtain the following theorem for the limiting behavior of the KS test statistic $D_n$ for negative binomial distributions. Recall that for this KS test, the null hypothesis $H_0$ is that $Y_1, \ldots, Y_n$ are iid from an NB distribution with unknown parameters, and the alternative hypothesis $H_1$ is that  $Y_1, \ldots, Y_n$ are iid from a Poisson distribution with unknown parameters.

\begin{theorem}\label{thm: KStest}
Given $Y_1, \ldots, Y_n$ iid $\sim {\rm Poisson}(\lambda)$ with $\lambda > 0$, we denote $D_n = \sup_y|F_n(y) - F_{\hat{\boldsymbol{\theta}}}(y)|$ as the KS test statistic for testing a negative binomial distribution with unknown parameters, where $F_n(y)$ is the EDF of the data, and $F_{\hat{\boldsymbol{\theta}}}(y)$ is the CDF of NB$(\hat{\boldsymbol{\theta}})$ with the MLE $\hat{\boldsymbol{\theta}} = (\hat\nu, \hat{p})$ obtained under the constraint $0<\nu\leq \nu_{\rm max}(n)$, where $\lim_{n\rightarrow \infty}\nu_{\rm max}(n) = \infty$. Then $D_n \rightarrow 0$ almost surely, as $n \rightarrow \infty$.
\end{theorem}

\medskip\noindent
{\bf Proof of Theorem~\ref{thm: KStest}:}
We denote $F_\lambda(y)$ as the CDF of Poisson($\lambda$). For any given $\epsilon>0$, there exists an integer $M>0$, such that, $1-F_\lambda(y) \leq \frac{\epsilon}{6}$ for all $y \geq M$. Then
\begin{eqnarray*}
D_n
&=& \sup_y|F_n(y) - F_{\hat{\boldsymbol{\theta}}}(y)|\\
&=& \sup_y\left|F_n(y) - F_{\lambda}(y) + F_{\lambda}(y) - F_{\hat{\boldsymbol{\theta}}}(y)\right|\\
&\leq & \sup_y\left|F_n(y) - F_{\lambda}(y)\right| + \sup_{y<M}\left|F_{\lambda}(y) - F_{\hat{\boldsymbol{\theta}}}(y)\right|
 + \sup_{y \geq M}\left|F_{\lambda}(y) - F_{\hat{\boldsymbol{\theta}}}(y)\right|\\
&= & \sup_y\left|F_n(y) - F_{\lambda}(y)\right| + \sup_{y<M}\left|F_{\lambda}(y) - F_{\hat{\boldsymbol{\theta}}}(y)\right|\\
& &  +\ \sup_{y \geq M}\left|\left[1-F_{\lambda}(y)\right] - \left[1-F_{\hat{\boldsymbol{\theta}}}(y)\right]\right|\\
&\leq & \sup_y\left|F_n(y) - F_{\lambda}(y)\right| + \sup_{y<M}\left|F_{\lambda}(y) - F_{\hat{\boldsymbol{\theta}}}(y)\right|\\
& &  +\ \sup_{y \geq M}\left[1-F_{\lambda}(y)\right] + \sup_{y \geq M} \left[1-F_{\hat{\boldsymbol{\theta}}}(y)\right]\\
&\leq & \sup_y\left|F_n(y) - F_{\lambda}(y)\right| + \sup_{y<M}\left|F_{\lambda}(y) - F_{\hat{\boldsymbol{\theta}}}(y)\right| + \left[1-F_{\lambda}(M)\right] + \left[1-F_{\hat{\boldsymbol{\theta}}}(M)\right]\ .
\end{eqnarray*}

According to the Glivenko-Cantelli theorem (see, e.g., Theorem~7.5.2 in \cite{resnick1999probpath}), 
$\sup_y|F_n(y) - F_{\lambda}(y)| \rightarrow 0$
almost surely, as $n \rightarrow \infty$. 

Let $A=\{\omega\in \Omega \mid \lim_{n\rightarrow \infty} \hat\nu_n(\omega) = \infty \mbox{ and } \lim_{n\rightarrow\infty} \sup_y|F_n(y)(\omega) - F_{\lambda}(y)| = 0\}$. Then according to Theorem~\ref{thm:nu_max_n} and the Glivenko-Cantelli theorem, $P(A)=1$.

For each $\omega \in A$ and the given $\epsilon>0$, there exists a constant $N_e(\omega)>0$, such that, $\sup_y|F_n(y)(\omega) - F_{\lambda}(y)| \leq \frac{\epsilon}{3}$ for all $n \geq N_e(\omega)$. 

According to Theorem~5.13.1 in \cite{fisz1963probability}, $f_{{\rm NB}(\nu, p(\nu))}(k) \rightarrow f_{\lambda}(k)$ as $\nu \rightarrow \infty$ for each $k$, where $p(\nu)=\frac{\nu}{\nu+\lambda}$. Then, $F_{{\rm NB}(\nu, p(\nu))}(k) \rightarrow F_{\lambda}(k)$ as $\nu \rightarrow \infty$, for each $k$. For each $k=0, 1, \ldots, M$, there exists a $V_k>0$ such that $|F_{{\rm NB}(\nu, p(\nu))}(k) - F_{\lambda}(k)| < \frac{\epsilon}{6}$ for all $\nu\geq V_k$. Let $V_{\rm max} = \max\{V_0, V_1, \ldots, V_{M}\}$. Then for all $\nu\geq V_{\rm max}$, $\sup_{k < M}|F_{{\rm NB}(\nu, p(\nu))}(k) - F_{\lambda}(k)| \leq \frac{\epsilon}{6}$ and $|F_{{\rm NB}(\nu, p(\nu))}(M) - F_{\lambda}(M)| \leq \frac{\epsilon}{6}$.

In addition, for each $k$,
\[
\lim_{n\rightarrow \infty} F_{\hat{\boldsymbol{\theta}}(\omega)}(k) = \lim_{n\rightarrow \infty} F_{{\rm NB}(\hat\nu_n(\omega), \hat{p}_n(\hat\nu_n(\omega)))}(k) = \lim_{\nu \rightarrow \infty} F_{{\rm NB}(\nu, p(\nu))}(k) = F_\lambda(k) < \infty\ .
\]
There exists an $N_k(\omega)$ such that for all $n\geq N_k(\omega)$, $\hat\nu_n(\omega)>V_{\rm max}$, $|F_{\hat{\boldsymbol{\theta}}(\omega)}(k) - F_{\lambda}(k)| \leq \frac{\epsilon}{6}$. Then for each $n\geq N_{\rm max}(\omega) = \max\{N_0(\omega), N_1(\omega), \ldots, N_M(\omega), N_e(\omega)\}$, $\sup_{k<M}|F_{\hat{\boldsymbol{\theta}}(\omega)}(k) - F_{\lambda}(k)| \leq \frac{\epsilon}{6}$ and $0 < 1-F_{\hat{\boldsymbol{\theta}}(\omega)}(M) \leq 1-F_\lambda(M) + |F_{\hat{\boldsymbol{\theta}}(\omega)}(M) - F_{\lambda}(M)| \leq \frac{\epsilon}{6} + \frac{\epsilon}{6} = \frac{\epsilon}{3}$, which leads to $D_n(\omega) \leq \epsilon$.

Therefore, $D_n(\omega) $ goes to zero as $n$ goes to infinity.
\hfill{$\Box$}

\medskip

Theorem~\ref{thm: KStest} theoretically justifies the vanishing phenomenon on $D_n$ observed in Table~\ref{table:Poisson-NB}.

As a direct conclusion of Theorem~\ref{thm: KStest} and the Glivenko-Cantelli theorem (see, e.g., Theorem~7.5.2 in \cite{resnick1999probpath}), we have the following corollary to conclude the convergence of the estimated NB distribution to the true Poisson($\lambda$) distribution.

\begin{corollary}\label{cor:NB_to_Poisson_nu_p}
Under the circumstance of Theorem~\ref{thm: KStest}, we must have $\sup_y|F_{\hat{\boldsymbol{\theta}}}(y) - F_\lambda(y)| \rightarrow 0$ almost surely, as $n \rightarrow \infty$.    
\end{corollary}

\medskip\noindent
{\bf Proof of Corollary~\ref{cor:NB_to_Poisson_nu_p}:}
Under the circumstance of Theorem~\ref{thm: KStest}, we have $\sup_y|F_n(y) - F_{\hat{\boldsymbol{\theta}}}(y)| \rightarrow 0$ almost surely, as $n \rightarrow \infty$. According to the Glivenko-Cantelli theorem (see, e.g., Theorem~7.5.2 in \cite{resnick1999probpath}), 
$\sup_y|F_n(y) - F_{\lambda}(y)| \rightarrow 0$
almost surely, as $n \rightarrow \infty$. We conclude because $|F_{\hat{\boldsymbol{\theta}}}(y) - F_\lambda(y)| \leq |F_n(y) - F_{\hat{\boldsymbol{\theta}}}(y)| + |F_n(y) - F_{\lambda}(y)|$.
\hfill{$\Box$}

\medskip

Similarly to Theorem~\ref{thm: KStest}
and Corollary~\ref{cor:NB_to_Poisson_nu_p}, we have the following theorem for the extended negative binomial distribution described via \eqref{eq:pmf_mu_p}, \eqref{eq:pmf_mu_0_p}, \eqref{eq:pmf_mu_p_1}, and \eqref{eq:pmf_mu_0_p_1} with parameters $\mu\in [0, \infty)$ and $p\in (0, 1]$, as in Section~\ref{sec:new_parameterization}.

\begin{theorem}\label{thm: KStest_extended_NB}
Given $Y_1, \ldots, Y_n$ iid $\sim {\rm Poisson}(\lambda)$ with $\lambda > 0$, we denote $D_n = \sup_y|F_n(y) - F_{\hat\mu, \hat{p}}(y)|$ as the KS test statistic for testing an extended negative binomial distribution with unknown parameters, where $F_{\hat\mu, \hat{p}}$ is the CDF of an extended NB distribution with parameters equal to their MLEs $\hat\mu, \hat{p}$. Then (i) $D_n \rightarrow 0$ almost surely, as $n \rightarrow \infty$; (ii) $\sup_y|F_{\hat\mu, \hat{p}}(y) - F_\lambda(y)| \rightarrow 0$ almost surely, as $n \rightarrow \infty$.
\end{theorem}

\medskip\noindent
{\bf Proof of Theorem~\ref{thm: KStest_extended_NB}:}
Similarly to the proof of Theorem~\ref{thm: KStest}, we denote $F_\lambda(y)$ as the CDF of Poisson($\lambda$). For any given $\epsilon>0$, there exists an integer $M>0$, such that, $1-F_\lambda(y) \leq \frac{\epsilon}{6}$ for all $y \geq M$. Then
\begin{eqnarray*}
D_n
&=& \sup_y|F_n(y) - F_{\hat\mu, \hat{p}}(y)|\\
&=& \sup_y\left|F_n(y) - F_{\lambda}(y) + F_{\lambda}(y) - F_{\hat\mu, \hat{p}}(y)\right|\\
&\leq & \sup_y\left|F_n(y) - F_{\lambda}(y)\right| + \sup_{y<M}\left|F_{\lambda}(y) - F_{\hat\mu, \hat{p}}(y)\right|
 + \sup_{y \geq M}\left|F_{\lambda}(y) - F_{\hat\mu, \hat{p}}(y)\right|\\
&= & \sup_y\left|F_n(y) - F_{\lambda}(y)\right| + \sup_{y<M}\left|F_{\lambda}(y) - F_{\hat\mu, \hat{p}}(y)\right|\\
& &  +\ \sup_{y \geq M}\left|\left[1-F_{\lambda}(y)\right] - \left[1-F_{\hat\mu, \hat{p}}(y)\right]\right|\\
&\leq & \sup_y\left|F_n(y) - F_{\lambda}(y)\right| + \sup_{y<M}\left|F_{\lambda}(y) - F_{\hat\mu, \hat{p}}(y)\right|\\
& &  +\ \sup_{y \geq M}\left[1-F_{\lambda}(y)\right] + \sup_{y \geq M} \left[1-F_{\hat\mu, \hat{p}}(y)\right]\\
&\leq & \sup_y\left|F_n(y) - F_{\lambda}(y)\right| + \sup_{y<M}\left|F_{\lambda}(y) - F_{\hat\mu, \hat{p}}(y)\right| + \left[1-F_{\lambda}(M)\right] + \left[1-F_{\hat\mu, \hat{p}}(M)\right]\ .
\end{eqnarray*}

According to the Glivenko-Cantelli theorem (see, e.g., Theorem~7.5.2 in \cite{resnick1999probpath}), 
$\sup_y|F_n(y) - F_{\lambda}(y)| \rightarrow 0$
almost surely, as $n \rightarrow \infty$. 

Let $A=\{\omega\in \Omega \mid \hat\mu_n(\omega) = \bar{Y}_n(\omega), \lim_{n\rightarrow \infty} \bar{Y}_n(\omega) = \lambda, \lim_{n\rightarrow \infty} \hat{p}_n(\omega) = 1,$ $ \lim_{n\rightarrow\infty} \sup_y$ $|F_n(y)(\omega) - F_{\lambda}(y)| = 0\}$. Then according to Theorem~\ref{thm:mu_p_poisson}, Kolmogorov's Strong Law of Large Numbers, and the Glivenko-Cantelli theorem, $P(A)=1$.

For each $\omega \in A$ and the given $\epsilon>0$, there exists a constant $N_e(\omega)>0$, such that, $\sup_y|F_n(y)(\omega) - F_{\lambda}(y)| \leq \frac{\epsilon}{3}$ for all $n \geq N_e(\omega)$. 

According to Theorem~5.13.1 in \cite{fisz1963probability}, $f_{{\rm NB}(\nu, p(\nu))}(k) \rightarrow f_{\lambda}(k)$ as $\nu \rightarrow \infty$ for each $k$, where $p(\nu)=\frac{\nu}{\nu+\lambda}$. Then, $F_{{\rm NB}(\nu, p(\nu))}(k) \rightarrow F_{\lambda}(k)$ as $\nu \rightarrow \infty$, for each $k$. For each $k=0, 1, \ldots, M$, there exists a $V_k>0$ such that $|F_{{\rm NB}(\nu, p(\nu))}(k) - F_{\lambda}(k)| < \frac{\epsilon}{6}$ for all $\nu\geq V_k$~. Let $V_{\rm max} = \max\{V_0, V_1, \ldots, V_{M}\}$. Then for all $\nu\geq V_{\rm max}$, $\sup_{k < M}|F_{{\rm NB}(\nu, p(\nu))}(k) - F_{\lambda}(k)| \leq \frac{\epsilon}{6}$ and $|F_{{\rm NB}(\nu, p(\nu))}(M) - F_{\lambda}(M)| \leq \frac{\epsilon}{6}$.

In addition, for each $k\in \{0, 1, \ldots, M\}$,
\begin{eqnarray*}
\left|F_{\hat\mu_n(\omega), \hat{p}_n(\omega)}(k) - F_\lambda(k)\right| &=&
\left\{\begin{array}{cl}
\left|F_{{\rm NB}(\hat\nu_n(\omega), \hat{p}_n(\omega))}(k) - F_\lambda(k)\right|\ , & \mbox{ if }\hat{p}_n(\omega)<1;\\
\left|F_{\hat\lambda = \bar{Y}_n(\omega)}(k) - F_\lambda(k)\right|\ , & \mbox{ if }\hat{p}_n(\omega)=1.
\end{array}\right.\\
&\longrightarrow & 0\ , 
\end{eqnarray*}
as $n$ goes to infinity, where $\hat\nu_n(\omega) = \bar{Y}_n(\omega)\hat{p}_n(\omega)/[1-\hat{p}_n(\omega)]$.

Then there exists an $N_k(\omega)$ such that for all $n\geq N_k(\omega)$, if we still have $\hat{p}_n(\omega) < 1$, then we must have $\hat\nu_n(\omega) \geq V_{\rm max}$ and $|F_{{\rm NB}(\hat\nu_n(\omega), \hat{p}_n(\omega))}(k) - F_{\lambda}(k)| \leq \frac{\epsilon}{6}$; if  $\hat{p}_n(\omega) = 1$, we still have $\left|F_{\hat\lambda = \bar{Y}_n(\omega)}(k) - F_\lambda(k)\right| \leq \frac{\epsilon}{6}$~. Then for each $n\geq N_{\rm max}(\omega)$ $=$ $\max\{N_0(\omega), N_1(\omega), \ldots, N_M(\omega),$ $N_e(\omega)\}$, $\sup_{k<M}|F_{\hat\mu_n(\omega), \hat{p}_n(\omega)}(k) - F_{\lambda}(k)| \leq \frac{\epsilon}{6}$ and $0 < 1 $ $-$ $F_{\hat\mu_n(\omega), \hat{p}_n(\omega)}(M) \leq 1-F_\lambda(M) + |F_{\hat\mu_n(\omega), \hat{p}_n(\omega)}(M) - F_{\lambda}(M)| \leq \frac{\epsilon}{6} + \frac{\epsilon}{6} = \frac{\epsilon}{3}$, which leads to $D_n(\omega) \leq \epsilon$.

Therefore, $D_n(\omega) $ goes to zero as $n$ goes to infinity.

Similarly to the proof for Corollary~\ref{cor:NB_to_Poisson_nu_p}, we obtain $\sup_y|F_{\hat\mu, \hat{p}}(y) - F_\lambda(y)| \rightarrow 0$ almost surely, as $n \rightarrow \infty$.
\hfill{$\Box$}

\section{More tables for Section~\ref{sec:numerical_studies}}\label{sec:supp_tables}

In this section, we provide more tables for supporting the discussions in Sections~\ref{sec:comparison_MLE} and \ref{sec:real_data}. 

We first extend the comparison analysis by considering 25 parameter pairs $(\nu, p)$ for negative binomial distributions, including $\nu = 0.01, 0.1, 1, 10, 100$ and $p = 0.99, 0.9, 0.5,$ $0.1, 0.01$. In Table~\ref{tab:NB_mean_sd}, we report the mean and variance derived from the 25 parameter pairs $\mathrm{NB}(\nu, p)$ used in our simulation study. 

\begin{table}[ht]
\caption{Mean and variance of $\mathrm{NB}(\nu, p)$}
\label{tab:NB_mean_sd}
\begin{center}
\begin{tabular}{@{}lccccc}
\toprule
& \multicolumn{5}{@{}c}{\textbf{Mean of $\mathrm{NB}(\nu, p)$}} \\\cmidrule{2-6}
$\bm{\nu}\backslash \bm{p}$ & \textbf{0.99} & \textbf{0.90} & \textbf{0.50} & \textbf{0.10} & \textbf{0.01} \\
\midrule
\textbf{0.01} & 0.0001 & 0.0011 & 0.01   & 0.09   & 0.99 \\
\textbf{0.1}  & 0.0010 & 0.0111 & 0.1   & 0.9   & 9.9 \\
\textbf{1}    & 0.0101 & 0.1111 & 1   & 9   & 99 \\
\textbf{10}   & 0.1010 & 1.1111 & 10  & 90  & 990 \\
\textbf{100}  & 1.0101 & 11.1111 & 100 & 900 & 9900 \\
\bottomrule

\toprule
& \multicolumn{5}{@{}c}{\textbf{Variance of $\mathrm{NB}(\nu, p)$}} \\\cmidrule{2-6}
$\bm{\nu}\backslash \bm{p}$ & \textbf{0.99} & \textbf{0.90} & \textbf{0.50} & \textbf{0.10} & \textbf{0.01} \\
\midrule
\textbf{0.01} & 0.0001 & 0.0012 & 0.02   & 0.9   & 99 \\
\textbf{0.1}  & 0.0010 & 0.0123 & 0.2   & 9   & 990 \\
\textbf{1}    & 0.0102 & 0.1235 & 2   & 90  & 9900 \\
\textbf{10}   & 0.1020 & 1.2346 & 20  & 900 & 99000 \\
\textbf{100}  & 1.0203 & 12.3457 & 200 & 9000 & 990000 \\
\bottomrule
\end{tabular}
\end{center}
\end{table}

For each pair $(\nu, p)$, we generate $B=100$ random samples of size $n=100$ and $1,000$, respectively. Table~\ref{tab:Error_rate_supplementary_100} and Table~\ref{tab:Error_rate_supplementary_1000} report the failure rate out of 100 samples of size $n=100$ and $1,000$, respectively, for which the corresponding algorithm fails to return an MLE for NB$(\nu,p)$ distributions. Table~\ref{tab:computational_time_supplementary_100} and Table~\ref{tab:computational_time_supplementary_1000} report the average computational time in seconds over 100 random samples of size $n=100$ and $1,000$, respectively. Table~\ref{tab:ratio_maximum_likelihood_supplementary_100} and Table~\ref{tab:ratio_maximum_likelihood_supplementary_1000} present the average ratio of maximum likelihood over 100 random samples of size $n=100$ and $1,000$, respectively.

\begin{table}[ht]
\caption{\small{Failure rate out of 100 samples for which algorithms fail to return an MLE for NB$(\nu, p)$, with sample size $n=100$}}\label{tab:Error_rate_supplementary_100}
\begin{tabular*}{\textwidth}{@{\extracolsep\fill}lcccccccccc}
\toprule%
& \multicolumn{5}{@{}c@{}}{APMA} & \multicolumn{5}{@{}c@{}}{AZIAD} \\\cmidrule{2-6}\cmidrule{7-11}%
$\bm{\nu \ \backslash \ p}$ & \textbf{0.99} & \textbf{0.9} & \textbf{0.5} & \textbf{0.1} & \textbf{0.01} & \textbf{0.99} & \textbf{0.9} & \textbf{0.5} & \textbf{0.1} & \textbf{0.01} \\
\midrule
\textbf{0.01} & 0   &  0   &  0   &  0   &  0   &  0   &  0   &  0   &  0   &  0 \\ 
\textbf{0.1} & 0   &  0   &  0   &  0   &  0   &  0   &  0   &  0   &  0   &  0 \\ 
\textbf{1} & 0   &  0   &  0   &  0   &  0   &  0   &  0   &  0   &  0   &  0 \\ 
\textbf{10} & 0   &  0   &  0   &  0   &  0   &  0   &  0   &  0   &  0   &  0 \\ 
\textbf{100} & 0   &  0   &  0   &  0   &  0   &  0   &  0   &  0   &  0   &  0 \\
\bottomrule

\toprule%
& \multicolumn{5}{@{}c@{}}{MASS} & \multicolumn{5}{@{}c@{}}{Nelder-Mead} \\\cmidrule{2-6}\cmidrule{7-11}%
$\bm{\nu \ \backslash \ p}$ & \textbf{0.99} & \textbf{0.9} & \textbf{0.5} & \textbf{0.1} & \textbf{0.01} & \textbf{0.99} & \textbf{0.9} & \textbf{0.5} & \textbf{0.1} & \textbf{0.01} \\
\midrule
\textbf{0.01}   &  0.99   &  0.86   &  0.46   &  0.10   &  0.02   &  0.99   &  0.86   &  0.46   &  0.09   &  0.02 \\ 
\textbf{0.1}   &  0.95   &  0.34   &  0   &  0   &  0   &  0.95   &  0.34   &  0   &  0   &  0 \\ 
\textbf{1}   &  0.35   &  0   &  0   &  0   &  0   &  0.35   &  0   &  0   &  0   &  0 \\ 
\textbf{10}   &  0   &  0.23   &  0   &  0   &  0   &  0   &  0   &  0   &  0   &  0 \\ 
\textbf{100}   &  0.56   &  0.20   &  0   &  0   &  0   &  0   &  0   &  0   &  0   &  0 \\
\bottomrule

\toprule%
& \multicolumn{5}{@{}c@{}}{BFGS} & \multicolumn{5}{@{}c@{}}{CG} \\\cmidrule{2-6}\cmidrule{7-11}%
$\bm{\nu \ \backslash \ p}$ & \textbf{0.99} & \textbf{0.9} & \textbf{0.5} & \textbf{0.1} & \textbf{0.01} & \textbf{0.99} & \textbf{0.9} & \textbf{0.5} & \textbf{0.1} & \textbf{0.01} \\
\midrule
\textbf{0.01}   &  0.99   &  0.86   &  0.51   &  0.28   &  0.31   &  0.99   &  0.86   &  0.60   &  0.83   &  0.45 \\ 
\textbf{0.1}   &  0.95   &  0.34   &  0.36   &  0.97   &  0   &  0.95   &  0.34   &  0.90   &  0.95   &  0 \\ 
\textbf{1}   &  0.35   &  0.23   &  1   &  0.06   &  0   &  0.35   &  0.81   &  1   &  0.06   &  0 \\ 
\textbf{10}   &  0.15   &  0.37   &  0   &  0   &  0   &  0.82   &  1   &  1   &  1   &  1 \\ 
\textbf{100}   &  0.17   &  0   &  0   &  0   &  0   &  1   &  1   &  1   &  1   &  1 \\ 
\bottomrule

\toprule%
& \multicolumn{5}{@{}c@{}}{L-BFGS-B} & \multicolumn{5}{@{}c@{}}{SANN} \\\cmidrule{2-6}\cmidrule{7-11}%
$\bm{\nu \ \backslash \ p}$ & \textbf{0.99} & \textbf{0.9} & \textbf{0.5} & \textbf{0.1} & \textbf{0.01} & \textbf{0.99} & \textbf{0.9} & \textbf{0.5} & \textbf{0.1} & \textbf{0.01} \\
\midrule
\textbf{0.01}   &  0.99   &  0.86   &  0.51   &  0.76   &  0.96   &  0.99   &  0.86   &  0.49   &  0.11   &  0.02 \\ 
\textbf{0.1}   &  0.95   &  0.34   &  0.72   &  1   &  1   &  0.95   &  0.34   &  0   &  0   &  0 \\ 
\textbf{1}   &  0.35   &  0.35   &  0.99   &  1   &  1   &  0.35   &  0   &  0   &  0   &  0 \\ 
\textbf{10}   &  0.24   &  0.55   &  0   &  0   &  0   &  0   &  0   &  0   &  0   &  0 \\ 
\textbf{100}   &  0.57   &  0.16   &  0   &  0   &  0   &  0   &  0   &  0   &  0   &  0 \\ 
\bottomrule
\end{tabular*}
\end{table}

\begin{table}[ht]
\caption{\small{Failure rate out of 100 samples for which algorithms fail to return an MLE for NB$(\nu, p)$, with sample size $n=1,000$}}\label{tab:Error_rate_supplementary_1000}
\begin{tabular*}{\textwidth}{@{\extracolsep\fill}lcccccccccc}
\toprule%
& \multicolumn{5}{@{}c@{}}{APMA} & \multicolumn{5}{@{}c@{}}{AZIAD} \\\cmidrule{2-6}\cmidrule{7-11}%
$\bm{\nu \ \backslash \ p}$ & \textbf{0.99} & \textbf{0.9} & \textbf{0.5} & \textbf{0.1} & \textbf{0.01} & \textbf{0.99} & \textbf{0.9} & \textbf{0.5} & \textbf{0.1} & \textbf{0.01} \\
\midrule
\textbf{0.01} & 0   &  0   &  0   &  0   &  0   &  0   &  0   &  0   &  0   &  0 \\ 
\textbf{0.1} & 0   &  0   &  0   &  0   &  0   &  0   &  0   &  0   &  0   &  0 \\ 
\textbf{1} & 0   &  0   &  0   &  0   &  0   &  0   &  0   &  0   &  0   &  0 \\ 
\textbf{10} & 0   &  0   &  0   &  0   &  0   &  0   &  0   &  0   &  0   &  0 \\ 
\textbf{100} & 0   &  0   &  0   &  0   &  0   &  0   &  0   &  0   &  0   &  0 \\
\bottomrule

\toprule%
& \multicolumn{5}{@{}c@{}}{MASS} & \multicolumn{5}{@{}c@{}}{Nelder-Mead} \\\cmidrule{2-6}\cmidrule{7-11}%
$\bm{\nu \ \backslash \ p}$ & \textbf{0.99} & \textbf{0.9} & \textbf{0.5} & \textbf{0.1} & \textbf{0.01} & \textbf{0.99} & \textbf{0.9} & \textbf{0.5} & \textbf{0.1} & \textbf{0.01} \\
\midrule
\textbf{0.01}   &  0.99   &  0.71   &  0   &  0.01   &  0   &  0.99   &  0.74   &  0.03   &  0   &  0 \\ 
\textbf{0.1}   &  0.76   &  0   &  0   &  0   &  0   &  0.76   &  0   &  0   &  0   &  0 \\ 
\textbf{1}   &  0   &  0   &  0   &  0   &  0   &  0   &  0   &  0   &  0   &  0 \\ 
\textbf{10}   &  0   &  0   &  0   &  0   &  0   &  0   &  0   &  0   &  0   &  0 \\ 
\textbf{100}   &  0.01   &  0.46   &  0   &  0   &  0   &  0   &  0   &  0   &  0   &  0 \\ 
\bottomrule

\toprule%
& \multicolumn{5}{@{}c@{}}{BFGS} & \multicolumn{5}{@{}c@{}}{CG} \\\cmidrule{2-6}\cmidrule{7-11}%
$\bm{\nu \ \backslash \ p}$ & \textbf{0.99} & \textbf{0.9} & \textbf{0.5} & \textbf{0.1} & \textbf{0.01} & \textbf{0.99} & \textbf{0.9} & \textbf{0.5} & \textbf{0.1} & \textbf{0.01} \\
\midrule
\textbf{0.01}   &  0.99   &  0.71   &  0   &  0.31   &  0.42   &  0.99   &  0.71   &  0   &  1   &  0.25 \\ 
\textbf{0.1}   &  0.76   &  0   &  0.62   &  1   &  0   &  0.76   &  0.01   &  1   &  1   &  0 \\ 
\textbf{1}   &  0   &  0.10   &  1   &  0   &  0   &  0   &  0.87   &  1   &  0   &  0 \\ 
\textbf{10}   &  0   &  0.11   &  0   &  0   &  0   &  0.40   &  0.99   &  1   &  1   &  1 \\ 
\textbf{100}   &  0   &  0   &  0   &  0   &  0   &  1   &  1   &  1   &  1   &  1 \\  
\bottomrule

\toprule%
& \multicolumn{5}{@{}c@{}}{L-BFGS-B} & \multicolumn{5}{@{}c@{}}{SANN} \\\cmidrule{2-6}\cmidrule{7-11}%
$\bm{\nu \ \backslash \ p}$ & \textbf{0.99} & \textbf{0.9} & \textbf{0.5} & \textbf{0.1} & \textbf{0.01} & \textbf{0.99} & \textbf{0.9} & \textbf{0.5} & \textbf{0.1} & \textbf{0.01} \\
\midrule
\textbf{0.01}   &  1   &  0.98   &  0.31   &  0.99   &  1   &  1   &  0.84   &  0.10   &  0.02   &  0 \\ 
\textbf{0.1}   &  1   &  0.49   &  0.99   &  1   &  1   &  0.97   &  0.01   &  0   &  0   &  0 \\ 
\textbf{1}   &  0.60   &  0.14   &  1   &  1   &  1   &  0   &  0   &  0   &  0   &  0 \\ 
\textbf{10}   &  0.01   &  0.19   &  0   &  0   &  0   &  0   &  0   &  0   &  0   &  0 \\ 
\textbf{100}   &  0.40   &  0.03   &  0   &  0   &  0   &  0   &  0   &  0   &  0   &  0 \\ 
\bottomrule
\end{tabular*}
\end{table}

\begin{table}[ht]
\caption{Average computational time (sec) over 100 simulations with $n=100$}\label{tab:computational_time_supplementary_100}
\tiny
\begin{tabular*}{\textwidth}{@{\extracolsep\fill}lcccccccccc}
\toprule%
& \multicolumn{5}{@{}c@{}}{APMA} & \multicolumn{5}{@{}c@{}}{AZIAD} \\\cmidrule{2-6}\cmidrule{7-11}%
$\bm{\nu \ \backslash \ p}$ & \textbf{0.99} & \textbf{0.9} & \textbf{0.5} & \textbf{0.1} & \textbf{0.01} & \textbf{0.99} & \textbf{0.9} & \textbf{0.5} & \textbf{0.1} & \textbf{0.01} \\
\midrule
\textbf{0.01}   &  0.0005   &  0.0003   &  0.0011   &  0.0011   &  0.0014   &  0.0009   &  0.0011   &  0.0061   &  0.0037   &  0.0039 \\ 
\textbf{0.1}   &  0.0005   &  0.0008   &  0.0006   &  0.0011   &  0.0014   &  0.0011   &  0.0078   &  0.0033   &  0.0022   &  0.0023 \\ 
\textbf{1}   &  0.0007   &  0.0006   &  0.0004   &  0.0007   &  0.0004   &  0.0075   &  0.0059   &  0.0011   &  0.0007   &  0.0006 \\ 
\textbf{10}   &  0.0003   &  0.0001   &  0.0003   &  0.0006   &  0.0006   &  0.0045   &  0.0024   &  0.0005   &  0.0006   &  0.0005 \\ 
\textbf{100}   &  0.0002   &  0.0003   &  0.0008   &  0.0007   &  0.0004   &  0.0031   &  0.0024   &  0.0009   &  0.0007   &  0.0006 \\ 
\bottomrule

\toprule%
& \multicolumn{5}{@{}c@{}}{MASS} & \multicolumn{5}{@{}c@{}}{Nelder-Mead} \\\cmidrule{2-6}\cmidrule{7-11}%
$\bm{\nu \ \backslash \ p}$ & \textbf{0.99} & \textbf{0.9} & \textbf{0.5} & \textbf{0.1} & \textbf{0.01} & \textbf{0.99} & \textbf{0.9} & \textbf{0.5} & \textbf{0.1} & \textbf{0.01} \\
\midrule
\textbf{0.01}   &  0.0005   &  0.0005   &  0.0031   &  0.0082   &  0.0111   &  0.0081   &  0.0125   &  0.0202   &  0.0162   &  0.0184 \\ 
\textbf{0.1}   &  0.0009   &  0.0013   &  0.0074   &  0.0073   &  0.0070   &  0.0166   &  0.0152   &  0.0145   &  0.0139   &  0.0142 \\ 
\textbf{1}   &  0.0015   &  0.0049   &  0.0042   &  0.0025   &  0.0026   &  0.0156   &  0.0142   &  0.0051   &  0.0059   &  0.0059 \\ 
\textbf{10}   &  0.0013   &  0.0061   &  0.0017   &  0.0022   &  0.0012   &  0.0084   &  0.0050   &  0.0038   &  0.0068   &  0.0046 \\ 
\textbf{100}   &  0.0078   &  0.0067   &  0.0024   &  0.0029   &  0.0022   &  0.0055   &  0.0057   &  0.0073   &  0.0042   &  0.0046 \\ 
\bottomrule

\toprule%
& \multicolumn{5}{@{}c@{}}{BFGS} & \multicolumn{5}{@{}c@{}}{CG} \\\cmidrule{2-6}\cmidrule{7-11}%
$\bm{\nu \ \backslash \ p}$ & \textbf{0.99} & \textbf{0.9} & \textbf{0.5} & \textbf{0.1} & \textbf{0.01} & \textbf{0.99} & \textbf{0.9} & \textbf{0.5} & \textbf{0.1} & \textbf{0.01} \\
\midrule
\textbf{0.01}   &  0.0025   &  0.0034   &  0.0089   &  0.0234   &  0.0344   &  0.0024   &  0.0061   &  0.0192   &  0.0369   &  0.0447 \\ 
\textbf{0.1}   &  0.0061   &  0.0056   &  0.0256   &  0.0391   &  0.0185   &  0.0058   &  0.0148   &  0.0453   &  0.0478   &  0.0334 \\ 
\textbf{1}   &  0.0065   &  0.0158   &  0.0192   &  0.0182   &  0.0082   &  0.0141   &  0.0428   &  0.0215   &  0.0202   &  0.0098 \\ 
\textbf{10}   &  0.0061   &  0.0104   &  0.0039   &  0.0055   &  0.0033   &  0.0213   &  0.0230   &  0.0258   &  0.0265   &  0.0242 \\ 
\textbf{100}   &  0.0059   &  0.0041   &  0.0078   &  0.0043   &  0.0047   &  0.0239   &  0.0229   &  0.0267   &  0.0257   &  0.0313 \\ 
\bottomrule

\toprule%
& \multicolumn{5}{@{}c@{}}{L-BFGS-B} & \multicolumn{5}{@{}c@{}}{SANN} \\\cmidrule{2-6}\cmidrule{7-11}%
$\bm{\nu \ \backslash \ p}$ & \textbf{0.99} & \textbf{0.9} & \textbf{0.5} & \textbf{0.1} & \textbf{0.01} & \textbf{0.99} & \textbf{0.9} & \textbf{0.5} & \textbf{0.1} & \textbf{0.01} \\
\midrule
\textbf{0.01}   &  0.0023   &  0.0092   &  0.0115   &  0.0075   &  0.0067   &  0.3050   &  0.5739   &  0.6874   &  0.7933   &  0.7272 \\ 
\textbf{0.1}   &  0.0057   &  0.0095   &  0.0136   &  0.0075   &  0.0067   &  0.6778   &  0.7015   &  0.7122   &  0.7438   &  0.8661 \\ 
\textbf{1}   &  0.0094   &  0.0084   &  0.0034   &  0.0037   &  0.0038   &  0.6739   &  0.4689   &  0.3559   &  0.4134   &  0.4211 \\ 
\textbf{10}   &  0.0046   &  0.0042   &  0.0043   &  0.0039   &  0.0036   &  0.2910   &  0.3650   &  0.4043   &  0.3748   &  0.3831 \\ 
\textbf{100}   &  0.0064   &  0.0071   &  0.0064   &  0.0053   &  0.0056   &  0.3992   &  0.4124   &  0.3823   &  0.3990   &  0.3802 \\  
\bottomrule
\end{tabular*}
\end{table}

\begin{table}[ht]
\caption{Average computational time (sec) over 100 simulations with $n=1,000$}\label{tab:computational_time_supplementary_1000}
\tiny
\begin{tabular*}{\textwidth}{@{\extracolsep\fill}lcccccccccc}
\toprule%
& \multicolumn{5}{@{}c@{}}{APMA} & \multicolumn{5}{@{}c@{}}{AZIAD} \\\cmidrule{2-6}\cmidrule{7-11}%
$\bm{\nu \ \backslash \ p}$ & \textbf{0.99} & \textbf{0.9} & \textbf{0.5} & \textbf{0.1} & \textbf{0.01} & \textbf{0.99} & \textbf{0.9} & \textbf{0.5} & \textbf{0.1} & \textbf{0.01} \\
\midrule
\textbf{0.01}   &  0.0009   &  0.0008   &  0.0009   &  0.0008   &  0.0012   &  0.0080   &  0.0294   &  0.0183   &  0.0104   &  0.0135 \\ 
\textbf{0.1}   &  0.0007   &  0.0010   &  0.0009   &  0.0010   &  0.0050   &  0.0280   &  0.0355   &  0.0083   &  0.0106   &  0.0190 \\ 
\textbf{1}   &  0.0011   &  0.0006   &  0.0008   &  0.0019   &  0.0052   &  0.0426   &  0.0113   &  0.0089   &  0.0067   &  0.0053 \\ 
\textbf{10}   &  0.0007   &  0.0005   &  0.0006   &  0.0036   &  0.0069   &  0.0282   &  0.0125   &  0.0050   &  0.0040   &  0.0123 \\ 
\textbf{100}   &  0.0012   &  0.0015   &  0.0027   &  0.0031   &  0.0032   &  0.0564   &  0.0349   &  0.0136   &  0.0053   &  0.0070 \\ 
\bottomrule

\toprule%
& \multicolumn{5}{@{}c@{}}{MASS} & \multicolumn{5}{@{}c@{}}{Nelder-Mead} \\\cmidrule{2-6}\cmidrule{7-11}%
$\bm{\nu \ \backslash \ p}$ & \textbf{0.99} & \textbf{0.9} & \textbf{0.5} & \textbf{0.1} & \textbf{0.01} & \textbf{0.99} & \textbf{0.9} & \textbf{0.5} & \textbf{0.1} & \textbf{0.01} \\
\midrule
\textbf{0.01}   &  0.0009   &  0.0017   &  0.0033   &  0.0247   &  0.0201   &  0.0257   &  0.0341   &  0.0248   &  0.0236   &  0.0205 \\ 
\textbf{0.1}   &  0.0016   &  0.0042   &  0.0216   &  0.0199   &  0.0484   &  0.0287   &  0.0205   &  0.0195   &  0.0486   &  0.0553 \\ 
\textbf{1}   &  0.0038   &  0.0051   &  0.0235   &  0.0252   &  0.0241   &  0.0181   &  0.0155   &  0.0268   &  0.0370   &  0.0252 \\ 
\textbf{10}   &  0.0033   &  0.0110   &  0.0116   &  0.0094   &  0.0208   &  0.0155   &  0.0139   &  0.0144   &  0.0160   &  0.0411 \\ 
\textbf{100}   &  0.0349   &  0.1868   &  0.0394   &  0.0114   &  0.0114   &  0.0519   &  0.0609   &  0.0444   &  0.0144   &  0.0181 \\  
\bottomrule

\toprule%
& \multicolumn{5}{@{}c@{}}{BFGS} & \multicolumn{5}{@{}c@{}}{CG} \\\cmidrule{2-6}\cmidrule{7-11}%
$\bm{\nu \ \backslash \ p}$ & \textbf{0.99} & \textbf{0.9} & \textbf{0.5} & \textbf{0.1} & \textbf{0.01} & \textbf{0.99} & \textbf{0.9} & \textbf{0.5} & \textbf{0.1} & \textbf{0.01} \\
\midrule
\textbf{0.01}   &  0.0056   &  0.0078   &  0.0058   &  0.0441   &  0.0547   &  0.0039   &  0.0125   &  0.0261   &  0.0698   &  0.0689 \\ 
\textbf{0.1}   &  0.0075   &  0.0062   &  0.0522   &  0.1370   &  0.0731   &  0.0099   &  0.0230   &  0.0809   &  0.2087   &  0.1516 \\ 
\textbf{1}   &  0.0059   &  0.0102   &  0.0871   &  0.1094   &  0.0407   &  0.0249   &  0.0739   &  0.1028   &  0.1189   &  0.0545 \\ 
\textbf{10}   &  0.0050   &  0.0170   &  0.0115   &  0.0092   &  0.0241   &  0.0451   &  0.1052   &  0.1525   &  0.1353   &  0.3398 \\ 
\textbf{100}   &  0.0242   &  0.0491   &  0.0417   &  0.0160   &  0.0167   &  0.2782   &  0.3469   &  0.3312   &  0.1232   &  0.1241 \\  
\bottomrule

\toprule%
& \multicolumn{5}{@{}c@{}}{L-BFGS-B} & \multicolumn{5}{@{}c@{}}{SANN} \\\cmidrule{2-6}\cmidrule{7-11}%
$\bm{\nu \ \backslash \ p}$ & \textbf{0.99} & \textbf{0.9} & \textbf{0.5} & \textbf{0.1} & \textbf{0.01} & \textbf{0.99} & \textbf{0.9} & \textbf{0.5} & \textbf{0.1} & \textbf{0.01} \\
\midrule
\textbf{0.01}   &  0.0050   &  0.0153   &  0.0128   &  0.0074   &  0.0059   &  1.2529   &  1.6548   &  1.0786   &  1.2961   &  1.0961 \\ 
\textbf{0.1}   &  0.0139   &  0.0195   &  0.0103   &  0.0113   &  0.0114   &  1.1444   &  1.1223   &  1.2078   &  3.4429   &  4.3661 \\ 
\textbf{1}   &  0.0156   &  0.0126   &  0.0079   &  0.0084   &  0.0094   &  1.068   &  1.2124   &  1.9263   &  2.5329   &  2.7738 \\ 
\textbf{10}   &  0.0082   &  0.0150   &  0.0139   &  0.0108   &  0.0266   &  1.2004   &  2.4412   &  2.7708   &  4.7207   &  5.8589 \\ 
\textbf{100}   &  0.0335   &  0.0650   &  0.0556   &  0.0191   &  0.0167   &  5.4629   &  6.2008   &  2.4625   &  2.3386   &  2.1775 \\  
\bottomrule
\end{tabular*}
\end{table}

\begin{table}[ht]
\caption{Average ratio of maximum likelihood over 100 simulations with $n=100$}\label{tab:ratio_maximum_likelihood_supplementary_100}
\begin{tabular*}{\textwidth}{@{\extracolsep\fill}lcccccccccc}
\toprule%
& \multicolumn{5}{@{}c@{}}{APMA} & \multicolumn{5}{@{}c@{}}{AZIAD} \\\cmidrule{2-6}\cmidrule{7-11}%
$\bm{\nu \ \backslash \ p}$ & \textbf{0.99} & \textbf{0.9} & \textbf{0.5} & \textbf{0.1} & \textbf{0.01} & \textbf{0.99} & \textbf{0.9} & \textbf{0.5} & \textbf{0.1} & \textbf{0.01} \\
\midrule
\textbf{0.01}   &  1   &  1   &  1   &  1   &  1   &  1   &  1   &  1   &  1   &  1 \\ 
\textbf{0.1}   &  1   &  1   &  1   &  1   &  1   &  1   &  1   &  1   &  1   &  1 \\ 
\textbf{1}   &  1   &  1   &  1   &  1   &  1   &  1   &  1   &  1   &  1   &  1 \\ 
\textbf{10}   &  1   &  1   &  1   &  1   &  1   &  1   &  1   &  1   &  1   &  1 \\ 
\textbf{100}   &  1   &  1   &  1   &  1   &  1   &  0.99   &  0.99   &  1   &  1   &  1 \\ 
\bottomrule

\toprule%
& \multicolumn{5}{@{}c@{}}{MASS} & \multicolumn{5}{@{}c@{}}{Nelder-Mead} \\\cmidrule{2-6}\cmidrule{7-11}%
$\bm{\nu \ \backslash \ p}$ & \textbf{0.99} & \textbf{0.9} & \textbf{0.5} & \textbf{0.1} & \textbf{0.01} & \textbf{0.99} & \textbf{0.9} & \textbf{0.5} & \textbf{0.1} & \textbf{0.01} \\
\midrule
\textbf{0.01}   &  1   &  1   &  0.73   &  0.85   &  0.96   &  1   &  1   &  1   &  1   &  1 \\ 
\textbf{0.1}   &  1   &  0.91   &  0.92   &  1   &  1   &  1   &  1   &  1   &  1   &  0.99 \\ 
\textbf{1}   &  0.99   &  0.94   &  1   &  1   &  1   &  1   &  1   &  1   &  1   &  1 \\ 
\textbf{10}   &  0.91   &  1   &  1   &  1   &  1   &  1   &  1   &  1   &  1   &  1 \\ 
\textbf{100}   &  1   &  0.97   &  1   &  1   &  1   &  1   &  1   &  1   &  1   &  1 \\
\bottomrule

\toprule%
& \multicolumn{5}{@{}c@{}}{BFGS} & \multicolumn{5}{@{}c@{}}{CG} \\\cmidrule{2-6}\cmidrule{7-11}%
$\bm{\nu \ \backslash \ p}$ & \textbf{0.99} & \textbf{0.9} & \textbf{0.5} & \textbf{0.1} & \textbf{0.01} & \textbf{0.99} & \textbf{0.9} & \textbf{0.5} & \textbf{0.1} & \textbf{0.01} \\
\midrule
\textbf{0.01}   &  1   &  1   &  0.70   &  0.79   &  0.94   &  1   &  1   &  0.86   &  0.63   &  0.97 \\ 
\textbf{0.1}   &  1   &  0.91   &  0.89   &  1   &  1   &  1   &  0.91   &  0.91   &  1   &  1 \\ 
\textbf{1}   &  0.99   &  0.93   &  NA   &  1   &  1   &  0.99   &  0.99   &  NA   &  1   &  1 \\ 
\textbf{10}   &  0.91   &  0.95   &  1   &  1   &  1   &  0.98 &  NA   &  NA   &  NA   &  NA \\ 
\textbf{100}   &  0.91   &  0.86   &  0.93   &  1   &  1   &  NA   &  NA   &  NA   &  NA   &  NA \\ 
\bottomrule

\toprule%
& \multicolumn{5}{@{}c@{}}{L-BFGS-B} & \multicolumn{5}{@{}c@{}}{SANN} \\\cmidrule{2-6}\cmidrule{7-11}%
$\bm{\nu \ \backslash \ p}$ & \textbf{0.99} & \textbf{0.9} & \textbf{0.5} & \textbf{0.1} & \textbf{0.01} & \textbf{0.99} & \textbf{0.9} & \textbf{0.5} & \textbf{0.1} & \textbf{0.01} \\
\midrule
\textbf{0.01}   &  1   &  1   &  0.70   &  0.46   &  0.29   &  1   &  1   &  0.72   &  0.31   &  0.46 \\ 
\textbf{0.1}   &  1   &  0.91   &  0.76   &  NA   &  NA   &  1   &  0.94   &  0.46   &  0.36   &  0.83 \\ 
\textbf{1}   &  0.99   &  0.92   &  1   &  NA   &  NA   &  0.99   &  0.85   &  0.39   &  0.47   &  0.83 \\ 
\textbf{10}   &  0.91   &  1   &  1   &  1   &  1   &  0.91   &  0.88   &  0.92   &  0.90   &  0.98 \\ 
\textbf{100}   &  0.99   &  1   &  1   &  1   &  1   &  0.85   &  0.16   &  $<0.01$   &  $<0.01$   &  $<0.01$ \\ 
\bottomrule
\end{tabular*}
Note: ``NA'' indicates that the corresponding algorithm fails to find the MLE in all 100 simulations.
\end{table}

\begin{table}[ht]
\caption{Average ratio of maximum likelihood over 100 simulations with $n=1,000$}\label{tab:ratio_maximum_likelihood_supplementary_1000}
\begin{tabular*}{\textwidth}{@{\extracolsep\fill}lcccccccccc}
\toprule%
& \multicolumn{5}{@{}c@{}}{APMA} & \multicolumn{5}{@{}c@{}}{AZIAD} \\\cmidrule{2-6}\cmidrule{7-11}%
$\bm{\nu \ \backslash \ p}$ & \textbf{0.99} & \textbf{0.9} & \textbf{0.5} & \textbf{0.1} & \textbf{0.01} & \textbf{0.99} & \textbf{0.9} & \textbf{0.5} & \textbf{0.1} & \textbf{0.01} \\
\midrule
\textbf{0.01}   &  1   &  1   &  1   &  1   &  1   &  1   &  1   &  1   &  1   &  1 \\ 
\textbf{0.1}   &  1   &  1   &  1   &  1   &  1   &  1   &  1   &  1   &  1   &  1 \\ 
\textbf{1}   &  1   &  1   &  1   &  1   &  1   &  1   &  1   &  1   &  1   &  1 \\ 
\textbf{10}   &  1   &  1   &  1   &  1   &  1   &  0.99   &  1   &  1   &  1   &  1 \\ 
\textbf{100}   &  1   &  1   &  1   &  1   &  1   &  0.96   &  0.99   &  1   &  1   &  0.99 \\ 
\bottomrule

\toprule%
& \multicolumn{5}{@{}c@{}}{MASS} & \multicolumn{5}{@{}c@{}}{Nelder-Mead} \\\cmidrule{2-6}\cmidrule{7-11}%
$\bm{\nu \ \backslash \ p}$ & \textbf{0.99} & \textbf{0.9} & \textbf{0.5} & \textbf{0.1} & \textbf{0.01} & \textbf{0.99} & \textbf{0.9} & \textbf{0.5} & \textbf{0.1} & \textbf{0.01} \\
\midrule
\textbf{0.01}   &  1   &  0.69   &  0.15   &  0.99   &  0.99   &  1   &  1   &  1   &  1   &  0.91 \\ 
\textbf{0.1}   &  1   &  0.70   &  0.97   &  1   &  1   &  1   &  1   &  1   &  1   &  0.99 \\ 
\textbf{1}   &  0.93   &  0.46   &  1   &  1   &  1   &  1   &  1   &  1   &  1   &  1 \\ 
\textbf{10}   &  0.82   &  0.98   &  1   &  1   &  1   &  1   &  1   &  1   &  1   &  1 \\ 
\textbf{100}   &  0.83   &  0.99   &  1   &  1   &  1   &  1   &  1   &  1   &  1   &  1 \\ 
\bottomrule

\toprule%
& \multicolumn{5}{@{}c@{}}{BFGS} & \multicolumn{5}{@{}c@{}}{CG} \\\cmidrule{2-6}\cmidrule{7-11}%
$\bm{\nu \ \backslash \ p}$ & \textbf{0.99} & \textbf{0.9} & \textbf{0.5} & \textbf{0.1} & \textbf{0.01} & \textbf{0.99} & \textbf{0.9} & \textbf{0.5} & \textbf{0.1} & \textbf{0.01} \\
\midrule
\textbf{0.01}   &  1   &  0.69   &  0.15   &  0.97   &  1   &  1   &  0.69   &  0.15   &  NA   &  1 \\
\textbf{0.1}   &  1   &  0.70   &  0.92   &  NA   &  1   &  1   &  0.70   &  NA   &  NA   &  1 \\ 
\textbf{1}   &  0.93   &  0.39   &  NA   &  1   &  1   &  0.93   &  0.93   &  NA   &  1   &  1 \\ 
\textbf{10}   &  0.82   &  0.95   &  1   &  1   &  1   &  0.94   &  1   &  NA   &  NA   &  NA \\ 
\textbf{100}   &  0.68   &  0.52   &  0.79   &  1   &  1   &  NA   &  NA   &  NA   &  NA   &  NA \\ 
\bottomrule

\toprule%
& \multicolumn{5}{@{}c@{}}{L-BFGS-B} & \multicolumn{5}{@{}c@{}}{SANN} \\\cmidrule{2-6}\cmidrule{7-11}%
$\bm{\nu \ \backslash \ p}$ & \textbf{0.99} & \textbf{0.9} & \textbf{0.5} & \textbf{0.1} & \textbf{0.01} & \textbf{0.99} & \textbf{0.9} & \textbf{0.5} & \textbf{0.1} & \textbf{0.01} \\
\midrule
\textbf{0.01}   &  NA   &  0.02   &  0.03   &  $<0.01$   &  NA   &  NA   &  0.57   &  0.16   &  0.04   &  0.19 \\ 
\textbf{0.1}   &  NA   &  0.49   &  1   &  NA   &  NA   &  1   &  0.72   &  0.12   &  0.18   &  0.64 \\ 
\textbf{1}   &  0.91   &  0.40   &  NA   &  NA   &  NA   &  0.94   &  0.61   &  0.28   &  0.27   &  0.57 \\ 
\textbf{10}   &  0.82   &  0.99   &  1   &  1   &  1   &  0.89   &  0.89   &  0.87   &  0.89   &  0.97 \\ 
\textbf{100}   &  0.83   &  1   &  1   &  1   &  1   &  0.53   &  $<0.01$   &  $<0.01$   &  $<0.01$   &  $<0.01$ \\  
\bottomrule
\end{tabular*}
Note: ``NA'' indicates that the corresponding algorithm fails to find the MLE in all 100 simulations.
\end{table}

In this section, we also provide two supplementary tables to further support the empirical findings discussed in Section~\ref{sec:real_data}. Table~\ref{tab:estimation_prussian} summarizes the results for the \texttt{prussian} dataset, and Table~\ref{tab:estimation_dataCar} reports the results for the \texttt{dataCar} dataset. For each dataset, we apply eight competing algorithms to find the MLE of NB$(\nu, p)$. When an MLE is obtained by an algorithm, we report the estimates $\hat{\nu}$ and $\hat{p}$, the corresponding log-likelihood value, and the computational time in seconds. If the corresponding algorithm fails to return an MLE, we mark it by ``--''. The two supplementary tables provide detailed evidence supporting our conclusions in Section~\ref{sec:real_data} and highlight the stability and advantages of APMA under practical data settings.

\begin{table}[ht]
\caption{Finding MLE of NB$(\nu, p)$ for \texttt{prussian} dataset}
\label{tab:estimation_prussian}
\footnotesize
\begin{tabular}{@{}lrrrrrrrr@{}}
\toprule
$\bm{\nu}$&\textbf{APMA}&\textbf{AZIAD}&\textbf{MASS}&\textbf{Nelder-Mead}&\textbf{BFGS}&\textbf{CG}&\textbf{L-BFGS-B}&\textbf{SANN}\\
\midrule
$\hat{\nu}$ & 7.6072 & 7.6072 & 7.7416 & 7.6191 & 9.9999 & -- & 7.6078 & 10.0000 \\
$\hat{p}$   & 0.9157 & 0.9157 & 0.9171 & 0.9159 & 0.9346 & -- & 0.9157 & 0.9346 \\
Est. Mean & 0.7000 & 0.7000 & 0.6999 & 0.7001 & 0.7000 & -- & 0.7000 & 0.7000 \\
Est. Variance & 0.7644 & 0.7644 & 0.7631 & 0.7644 & 0.7490 & -- & 0.7644 & 0.7490 \\
Log-Likelihood & -313.65 & -313.65 & -313.65 & -313.65 & -313.68 & -- & -313.65 & -313.68 \\
Time Cost (sec) & $<0.01$ & 0.02   & $<0.01$   & 0.03   & $<0.01$  & -- & $<0.01$   & 0.64 \\
\bottomrule
\end{tabular}
\normalsize
Note: ``--'' indicates that the corresponding algorithm fails to find the MLE.
\end{table}

\begin{table}[ht]
\caption{Finding MLE of NB$(\nu, p)$ for \texttt{dataCar} dataset}
\label{tab:estimation_dataCar}
\footnotesize
\begin{tabular}{@{}lrrrrrrrr@{}}
\toprule
$\bm{\nu}$&\textbf{APMA}&\textbf{AZIAD}&\textbf{MASS}&\textbf{Nelder-Mead}&\textbf{BFGS}&\textbf{CG}&\textbf{L-BFGS-B}&\textbf{SANN}\\
\midrule
$\hat{\nu}$ & 1.1568 & 1.1568 & 10.0000 & 1.1562 & 10.0000 & 10.0000 & 10.0000 & 0.9739 \\
$\hat{p}$   & 0.9408 & 0.9408 & 0.9928 & 0.9408 & 0.9928 & 0.9928 & 0.9928 & 0.9276 \\
Log-Likelihood & -18050 & -18050 & -18088 & -18050 & -18088 & -18088 & -18088 & -18055 \\
Time Cost (sec) & 0.01 & 0.74 & 0.20 & 0.78 & 0.16 & 0.83 & 0.34 & 85.59 \\
\bottomrule
\end{tabular}
\normalsize
\end{table}

\section{R code for Algorithm~\ref{algo:MLE_nu_p} and Algorithm~\ref{algo:MLE_mu_p}}\label{sec:r_code}

In this section, we provide R code for implementing Algorithm~\ref{algo:MLE_nu_p} (see function {\tt fmle.nb.new1}) and Algorithm~\ref{algo:MLE_mu_p} (see function {\tt fmle.nb.new2}), as well as applications to three examples.

\medskip

\footnotesize
\begin{verbatim}
#### APMA: Algorithm 1 and Algorithm 2 ####
#### Algorithm 1 to find MLE of NB(nu, p) ####
fmle.nb.new1 <- function(x, nu.max=1e4, epsilon=1e-3, delta=0.1) {
  #### step0
  ny=length(x)
  y=y0=table(x)
  if(sum(x==0)>0) y=y[-1] # remove y=0
  Iset0 = as.numeric(names(y0)) # I, collection of distinct values in y0
  fyvec0 = as.vector(y0) #f_y, frequency of y in I
  Iset = as.numeric(names(y)) # I\{0}, collection of distinct values in y
  fyvec = as.vector(y) #f_y, frequency of y in I\{0}
  
  #### step1. calculate sample mean and sample variance
  ymean=mean(x)
  ysamplevar=var(x)
  
  #### step2. calculate the initial value nu0
  ## define the nu.max, epsilon, and delta
  ## setup if else statement
  if (ymean==0){
    nu.hat=1
    p.hat=1
    loglikelihood=0
  } else{
    nu0=min(nu.max, (ymean^2/max(epsilon, (ysamplevar-ymean))))
    nu=nu0
    
    #### step3 and step4. get the mle of nu, "nu.hat"
    # distinguish few-distinct-values scenario and many-distinct-values scenario
    threshold = delta 
    
    # few-distinct-values scenario
    if (length(unique(x))/ny < threshold) {
      h.self <- function(nu) {
        sum(fyvec0*lgamma(nu+Iset0)) - sum(fyvec0*lgamma(Iset0+1)) - ny*lgamma(nu) - 
          ny*(nu+ymean)*log(nu+ymean) + ny*nu*log(nu) + ny*ymean*log(ymean);
      }
      neg.h.self <- function(nu) {
        (-1)*h.self(nu);
      }
      
      ## define the gnu function "gnu.self"
      gnu.self <- function(nu) {
        # need ny, Iset, fyvec, ymean
        fy.temp <- function(y) {sum(1/(nu + (0:(y-1))));};
        fycoef = fyvec; # initialize coefficients of f_y
        for(i in 1:length(Iset)) fycoef[i]=fy.temp(Iset[i]);
        sum(fyvec*fycoef)/ny - log(1+ymean/nu);
      }
      #
      gr.neg.h <- function(nu) {
        (-1)*ny*gnu.self(nu);
      }
      #
      nuhat.ans=optim(par=nu, fn=neg.h.self, gr=gr.neg.h, method = "L-BFGS-B",
                      lower=epsilon, upper=nu.max)
    } else{
      # many-distinct-values scenario
      h.self <- function(nu) {
        sum(dnbinom(x, size=nu, prob=nu/(nu+ymean),log=T));
      }
      neg.h.self <- function(nu) {
        (-1)*h.self(nu);
      }
      gr.neg.h <- function(nu) {
        -sum(digamma(nu+x))+ny*digamma(nu)+ny*log(1+ymean/nu);
      }
      nuhat.ans=optim(par=nu, fn=neg.h.self, gr=gr.neg.h, method = "L-BFGS-B",
                      lower=epsilon, upper=nu.max)      
    }
    nu.hat=nuhat.ans$par # mle of nu
    
    #### step5. considering nu.max
    if(h.self(nu.max) > h.self(nu.hat)) {
      nu.hat=nu.max
    }
    
    #### step6. get the mle of p "p.hat"
    #### add likelihood of this algorithm
    if (nu.hat == nu.max) {
      warning("nu.hat has reached the upper bound nu.max. You may consider 
      increasing nu.max and rerunning the procedure if necessary.")
    }
    p.hat=nu.hat/(nu.hat+ymean)
    loglikelihood=h.self(nu.hat)
  }
  #### step7. report final answers
  list(nu.hat=nu.hat, p.hat=p.hat, nu.max=nu.max, loglikelihood=loglikelihood);
}

#### Algorithm 2 to find MLE of NB(mu, p) ####
fmle.nb.new2 <- function(x, nu.max=1e4, epsilon=1e-3, delta=0.1) {
  
  #### step0
  ny=length(x)
  y=y0=table(x)
  if(sum(x==0)>0) y=y[-1] # remove y=0
  Iset0 = as.numeric(names(y0)) # I, collection of distinct values in y0
  fyvec0 = as.vector(y0) #f_y, frequency of y in I
  Iset = as.numeric(names(y)) # I\{0}, collection of distinct values in y
  fyvec = as.vector(y) #f_y, frequency of y in I\{0}
  
  #### step1. calculate sample mean and sample variance
  ymean=mean(x)
  ysamplevar=var(x)
  
  #### step2. calculate the initial value nu0
  ## define the nu.max, epsilon, and delta
  ## setup if else statement
  if (ymean==0){
    mu.hat=0
    p.hat=1
    loglikelihood=0
  } else if (ymean >= ysamplevar*(ny-1)/ny) {
    mu.hat=ymean
    p.hat=1
    loglikelihood=sum(dpois(x,lambda=mu.hat,log=T))
  } else{
    nu0=min(nu.max, (ymean^2/max(epsilon, (ysamplevar-ymean))))
    nu=nu0
    
    #### step3 and step4. get the mle of nu, "nu.hat"
    # distinguish few-distinct-values scenario and many-distinct-values scenario
    threshold = delta 
    
    # few-distinct-values scenario
    if (length(unique(x))/ny < threshold) {
      h.self <- function(nu) {
        sum(fyvec0*lgamma(nu+Iset0)) - sum(fyvec0*lgamma(Iset0+1)) - ny*lgamma(nu) - 
          ny*(nu+ymean)*log(nu+ymean) + ny*nu*log(nu) + ny*ymean*log(ymean);
      }
      neg.h.self <- function(nu) {
        (-1)*h.self(nu);
      }
      
      ## define the gnu function "gnu.self"
      gnu.self <- function(nu) {
        # need ny, Iset, fyvec, ymean
        fy.temp <- function(y) {sum(1/(nu + (0:(y-1))));};
        fycoef = fyvec; # initialize coefficients of f_y
        for(i in 1:length(Iset)) fycoef[i]=fy.temp(Iset[i]);
        sum(fyvec*fycoef)/ny - log(1+ymean/nu);
      }
      #
      gr.neg.h <- function(nu) {
        (-1)*ny*gnu.self(nu);
      }
      #
      nuhat.ans=optim(par=nu, fn=neg.h.self, gr=gr.neg.h, method = "L-BFGS-B",
                      lower=epsilon, upper=nu.max)
    } else{
      # many-distinct-values scenario
      h.self <- function(nu) {
        sum(dnbinom(x, size=nu, prob=nu/(nu+ymean),log=T));
      }
      neg.h.self <- function(nu) {
        (-1)*h.self(nu);
      }
      gr.neg.h <- function(nu) {
        -sum(digamma(nu+x))+ny*digamma(nu)+ny*log(1+ymean/nu);
      }
      nuhat.ans=optim(par=nu, fn=neg.h.self, gr=gr.neg.h, method = "L-BFGS-B",
                      lower=epsilon, upper=nu.max)      
    }
    nu.hat=nuhat.ans$par # mle of nu
    
    #### step5. considering nu.max
    if(h.self(nu.max) > h.self(nu.hat)) {
      nu.hat=nu.max
    }
    
    #### step6. get the mle of p "p.hat"
    #### add likelihood of this algorithm
    if (nu.hat == nu.max) {
      warning("nu.hat has reached the upper bound nu.max. You may consider 
      increasing nu.max and rerunning the procedure if necessary.")
    }
    p.hat=nu.hat/(nu.hat+ymean)
    loglikelihood=h.self(nu.hat)
  }
  #### step7. report final answers
  list(mu.hat=ymean, p.hat=p.hat, nu.max=nu.max, loglikelihood=loglikelihood);
}

#### Example 1
library(datasets)
deaths_data <- as.numeric(UKDriverDeaths) 

fmle.nb.new1(deaths_data)
# $nu.hat
# [1] 34.99521
# 
# $p.hat
# [1] 0.02052141
# 
# $nu.max
# [1] 10000
# 
# $loglikelihood
# [1] -1356.043

fmle.nb.new2(deaths_data)
# $mu.hat
# [1] 1670.307
# 
# $p.hat
# [1] 0.02052141
# 
# $nu.max
# [1] 10000
# 
# $loglikelihood
# [1] -1356.043

#### Example 2
nu = 1
p = 0.99
set.seed(100)
data_temp = rnbinom(1000, size=nu, prob=p)
mean(data_temp) # 0.01
var(data_temp)*(1000-1)/1000  # 0.0099

fmle.nb.new1(data_temp)
# $nu.hat
# [1] 10000
# 
# $p.hat
# [1] 0.999999
# 
# $nu.max
# [1] 10000
# 
# $loglikelihood
# [1] -56.05171
# 
# Warning message:
#   In fmle.nb.new1(data_temp) :
#   nu.hat has reached the upper bound nu.max. You may consider 
# increasing nu.max and rerunning the procedure if necessary.

fmle.nb.new2(data_temp)
# $mu.hat
# [1] 0.01
# 
# $p.hat
# [1] 1
# 
# $nu.max
# [1] 10000
# 
# $loglikelihood
# [1] -56.0517

#### Example 3
set.seed(100)
data_temp = rpois(1000,lambda=5)

fmle.nb.new1(data_temp)
# $nu.hat
# [1] 10000
# 
# $p.hat
# [1] 0.9994865
# 
# $nu.max
# [1] 10000
# 
# $loglikelihood
# [1] -2219.935
# 
# Warning message:
#   In fmle.nb.new1(data_temp) :
#   nu.hat has reached the upper bound nu.max. You may consider 
# increasing nu.max and rerunning the procedure if necessary.

fmle.nb.new2(data_temp)
# $mu.hat
# [1] 5.138
# 
# $p.hat
# [1] 1
# 
# $nu.max
# [1] 10000
# 
# $loglikelihood
# [1] -2219.934
\end{verbatim}
\normalsize

\end{document}